\documentclass[final,onefignum,onetabnum,letterpaper,reqno]{siamonline190516}

\usepackage{lipsum}
\usepackage{amsfonts}
\usepackage{epstopdf}
\ifpdf
  \DeclareGraphicsExtensions{.eps,.pdf,.png,.jpg}
\else
  \DeclareGraphicsExtensions{.eps}
\fi

\usepackage{datetime}
\newdateformat{monthyeardate}{%
	\monthname[\THEMONTH] \THEDAY, \THEYEAR}

\usepackage{academicons}
\usepackage{xcolor}
\renewcommand{\orcid}[1]{\href{https://orcid.org/#1}{\textcolor[HTML]{A6CE39}{orcid.org/#1}}}

\usepackage{amsmath}
\allowdisplaybreaks
\usepackage{amssymb}
\usepackage{commath}
\usepackage{mathtools}
\newtheorem{manualtheoreminner}{Theorem}
\newenvironment{manualtheorem}[1]{%
  \renewcommand\themanualtheoreminner{#1}%
  \manualtheoreminner
}{\endmanualtheoreminner}
\newtheorem{manualcorollaryinner}{Corollary}
\newenvironment{manualcorollary}[1]{%
  \renewcommand\themanualcorollaryinner{#1}%
  \manualcorollaryinner
}{\endmanualcorollaryinner}
\usepackage{bbm}
\usepackage{tabularx}
\usepackage{bm}
\usepackage{physics}
\usepackage{array, ltablex, multirow}

\usepackage{color}
\usepackage{graphicx}
\usepackage[small]{caption}
\usepackage{subcaption}

\usepackage{relsize}
\usepackage{adjustbox}
\usepackage{algorithmic}
\usepackage{booktabs}
\usepackage{tikz}
\usepackage{pifont}
\usetikzlibrary{bayesnet} 

\DeclareMathOperator{\Ker}{Ker}

\DeclareMathOperator{\diag}{diag}
\DeclareMathOperator{\sgn}{sgn}
\DeclareMathOperator{\dB}{dB}

\usepackage{enumitem}
\setlist[enumerate]{leftmargin=.5in}
\setlist[itemize]{leftmargin=.5in}

\newtheorem{remark}{Remark}
\newsiamremark{example}{Example}
\newsiamremark{hypothesis}{Hypothesis}
\crefname{hypothesis}{Hypothesis}{Hypotheses}
\newsiamthm{claim}{Claim}

\makeatother

\headers{Complex-Valued Signal Recovery using the Bayesian LASSO}{D.\ Green, J.\ Lindbloom, and A.\ Gelb}

\title{Complex-Valued Signal Recovery using the Bayesian LASSO\footnote{This work is partially supported by NSF DMS grant \#1912685, AFOSR  grant \#F9550-22-1-0411, DOE ASCR grant \#DE-ACO5-000R22725, and DoD ONR MURI  grant \#N00014-20-1-2595.}
}
\author{ 
Dylan Green\thanks{Department of Mathematics, Dartmouth College, Hanover, NH 03755, USA (\email{Dylan.P.Green.GR@Dartmouth.edu}, \orcid{0000-0001-6184-0183} and \email{Jonathan.T.Lindbloom@Dartmouth.edu}, \orcid{0000-0002-1789-2629} and \email{Anne.E.Gelb@Dartmouth.edu}, \orcid{0000-0002-9219-4572})}
\and 
Jonathan Lindbloom\footnotemark[2]
\and 
Anne Gelb\footnotemark[2]
}

\usepackage{amsopn}

\usepackage{comment}

\usepackage[normalem]{ulem}
\usepackage{upgreek}

\usepackage[normalem]{ulem}

\DeclareMathOperator*{\argmin}{arg\,min}

\renewcommand{\d}{\mathrm{d}} 
 
\renewcommand{\vec}{\, \mathrm{vec}}

\newcommand{\C}{\mathbb{C}}
\renewcommand{\Re}{\mathrm{Re}}
\renewcommand{\Im}{\mathrm{Im}}
\input{jlmacros}

\usepackage{lineno}

\ifpdf
\hypersetup{
	pdftitle={2022EmpBayesComplex},
 	pdfauthor={Dylan Green}
}
\fi

\begin{document}

\maketitle


\begin{abstract}
	\normalsize
Recovering complex-valued image recovery from noisy indirect data is important in applications such as ultrasound imaging and synthetic aperture radar. While there are many effective algorithms to recover point estimates of the magnitude, fewer are designed to recover the phase. Quantifying uncertainty in the estimate can also provide valuable   information for real-time decision making. This investigation therefore proposes a new Bayesian inference method that recovers point estimates while also quantifying the uncertainty for  complex-valued signals or images given noisy and indirect observation data.   Our method is motivated by the Bayesian LASSO approach for real-valued sparse signals, and here we demonstrate that the Bayesian LASSO can be effectively adapted to recover complex-valued images whose magnitude is sparse in some (e.g.~the gradient) domain. Numerical examples demonstrate our algorithm's robustness to noise as well as its computational efficiency.
\end{abstract}

\begin{keywords}
  Bayesian inference,  Bayesian LASSO, complex-valued image recovery
\end{keywords}

\begin{AMS}
	62C10, 62F15, 65C20
\end{AMS}


\section{Introduction}
\label{sec:intro}
 Recovering complex-valued images or signals\footnote{We  use the terms  image and signal interchangeably throughout our manuscript.} from noisy and/or under-sampled data is important in coherent imaging applications such as synthetic aperture radar (SAR) \cite{cheney2001mathematical,jakowatz2012spotlight}, ultrasound \cite{ylitalo1994ultrasound}, and digital holography \cite{yaroslavskii1980methods}. The problem is often modeled as
\begin{align}\label{eqn:observation_model}
    \bm y=F\bm z+\bm\varepsilon.
\end{align}
Here  $\bm z = \bm g \odot e^{i{\bm \phi}} \in\mathbb{C}^n$ with $\odot$ indicating componentwise multiplication is the unknown signal of interest decomposed into magnitude $\bm g\in\left(\mathbb{R^+}\right)^n$  and phase $\bm \phi\in[-\pi,\pi)^n$, $\bm y\in\mathbb{C}^m$ represents the observable data, $F\in\mathbb{C}^{m\cross n}$ is the known forward linear operator, and $\bm\varepsilon\in\mathbb{C}^m$ is centered Gaussian noise with covariance $\sigma^2\mathbb{I}_m$ and probability density
\begin{align*}
    f_\mathcal{E}(\bm\varepsilon)=\frac{1}{\pi^m\sigma^{2m}}\exp(-\frac{1}{\sigma^2}\norm{\bm\varepsilon}_2^2).
\end{align*}

Many techniques seek to exclusively recover the magnitude ${\bm g}$, that is, without considering the phase $\bm{\phi}$,  and as such are designed to promote sparsity in $\bm g$ or in some known transformation of $\bm g$ (e.g., its gradient) \cite{lustig2007sparse,sanders2017composite}. However  phase information is also useful for some applications, including  coherent change detection \cite{barber2015generalized,jakowatz2012spotlight,novak2005coherent,wahl2016new}, high-resolution imaging \cite{touzi1999coherence}, and interferometry \cite{gens1996review,jakowatz2012spotlight}. In this regard, some methods have been developed to recover point estimates of complex-valued signals to promote sparsity of $\bm g$ or of some transformation of $\bm g$ to a sparse domain, see e.g.~\cite{cetin2006feature,green2023leveraging,rajagopal2023enhanced,scarnati2018joint,yu2012compressed}. 

More recent approaches have incorporated uncertainty quantification (UQ) into complex-valued signal recovery methods (see e.g.~\cite{churchill2022sampling,churchill2023sub,duan2015pattern}), thereby providing useful information for real-time decision making. These methods are generally designed to incorporate prior sparse knowledge of the magnitude or some linear transform of the magnitude of the signal into the recovery of the posterior. In the case where sparsity is in some linear transform of the magnitude, however, these methods do not \emph{infer} information regarding the phase. Instead the phase is approximated as part of an optimization step, limiting the uncertainty information available.

The method developed in this investigation addresses some of these issues and  brings a more comprehensive approach to recovering complex-valued images. Specifically, we build on the {\em Bayesian LASSO} technique \cite{park2008bayesian}, which was originally designed to recover sparse {\it real}-valued signals. While other inference methods, such as generalized sparse Bayesian learning \cite{glaubitz2023generalized}, recover uncertainty information for the signal itself, the Bayesian LASSO method also provides UQ for hyperparameters that describe the structure and overall sparsity of the image. The idea in Bayesian LASSO is to recast the $\ell_1$-regularized optimization problem (commonly employed in compressed sensing applications \cite{candes2008enhancing,donoho2006compressed}) in a Bayesian framework by treating both the data and the unknown signal as random variables, yielding the linear system
\begin{equation}
\mathcal{Y}=F\mathcal{X}+\mathcal{E}.
\label{eq:linearinverse_real}
\end{equation}
The forward operator $F\in\mathbb{R}^{m\cross n}$ is known, and $\mathcal{X}\in\mathbb{R}^n$, $\mathcal{E},\mathcal{Y}\in\mathbb{R}^m$ are random variables with respective realizations $\bm x$, $\bm\varepsilon$, $\bm y =F\bm x+\bm\varepsilon$. An estimate of the full posterior density function of the real-valued signal can then be recovered using Bayes' theorem, expressed as
\begin{equation}\label{eqn:bayes}
    f_{\mathcal{X}|\mathcal{Y}}(\bm x|\bm y) \propto f_{\mathcal{Y}|\mathcal{X}}(\bm y|\bm x) f_{\mathcal{X}|\upeta}(\bm x|\eta)f_{\upeta}(\eta).
\end{equation}
Here $f_{\mathcal{Y}|\mathcal{X}}(\bm y|\bm x)$ is the likelihood density function determined by $F$ and assumptions on $\mathcal{E}$,  $f_{\mathcal{X}|\upeta}(\bm x|\eta)$ is the prior density function that encodes a priori assumptions about the unknown, and $f_\upeta(\eta)$ is the hyperprior density function on the scale parameter $\eta$ of the prior density. For example, for sparsifying transform matrix  $L \in \mathbb{R}^{k\cross n}$,   $f_{\mathcal{X}|\eta}(\bm x|\eta)$ may be defined as the product of Laplace probability densities yielding
\begin{equation}\label{eq:prior_real}
    f_{\mathcal{X}|\upeta}(\bm x|\eta)\propto\prod_{i=1}^{k}\frac{1}{2\eta}\exp\left(-\frac{1}{\eta}|L_i\bm x|\right),
\end{equation}
where $L_i$ is the $i$th row of $L$ for $i=1,\dots,k$. 

Our new method extends the Bayesian LASSO technique to {\it complex}-valued signals whose magnitude is sparse in some domain by treating $\bm y$, $\bm z$, and $\bm\varepsilon$  in \eqref{eqn:observation_model} as realizations of respective random variables $\mathcal{Y}$, $\mathcal{Z}$, and $\mathcal{E}$, leading to the probabilistic forward model
\begin{align}
\mathcal{Y}=F\mathcal{Z}+\mathcal{E}.
\label{eq:linearinverse_complex}
\end{align}
The corresponding estimate of the full posterior density function of the complex-valued signal can then be recovered according to 
\begin{equation}\label{eqn:bayescomplex}
    f_{\mathcal{Z}|\mathcal{Y}}(\bm z|\bm y) \propto f_{\mathcal{Y}|\mathcal{Z}}(\bm y|\bm z) f_{\mathcal{Z}|\upeta}(\bm z|\eta)f_\upeta(\eta),
\end{equation}
where similarly $f_{\mathcal{Y}|\mathcal{Z}}(\bm y|\bm z)$ is the likelihood density function determined by $F$ and assumptions on $\mathcal{E}$, $f_{\mathcal{Z}|\upeta}(\bm z|\eta)$ is the prior density function encoding a priori assumptions about the unknown, and $f_\upeta(\eta)$ is the hyperprior density function on $\eta$, the scale parameter of the prior density. We call our resulting method the complex-valued Bayesian LASSO (CVBL). As a primary benefit, the CVBL allows us to fully exploit the sparsity of the underlying complex-valued signal in the recovery without sacrificing the phase, all while quantifying the uncertainty regarding the entire complex-valued signal and the hyperparameters that describe its structure and sparsity.

\subsection*{Our contribution}
Given noisy indirect observable data, we  introduce a new Bayesian model that uses \emph{a priori} assumptions regarding the magnitude of the underlying complex-valued image of interest to recover its posterior distribution as well as to quantify the uncertainty for both the magnitude and the phase of the unknown signal. Adapting the Bayesian LASSO approach  allows us to develop efficient sampling techniques to simulate random draws from these resulting posterior distributions.  
\subsection*{Paper organization}
The rest of this paper is organized as follows.   \Cref{sec:preliminaries} details the construction of the likelihood, prior, and hyperprior densities for our method. \Cref{sec:gibbssamplelaplace} discusses the real-valued Bayesian LASSO approach for sparse signal recovery, which we then expand to include recovery of signals that are sparse in some transform domain.  The corresponding {\em complex-valued} Bayesian LASSO is then proposed in \Cref{sec:proposed_methods}.  Numerical experiments in  \Cref{sec:numerics} consider different forward operators $F$ in \eqref{eqn:observation_model} as well as various signal to noise (SNR) values, demonstrating our method's utility and robustness to noise.  We also compare our results to those obtained using the more classical LASSO maximum a posteriori (MAP) estimate approach \cite{tibshirani1996regression}. We provide some concluding remarks and ideas for future work in \Cref{sec:conclusion}.

\section{Bayesian Formulation}\label{sec:preliminaries}
Since \eqref{eqn:observation_model} is easily understood for one-dimensional problems, we develop our method for $\bm z \in \mathbb{C}^n$ and note that higher-dimensional signals can be readily vectorized to fit this form.  Section \ref{sec:numerics} includes both one and two-dimensional examples.
Below we collect the ingredients needed for our new method, including a description of the likelihood density function in Section \ref{subsec:likelihood}, the construction of sparsity-promoting priors in Section \ref{sec:prior}, and a review of commonly used sparsifying transform operators in Section \ref{sec:sparseoperator}.

\subsection{The likelihood}\label{subsec:likelihood}
The likelihood density function $f_{\mathcal{Y}|\mathcal{Z}}(\bm y|\bm z)$ in \eqref{eqn:bayescomplex} is determined from the density function  of the noise present in the system. Following  \cite{park2008bayesian} and adjusting accordingly for the complex-valued signal case, here we assume that $\mathcal{E}$ follows a central complex normal distribution with covariance $\sigma^{2}\mathbb{I}_m$, yielding
\begin{equation*}
f_\mathcal{E}(\bm\varepsilon)=\mathcal{CN}(\bm0,\sigma^{2}\mathbb{I}_m)
    =\frac{1}{\pi^m\sigma^{2m}}\exp(-\frac{1}{\sigma^2}\norm{\bm \varepsilon}_2^2).
\end{equation*}
The law of total probability provides that
\begin{equation}\label{eq:totalprob}
    f_{\mathcal{Y}|\mathcal{Z}}(\bm y|\bm z)=\int_{\mathbb{C}^m} f_{\mathcal{Y}|\mathcal{Z},\mathcal{E}}(\bm y|\bm z,\bm\varepsilon)f_\mathcal{E}(\bm\varepsilon)\mathrm{d}\bm\varepsilon.
\end{equation}
When conditioned on $\bm z=\mathcal{Z}$ and $\bm\varepsilon=\mathcal{E}$, \cref{eq:linearinverse_complex} implies that $\mathcal{Y}$ is entirely determined according to the density function
\begin{equation}\label{eq:y_delta}
    f_{\mathcal{Y}|\mathcal{Z},\mathcal{E}}(\bm y|\bm z,\bm\varepsilon)=\delta(\bm y-F\bm z-\bm\varepsilon).
\end{equation}
Combining \eqref{eq:totalprob} and \cref{eq:y_delta}, the likelihood density function is therefore
\begin{equation}\label{eq:likelihood}
    f_{\mathcal{Y}|\mathcal{Z}}(\bm y|\bm z)=\int_{\mathbb{C}^m} \delta(\bm y-F\bm z-\bm\varepsilon)f_\mathcal{E}(\bm\varepsilon)\mathrm{d}\bm\varepsilon
    =f_\mathcal{E}(\bm y-F\bm z)
    =\frac{1}{\pi^m\sigma^{2m}}\exp(-\frac{1}{\sigma^2}\norm{\bm y-F\bm z}_2^2).
\end{equation}

\begin{remark}
\label{rem:knownvariance}
We note that an improper hyperprior was placed on $\sigma^2$ in \cite{park2008bayesian}.  For simplicity here we assume that $\sigma^2$ is either explicitly known or easily estimated. In general, $\sigma^2$ may be provided \emph{a priori} when the errors present in the forward model are well understood.  In other situations, an estimate of $\sigma^2$ may be acquired when an area of low intensity in $\mathcal{Z}$ is known or when multiple measurement vectors (MMV) are available, see e.g.~\cite{immerkaer1996fast,sanders2020effective,zhang2022empirical}.
\end{remark}

\subsection{The prior density function}\label{sec:prior}
A primary goal for our new method is to ensure that the prior density function used in \cref{eqn:bayescomplex} enables efficient computation for the complex-valued Bayesian LASSO (CVBL).  In this regard we assume that either the magnitude or some linear transform of the magnitude of the unknown signal is sparse. 

\subsubsection{Sparse magnitude}
\label{sec:sparsemagnitude}
Let $\mathcal{Z}$ be decomposed as $\mathcal{Z}=\mathcal{A}+i\mathcal{B}$, where $\mathcal{A}\in\mathbb{R}^n$ and $\mathcal{B}\in\mathbb{R}^n$ are mutually independent with respective realizations $\bm a$ and $\bm b$. When $|\mathcal{Z}|$ is sparse we can simply adapt the approach from \cite{zhang2022empirical} used for real-valued signals and define the conditional prior density function as
\begin{align}
f_{\mathcal{Z}|\upeta}(\bm z|\eta) =
    f_{\mathcal{A},\mathcal{B}|\upeta}(\bm a,\bm b|\eta)=\frac{1}{(2\eta)^{n}}\exp(-\frac1\eta\sum_{j=1}^n\sqrt{a_j^2+b_j^2}),
    \label{eq:prior_signal}
\end{align}
where $\upeta\in\mathbb{R}^+$ is a random variable with realization  $\eta$.

\subsubsection{Sparse transform of magnitude} \label{sec:cond_Gaussian} 
When the magnitude is presumed sparse in a domain other than the imaging domain, we first reformulate the model as 
\begin{align}
	\mathcal{Y}=F\mathcal{Z}+\mathcal{E}=F\left(\mathcal{G}\odot e^{i{\Phi}}\right)+\mathcal{E}.
 \label{eq:splitmagphase}
\end{align}
Here $\mathcal{Z}=\mathcal{G}\odot e^{i\Phi}$, where $\mathcal{G}\in\left(\mathbb{R}^+\right)^{n}$ and $\Phi\in[-\pi,\pi)^{n}$, with $\mathcal{G}$ and $\Phi$  assumed to be mutually independent.  The componentwise decomposition of $\mathcal{Z}$ allows us to rewrite the likelihood and prior density functions in \eqref{eq:likelihood} respectively as
\begin{align*}
	f_{\mathcal{Y}|{\mathcal{G},\Phi}}(\bm y|\bm g,\bm\phi)=\frac{1}{\pi^n\sigma^{2n}}\exp\left(-\frac{1}{\sigma^2}\norm{\bm y- F\left(\bm g\odot e^{i\bm \phi}\right)}^2_2\right), \quad f_{\mathcal{Z}|\upeta}(\bm z|\eta)=f_{\mathcal{G}|\upeta}(\bm g|\eta)f_\Phi(\bm\phi).
\end{align*}

Let $L \in \mathbb{R}^{k \times n}$ be the rank $n$ operator that transforms $\mathcal{G}$ to the sparse domain.\footnote{The rank $n$ requirement is needed for the particular computational implementation  in \Cref{alg:gibbs_laplace} and \Cref{alg:propMethod}.  A rank deficient sparse transform operator is commonly augmented by imposing boundary constraints \cite{kaipio2000statistical}.} Denoting $\bm1_{\mathbb{R}^+}$ as the indicator function for positive real vectors,
the conditional prior probability density for $\mathcal{G}$ is then
\begin{align}
	f_{{\mathcal{G}}|{\upeta}}(\bm g|\eta) \propto \frac{1}{(2\eta)^{k}}\exp\left(-\frac{1}{\eta}\norm{L\bm g}_1\right)\bm1_{\mathbb{R}^+}(\bm g),
\label{eq:prior_edge}
\end{align}
where again $\upeta\in\mathbb{R}^+$ is  a random variable with realization $\eta$.   

As there is often no prior information regarding the phase $\Phi$ in coherent imaging systems,  it is reasonable to impose the uniform prior  as
\begin{equation}
\label{eqn:phaseprior}
f_\Phi(\bm\phi)=\frac{1}{(2\pi)^n}\bm1_{[-\pi,\pi)}(\bm\phi),
\end{equation} 
where $\bm1_{[-\pi,\pi)}$ denotes the indicator function for vectors whose elements are in $[-\pi,\pi)$.

\subsubsection{The hyperprior}\label{subsub:hyperprior}
There are, of course, various options for the hyperprior on $\upeta$ in \eqref{eq:prior_signal} and \eqref{eq:prior_edge}. In some cases it may be appropriate to simply use the delta density function
\begin{align}\label{eq:delta_hyper}
    f_{\upeta}(\eta)=\delta(\eta-\hat{\eta}),
\end{align}
with point estimate $\hat{\eta}$ for $\upeta$. Such an estimate may be available when considering signals whose sparsity is known or well-approximated, see e.g.~\cite{sanders2020effective,zhang2022empirical}. Here we follow the original Bayesian LASSO method \cite{park2008bayesian} and place a gamma hyperprior on $\upeta^{-2}$ in order to maintain conjugacy,  yielding the probability density function\footnote{See discussion surrounding \cref{eq:RVBL_prior} for explanation regarding using $\upeta^{-2}$ in place of $\upeta$ in \eqref{eqn:bayes}.}
\begin{align}\label{eq:hyperprior}
    f_{\upeta^{-2}}(\eta^{-2})=\frac{\delta^r}{\Gamma(r)}(\eta^{-2})^{r-1}\exp(-\delta\eta^{-2}), \quad r > 0, \delta > 0.
\end{align}
Note that choosing \emph{shape} parameter $r\leq1$ ensures that the mode of $f_{\upeta^{-2}}(\eta^{-2})$ is zero while also encouraging sparsity in the solution.  Furthermore, having \emph{rate} parameter $\delta\ll1$ gives $f_{\upeta^{-2}}(\eta^{-2})$ a fatter tail, making the hyperprior relatively uninformative.  To encourage sparsity while still allowing for small values of $\eta^{-2}$, unless otherwise noted our numerical experiments use hyperparameters
\begin{align}\label{eq:r_and_delta}
    r = 1,\quad \delta = 10^{-3},
\end{align}
with no additional tuning.

\subsection{The sparsifying transform operator}
\label{sec:sparseoperator}
This investigation assumes that the underlying signal is either sparse or that its magnitude is piecewise constant, in which case the first order differencing (TV) operator is well suited to suppress variation and noise in smooth regions.  The numerical experiment  in \Cref{subsec:numerics_sarsemagtransform} verifies that this assumption is reasonable even when the true signal has piecewise smooth  (but not constant) magnitude. We note, however, that such an assumption is not an inherent limitation to our new method, which can be straightforwardly adapted to other sparsifying transform operators, such as higher order TV (HOTV) \cite{archibald2016image} or wavelets \cite{aguilera2013wavelet,loris2007tomographic} as appropriate. 

\begin{remark}\label{rem:flat_post}
    When using $L\in\mathbb{R}^{k\cross n}$ with $k>n$ as the sparsifying transform, for instance when $L=[L_1^T \ L_2^T]^T$ is the 2D gradient operator with $L_1\in\mathbb{R}^{n\cross n}$ computing vertical differences and $L_2\in\mathbb{R}^{n\cross n}$ computing horizontal differences, the posterior magnitude may become overly smoothed. We conjecture that when $k\gg n$, the resulting posterior density is multimodal and has significant mass concentrated where $\eta^{-2}$ is large and the magnitude is relatively flat. Future investigations will explore  methods to overcome this undesirable result.
\end{remark}

Before introducing our approach for complex-valued signal recovery in \Cref{sec:proposed_methods},  \Cref{sec:gibbssamplelaplace}  first discusses its real-valued analog.  We use what we will call the {\em real-}valued Bayesian LASSO (RVBL) technique \cite{park2008bayesian}, which was initially developed to sample a posterior consisting of a Gaussian likelihood and Laplace prior on the signal itself, that is, where $L$ is the identity matrix. The technique was extended to include the anisotropic TV operator in \cite{jarvenpaa2014bayesian}. We now adapt these ideas to accommodate \emph{any} sparsifying rank $n$ linear operator $L$ for  real-valued signal recovery, which to our knowledge has not been previously done.

\section{Bayesian LASSO for a real-valued signal (RVBL)}\label{sec:gibbssamplelaplace} The RVBL method is a blocked Gibbs sampling technique. Critical to the approach is the equivalent representation of each element of the product in \eqref{eq:prior_real} as a scale mixture of Gaussians with an exponential mixing density written for each component $L_i\bm x$, $i = 1,\dots,k$. In particular,
\begin{align}
    \frac{1}{2\eta}\exp\left(-\frac{1}{\eta}|L_i\bm x|\right) =\int_0^\infty\frac{1}{\sqrt{2\pi s_i}}\exp\left(-\frac{(L_i\bm x)^2}{2s_i}\right)\frac{1}{2\eta^2}\exp\left(-\frac{s_i}{2\eta^2}\right)\mathrm{d}s_i.
    \label{eq:scalemix}
\end{align}
Substituting \eqref{eq:scalemix} into \eqref{eq:prior_real} and denoting  $D(\bm s)=\diag(\bm s)$ yields
\begin{align}\label{eq:simplified_mix}
    f_{\mathcal{X}|\upeta}(\bm x|\eta)
    =\int_{(\mathbb{R}^+)^k}\frac{1}{(2\pi)^{\frac{k}{2}} |D(\bm s)|^{\frac12}}\exp\left(-\frac12(L\bm x)^T[D(\bm s)]^{-1}(L\bm x)\right)\frac{1}{2^{k}\eta^{2 k}}\exp\left(-\frac{\norm{\bm s}_1}{2\eta^2}\right)\mathrm{d}\bm s, 
\end{align}

Now consider the marginal density of $\mathcal{X}$ over the joint distribution of $\mathcal{X}$ and some random variable $\mathcal{T}^2\in\mathbb{R}^n$, which we refer to as the scale mixture parameter, given by 
\begin{align}
    f_{\mathcal{X}|\upeta}(\bm x|\eta) = \int_{(\mathbb{R}^+)^k} f_{\mathcal{X},\mathcal{T}^2|\upeta}(\bm x,\bm\tau^2|\eta)\mathrm{d}\bm\tau^2 = \int_{(\mathbb{R}^+)^k} f_{\mathcal{X}|\mathcal{T}^2,\upeta}(\bm x|\bm\tau^2,\eta)f_{\mathcal{T}^2|\upeta}(\bm\tau^2|\eta)\mathrm{d}\bm\tau^2.
    \label{eqn:joint_to_marginal}
\end{align}
Note that in \cref{eqn:joint_to_marginal} $\bm \tau^2\in\mathbb{R}^k$ is defined componentwise, with
\begin{align}\label{eq:vector_square}
    \bm\tau^2= \begin{bmatrix} \tau_1^2 & \tau_2^2 & \cdots & \tau_k^2\end{bmatrix}^T.
\end{align}
By making the substitution $\bm\tau^2$ for $\bm s$ and directly comparing the integrand terms in \eqref{eq:simplified_mix} to those in \eqref{eqn:joint_to_marginal},  we obtain the Gaussian density function
\begin{align}\label{eq:RVBL_prior}
    f_{\mathcal{X}|\mathcal{T}^2,\upeta^{-2}}(\bm x|\bm\tau^2,\eta^{-2})=f_{\mathcal{X}|\mathcal{T}^2}(\bm x|\bm\tau^2)\propto\frac{1}{|D(\bm\tau^2)|^{\frac12}}\exp\left(-\frac12(L\bm x)^T[D(\bm\tau^2)]^{-1}(L\bm x)\right)
\end{align}
as the prior along with the product of exponential densities 
\begin{align}\label{eq:RVBL_hyperprior}
    f_{\mathcal{T}^2|\upeta^{-2}}(\bm\tau^2|\eta^{-2})=\prod_{i=1}^{k}\frac{\eta^{-2}}{2}\exp\left(-\frac{\tau_i^2\eta^{-2}}{2}\right) 
\end{align}
as the hyperprior in \eqref{eqn:bayes}. Finally with $f_{\upeta^{-2}}(\eta^{-2})$ given in \eqref{eq:hyperprior},\footnote{For ease of presentation the rest of this manuscript uses the hyperparameter $\upeta^{-2}$ in place  of the original hyperparameter $\upeta$.} we are now able to completely characterize the hierarchical model from \eqref{eq:likelihood}, \eqref{eq:hyperprior}, \eqref{eq:RVBL_prior}, and \eqref{eq:RVBL_hyperprior}  as
\begin{subequations}
\label{eq:sparse_signal_variance_model}
\begin{align}
\bm y|\bm x &\sim \mathcal{N}(F\bm x,\sigma^2\mathbb{I})\label{eqn:hier_likelihood} \\
 \bm x|\bm\tau^2 &\sim \mathcal{N}(0,{ ( L^T[D(\bm\tau^{2})]^{-1}L)^{-1} }) \label{eqn:hier_prior} \\
\bm\tau^2|\eta^{-2} &\sim\prod_{j=1}^{k} \frac{\eta^{-2}}{2}\exp(-\frac{\tau_{j}^2\eta^{-2}}{2})\mathrm{d}\tau_{j}^2\label{eqn:hier_hyper} \\
\eta^{-2}&\sim\Gamma(r,\delta),
\label{eqn:hier_hypereta}
\end{align}
\end{subequations}
which yields the corresponding posterior density function
\begin{equation}
    f_{\mathcal{X},\mathcal{T}^2,\upeta^{-2}|\mathcal{Y}}(\bm x,\bm\tau^2,\eta^{-2}|\bm y) \propto f_{\mathcal{Y}|\mathcal{X}}(\bm y|\bm x)f_{\mathcal{X}|\mathcal{T}^2}(\bm x|\bm\tau^2)f_{\mathcal{T}^2|\upeta}(\bm\tau^2|\eta)f_{\upeta^{-2}}(\eta^{-2}),\label{eq:post_edge_L1Real}
    \end{equation}
where
\begin{eqnarray*}
f_{\mathcal{Y}|\mathcal{X}}(\bm y|\bm x)&\propto& \frac{1}{(\sigma^2)^\frac{n}{2}}\exp(-\frac{1}{2\sigma^2}\norm{ \bm y -F\bm x}_2^2)\\
f_{\mathcal{X}|\mathcal{T}^2}(\bm x|\bm\tau^2) & \propto  &\frac{1}{\sqrt{ |D(\bm\tau^2)|}}\exp\left(-\frac12(L\bm x)^T[D(\bm\tau^2)]^{-1}(L\bm x)\right)\\
f_{\mathcal{T}^2|\upeta}(\bm\tau^2|\eta) &\propto & \prod_{j=1}^{k}\eta^{-2}\exp(-\frac{\tau_{j}^2\eta^{-2}}{2}) \\
f_{\upeta^{-2}}(\eta^{-2}) &\propto &(\eta^{-2})^{r-1}\exp(-\delta\eta^{-2}).
\end{eqnarray*}
By isolating the parts of (\ref{eq:post_edge_L1Real}) that depend on $\bm x$ and defining
\begin{eqnarray*}
    G&:=&\frac{1}{\sigma^2}F^TF+L^T[D(\bm\tau^{2})]^{-1}L,\\
    \bar{\bm x}&:=&\frac{1}{\sigma^2}G^{-1}F^T\bm y,
\end{eqnarray*}
we obtain
\begin{align}
    f_{\mathcal{X}|\mathcal{Y},\mathcal{T}^2}(\bm x|\bm y,\bm\tau^2)
    &\propto\exp(-\frac12(\bm x-\bar{\bm x})^TG(\bm x-\bar{\bm x})).\label{eq:RVBL_x_cond}
\end{align}

\subsubsection*{Conditionals on $\mathcal{T}_{j}^2$}  When conditioned on all other variables {\it including} $\mathcal{T}_{-j}^2=\{\mathcal{T}_{j'}^2:j'=1,\dots,k, \ j'\neq j\}$, \eqref{eq:post_edge_L1Real} gives
\begin{align}
    f_{\mathcal{T}_{j}^2|\mathcal{T}_{-j}^2,\mathcal{X},\mathcal{Y},\upeta^{-2}}(\tau_{j}^2|\tau_{-j}^2,\bm x,\bm y,\eta^{-2})&\propto\left(\tau_{j}^2\right)^{-\frac12}\exp(-\frac{[L\bm x]_{j}^2}{2\tau_{j}^2}-\frac{\tau_{j}^2\eta^{-2}}{2}),\quad j = 1,\dots,k.
    \label{eq:tau_posterior_real}
\end{align}
Making the change of variables $\nu_{j}^2=1/\tau_{j}^2$, the conditional can alternatively be expressed as 
\begin{align}\label{eq:inversegauss_edgeL1}
f_{\mathcal{V}_{j}^2|\mathcal{T}_{-j}^2,\mathcal{X},\upeta^{-2}}(\nu_{j}^2|\tau_{-j}^2,\bm x,\eta^{-2})&\propto\left(\nu_{j}^2\right)^{-\frac32}\exp(-\frac{[L\bm x]_{j}^2\nu_{j}^2}{2}-\frac{\eta^{-2}}{2\nu_{j}^2}).
\end{align}
By defining mean parameter $\mu' = \sqrt{  \eta^{-2} / [L \bm x]_j^2 }$ and shape parameter $\lambda' = \eta^{-2}$, we observe that \cref{eq:inversegauss_edgeL1} fits the form of the density function of inverse Gaussian distribution given by
\begin{align*}
f_{\mu',\lambda'}(y) \propto y^{-3/2} \exp( - \frac{\lambda' y}{2 (\mu')^2} - \frac{\lambda'}{2 y}   )
\propto y^{-3/2} \exp( - \frac{\lambda' (y - \mu')^2}{2 (\mu')^2 y} )\quad  \mu', \lambda' > 0.
\end{align*}
Thus \cref{eq:inversegauss_edgeL1} can be directly sampled using a standard routine for sampling the inverse Gaussian or Wald distribution (e.g., see \cite{Michael1976InverseGaussian}).

\begin{remark}\label{rem:gammaforconditional}
In the case where $[L \bm x]_j \approx 0$ (where $L$ could be the identity matrix $\mathbb{I}_n$), we can refrain from performing the change of variables to obtain \eqref{eq:inversegauss_edgeL1} and note that the original conditional (\ref{eq:tau_posterior_real}) simplifies to
\begin{align}
\label{eq:Lis0condition}
f_{\mathcal{T}_{j}^2|\mathcal{T}_{-j}^2,\mathcal{X},\mathcal{Y},\upeta^{-2}}(\tau_{j}^2|\tau_{-j}^2,\bm x,\bm y,\eta^{-2})&\propto\left(\tau_{j}^2\right)^{-\frac12}\exp(-\frac{\tau_{j}^2\eta^{-2}}{2}),
\end{align}
 corresponding to the Gamma  distribution $\Gamma(1/2, \eta^{-2}/2)$. Thus whenever $|[L \bm x]_j| < 10^{-8}$ in our numerical experiments we instead sample the conditional density of each $\mathcal{T}_j^2$ according to \cref{eq:Lis0condition}.\footnote{Due to system noise, nearly all values $|[L \bm x]_j|$ are greater than $10^{-8}$ with typical values in $[.001,.1]$.}
\end{remark}

\subsubsection*{Conditional on $\upeta^{-2}$} Lastly, to sample the conditional on $\upeta^{-2}$, we note that
\begin{align}\label{eq:eta_conditional}
    f_{\upeta^{-2}|\mathcal{X}, \mathcal{Y}, \mathcal{T}^2}(\eta^{-2}|\bm x, \bm y, \bm \tau^2)
    &\propto (\eta^{-2})^{r + k - 1} \exp\left( - \eta^{-2} \left(\delta + \frac{1}{2} \sum_{j=1}^k \tau_j^2 \right) \right),
\end{align}
which due to conjugacy is the density for the Gamma distribution  $\Gamma(r+k, \delta + \frac12\sum_{j=1}^k \tau_j^2 )$. Algorithm \ref{alg:gibbs_laplace} outlines the full Gibbs sampling approach.

\begin{algorithm}[H]
	\caption{Real-valued Bayesian LASSO (RVBL)}
	\label{alg:gibbs_laplace}
	\begin{algorithmic}
		\item[] {\bf Input} Data vector $\bm y$, 
  parameters $\sigma$ and $\eta$, sparsifying operator $L \in \mathbb{R}^{k \times n}$, chain length $N_M$, and burn-in length $B$.
        \item[] {\bf Output} Samples $\bm x^{(s-B+1)}$ for $s=B,\dots,N_M$.
		\item[1] Set $\bm x^{(0)}=F^T\bm y$ and $\left(\tau_{j}^2\right)^{(0)}=\left(\eta^{-2}\right)^{(0)}=1$ for $j=1,\dots,k$.
		\item[2] {\bf For} $l=1,\dots,N_M$ {\bf do}
		\begin{itemize}
		    \item[i.] Sample $\bm x^{(l)}$ from $\mathcal{X}\Big|\mathcal{T}^2=\left(\bm \tau^{2}\right)^{(l-1)}$
            \item[ii.] Sample each $\left(\tau_{j}^{2}\right)^{(l)}$ using $\mathcal{T}_{j}^{-2}\Big|\mathcal{X}=\bm x^{(l)},\upeta^{-2}=\left(\eta^{-2}\right)^{(l-1)}$.
            \item[iii.] Sample $\left(\eta^{-2}\right)^{(l)}$ using 
            $\upeta^{-2}|\mathcal{T}^2=\left(\bm \tau^{2}\right)^{(l-1)}$.
		\end{itemize}
	\end{algorithmic}
\end{algorithm}

\begin{figure}[h!]
    \centering
    \begin{subcaptionblock}{.28\textwidth}
    \centering
    \includegraphics[width=.8\linewidth]{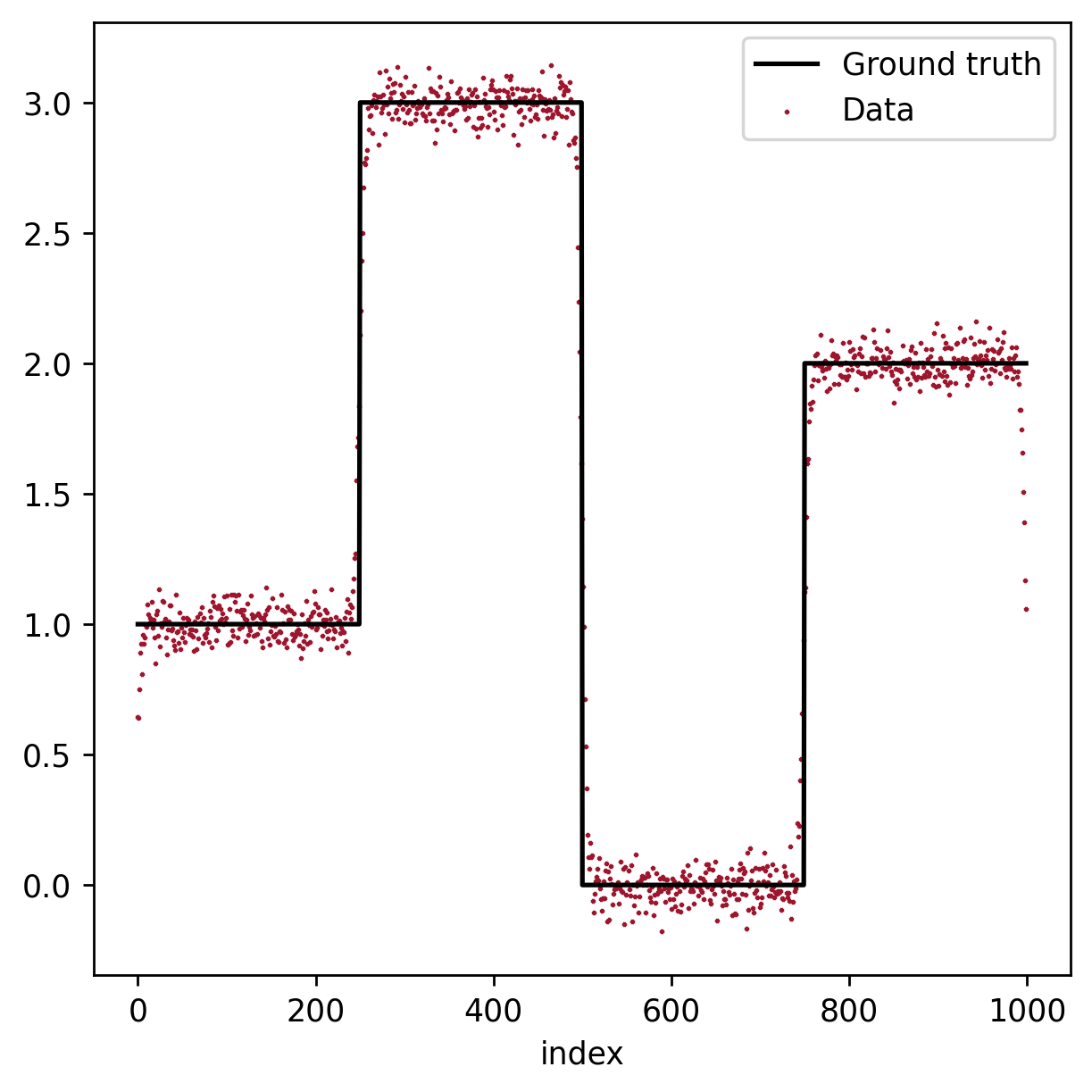}
    \caption{Deblurring problem data}\label{fig:rvbl_data}
    \end{subcaptionblock}
    \begin{subcaptionblock}{.28\textwidth}
    \centering
    \includegraphics[width=.8\linewidth]{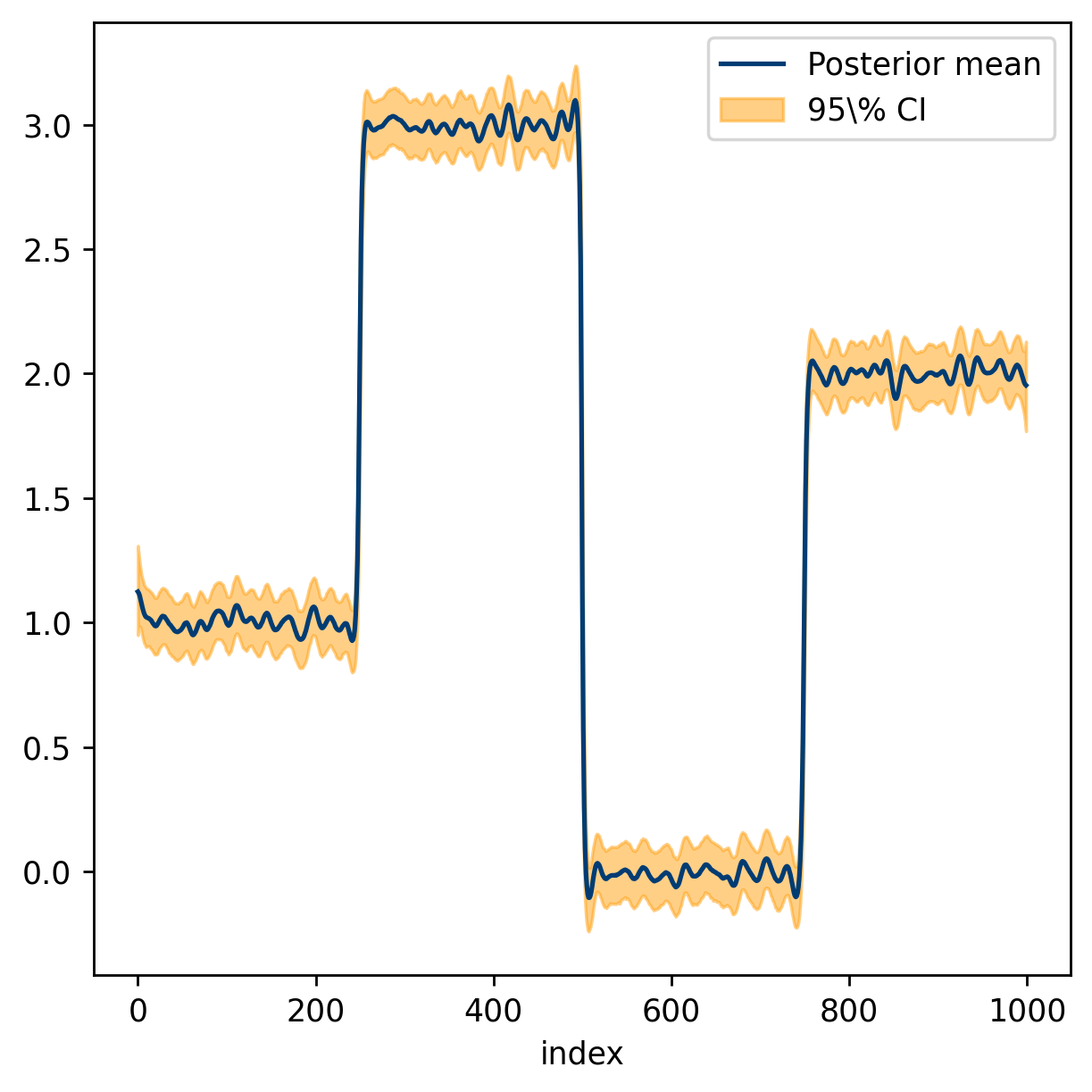}
    \caption{Estimate of $\bm x$}\label{fig:rvbl_x_result}
    \end{subcaptionblock}
    \begin{subcaptionblock}{.28\textwidth}
    \centering
    \includegraphics[width=.8\linewidth]{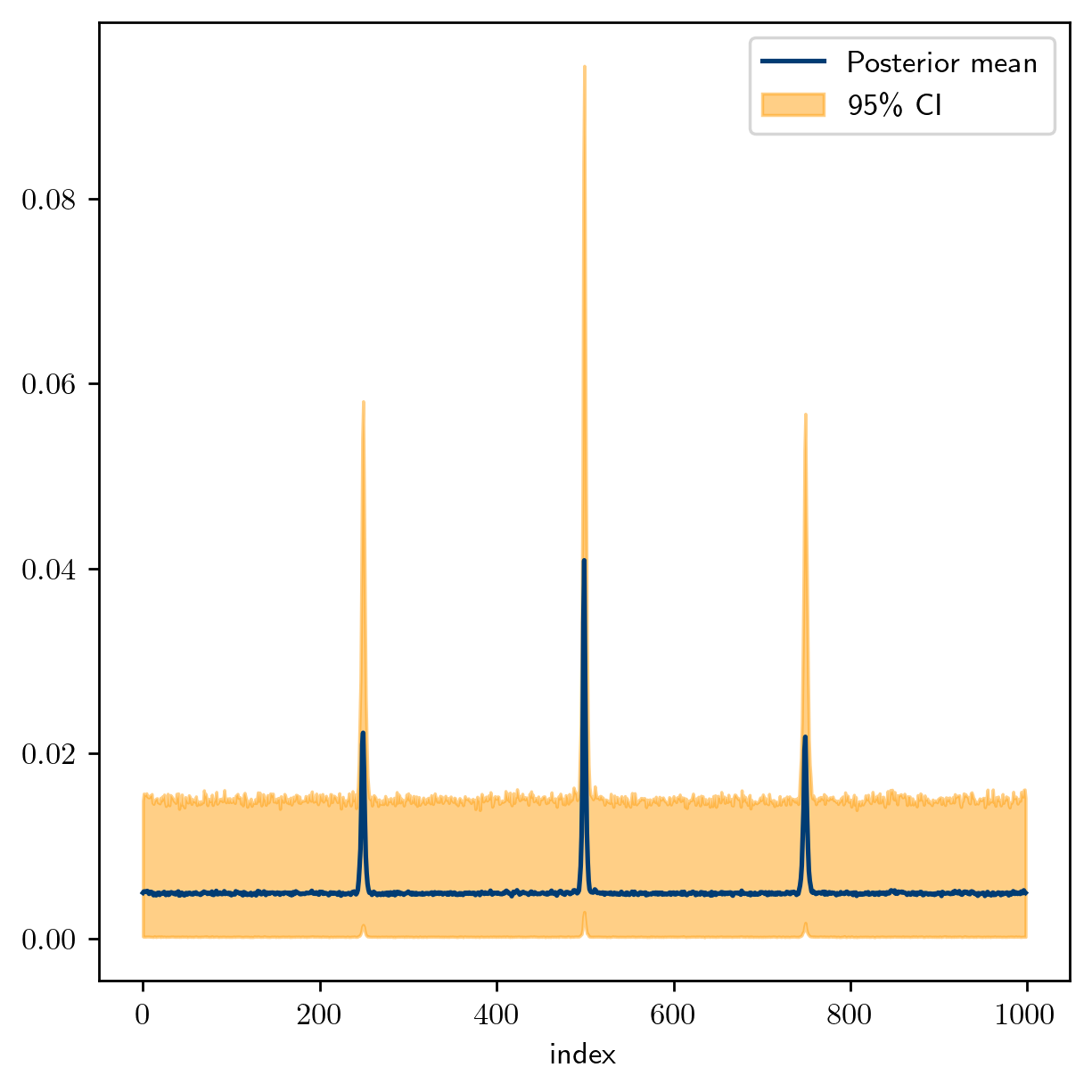}
    \caption{Estimate of $\bm \tau^2$}\label{fig:rvbl_theta_result}
    \end{subcaptionblock}
\caption{The RVBL method applied to a signal deblurring problem. The posterior mean estimate for the parameter $\upeta^{-2}$ is $\bar{\eta}^{-2} \approx 401$ 
and SNR = 30. in \eqref{eq:SNR}.}\label{fig:real_valued_lasso}
\end{figure}

We conclude this section with illustrative numerical example in \cref{fig:real_valued_lasso} of the RVBL method applied to a signal deblurring problem, where $F$ is a Gaussian blurring kernel with standard deviation $\sigma_{\text{blur}} = 5.0$ and with a zero boundary condition. We assume that the data $\bm y$ are generated by $\bm y = F\bm x  +\bm\varepsilon $, where $\bm\varepsilon$ is Gaussian white noise with SNR $= 30\dB$, where SNR is defined in \eqref{eq:SNR} for complex-valued signals. The problem data and ground truth signal are shown in \cref{fig:real_valued_lasso} (a). To apply the RVBL method, we take the sparsifying transformation $L$ to be the discrete gradient operator \eqref{eq:Ldiff} and adopt the uninformative hyperprior  given by \eqref{eq:hyperprior}. \cref{fig:real_valued_lasso} (b)-(c) show the resulting posterior mean estimates for $\bm x$ and $\bm \tau^2$, as well as their 95\% credible intervals.

\section{Complex-valued Bayesian LASSO (CVBL)} \label{sec:proposed_methods}
We now have all of the ingredients needed to compute and sample from the complex-valued posterior density functions using the likelihood given in \eqref{eq:likelihood}, the priors described in \eqref{eq:prior_signal} and \eqref{eq:prior_edge}, and the hyperprior \eqref{eq:hyperprior}. We will describe the CVBL  for two distinct cases:  \Cref{subsec:sparseL1} considers the magnitude itself to be sparse, while \Cref{subsec:edgeL1} examines the case for which the magnitude is sparse in the some transform domain.

The scale mixture \eqref{eq:scalemix} used in the RVBL relies on the likelihood being a real-valued probability density function. As such, for our algorithmic development we will make use of the equivalency that for given $\bm y\in\mathbb{C}^m$, $F\in\mathbb{C}^{m\cross n}$, and $\bm x\in\mathbb{R}^n$, we have
\begin{align}
\label{eq:mathfact}
    \norm{\bm y-F\bm x}_2^2 = \norm{\tilde{\bm y}-\tilde{F}\bm x}_2^2,
\end{align}
where $\tilde{\bm y} = [\Re(\bm y)^T \ \Im(\bm y)^T]^T$ and $\tilde{F} = [\Re(F)^T \ \Im(F)^T]^T$.

\subsection{CVBL for the Sparse Magnitude Case}\label{subsec:sparseL1}
We begin by writing $\mathcal{Z}$ as
\begin{align*}
\mathcal{Z}=\mathcal{A}+i\mathcal{B},
\end{align*}
where $\mathcal{Z}\in\mathbb{C}^n$ and $\mathcal{A},\mathcal{B}\in\mathbb{R}^n$, with realizations $\bm z=\bm a+i\bm b$.  When $|\mathcal{Z}| = |\mathcal{A}+i\mathcal{B}|$ is assumed to be sparse, we choose to impose a $1-$norm prior on $\sqrt{\mathcal{A}^2+\mathcal{B}^2}$ using \eqref{eq:prior_signal}. We now provide the details for how this is accomplished.

For observations given in \eqref{eqn:observation_model}, an analogous derivation that yielded \eqref{eq:sparse_signal_variance_model}  results in the hierarchical model for the {\em complex-valued} sparse signal recovery
\begin{align*}
\bm y|\bm a,\bm b&\sim \mathcal{CN}(F(\bm a+i\bm b),\sigma^2\mathbb{I}),\\
\bm a|\bm\tau^2&\sim \mathcal{N}(0,D(\bm\tau^{-2})),\\
\bm b|\bm\tau^2&\sim \mathcal{N}(0,D(\bm\tau^{-2})),\\
\bm\tau^2|\eta^{-2} &\sim\prod_{j=1}^k \frac{\eta^{-2}}{2}\exp(-\frac{\tau_{j}^2\eta^{-2}}{2})\mathrm{d}\tau_{j}^2,\label{eqn:hier_hyper_comp} \\
\eta^{-2} &\sim\Gamma(r,\delta),
\end{align*}
where $\bm\tau^2$ (defined componentwise as in \eqref{eq:vector_square}) is a realization of the scale mixture parameter $\mathcal{T}^2$ introduced in \eqref{eqn:joint_to_marginal}.  We can choose hyperparameters $r$ and $\delta$ as in \eqref{eq:r_and_delta} to promote sparsity.  The posterior is then written as
\[f_{\mathcal{Z},\mathcal{T}^2,\upeta^{-2}|\mathcal{Y}}(\bm z,\bm\tau^2,\eta^{-2}|\bm y)
= 
f_{\mathcal{A},\mathcal{B},\mathcal{T}^2,\upeta^{-2}|\mathcal{Y}}(\bm a,\bm b,\bm\tau^2,\eta^{-2}|\bm y),\]
with
\begin{equation}
f_{\mathcal{A},\mathcal{B},\mathcal{T}^2,\upeta^{-2}|\mathcal{Y}}(\bm a,\bm b,\bm\tau^2,\eta^{-2}|\bm y) \propto f_{\mathcal{Y}|\mathcal{A},\mathcal{B}}(\bm y|\bm a,\bm b)f_{\upeta^{-2}}(\eta^{-2})\prod_{j=1}^k f_{\mathcal{A}_{j},\mathcal{B}_j|\mathcal{T}_{j}^2}(a_j,b_j|\tau_{j}^2)f_{\mathcal{T}_{j}^2|\upeta^{-2}}(\tau_{j}^2|\eta^{-2}).\label{eq:posterior_signalL1}
\end{equation}
Here
\begin{eqnarray*}
 f_{\mathcal{Y}|\mathcal{A},\mathcal{B}}(\bm y|\bm a,\bm b)  & \propto & \frac{1}{(\sigma^2)^n}\exp(-\frac{1}{\sigma^2}\norm{\bm y-F\bm a-iF\bm b}_2^2),\\
 f_{\upeta^{-2}}(\eta^{-2}) &\propto& (\eta^{-2})^{r-1}\exp(-\delta\eta^{-2}),\\
f_{\mathcal{A}_{j},\mathcal{B}_j|\mathcal{T}_{j}^2}(a_j,b_j|\tau_{j}^2)  & \propto & \frac{1}{\sqrt{\tau_{j}^2}}\exp(-\frac{a_{j}^2+b_j^2}{2\tau_{j}^2}),\\
 f_{\mathcal{T}_{j}^2|\upeta^{-2}}(\tau_{j}^2|\eta^{-2}) &\propto& \eta^{-2}\exp(-\frac{\tau_{j}^2\eta^{-2}}{2}).
\end{eqnarray*}

For ease of notation we define
\begin{equation}
\label{eq:splitting}
\tilde{F}=[\Re(F)^T \ \Im(F)^T]^T, \quad \tilde{F}^*=[-\Im(F)^T \ \Re(F)^T]^T,\quad\tilde{\bm y}=[\Re(\bm y)^T \ \Im(\bm y)^T]^T, 
\end{equation}
so that the real and imaginary parts of the observations are respectively
\begin{equation}
    \bm y_1=\tilde{\bm y}-\tilde{F}^*\bm b \quad\text{ and }\quad
\bm y_2=\tilde{\bm y}-\tilde{F}\bm a. \label{eq:splitobservationdata}
\end{equation}
\subsubsection*{Sampling on $\mathcal{A}$ and $\mathcal{B}$}
Analogous to the procedure resulting in \eqref{eq:RVBL_x_cond}, we now sample from the conditionals on $\mathcal{A}$ and $\mathcal{B}$ by first isolating the parts of (\ref{eq:posterior_signalL1}) that depend on $\bm a$ and $\bm b$ respectively, yielding
\begin{subequations}
\begin{equation}
f_{\mathcal{A}|\mathcal{Y},\mathcal{B},\mathcal{T}^2}(\bm a|\bm y,\bm b,\bm\tau^2)\propto\exp(-\frac12(\bm a-\bar{\bm a})^TG(\bm a-\bar{\bm a})), \label{eq:condA_signalL1}
\end{equation}
\begin{equation}
f_{\mathcal{B}|\mathcal{Y},\mathcal{A},\mathcal{T}^2}(\bm b|\bm y,\bm a,\bm\tau^2)\propto\exp(-\frac12(\bm b-\bar{\bm b})^TG(\bm b-\bar{\bm b})), \label{eq:condB_signalL1}
\end{equation}
\end{subequations}
where $G=\frac{2}{\sigma^2}\tilde{F}^T\tilde{F}+[D(\bm\tau^{2})]^{-1}$, $\bar{\bm a}=\frac{2}{\sigma^2}G^{-1}\tilde{F}^T\bm y_1$ and  $\bar{\bm b}=\frac{2}{\sigma^2}G^{-1}\tilde{F}^{*T}\tilde{\bm y}_2$. Observe that both  $\mathcal{A}$ and $\mathcal{B}$ are  conditionally Gaussian.

\begin{remark}
    If $F$ is unitary, for example when it represents the $n \times n$ normalized discrete Fourier transform, then \eqref{eq:condA_signalL1} and \eqref{eq:condB_signalL1} can be simplified to multivariate Gaussian densities with respective means
    \begin{gather*}
        \bm a = \frac{2}{\sigma^2}G^{-1}\Re(F^H(\bm y-iF\bm b)), \quad
        \bm b = \frac{2}{\sigma^2}G^{-1}\Im(F^H(\bm y-F\bm a)),
    \end{gather*}
    where
    $G=\frac{2}{\sigma^2}\mathbb{I}+[D(\bm \tau^2)]^{-1}.$
    The splitting in \eqref{eq:splitting} is needed in applications where the input data are under-sampled, or when portions of the data must be discarded.
\end{remark}

\subsubsection*{Sampling on scale mixing parameter $\mathcal{T}_j^2$}
Obtaining the conditional density of $\{\mathcal{T}_{j}^2\}_{j = 1}^k$ from \eqref{eq:posterior_signalL1} yields
\begin{align}
    f_{\mathcal{T}_{j}^2|\mathcal{A},\mathcal{B},\mathcal{T}_{-j}^2,\upeta^{-2}}(\tau_{j}^2|\bm a,\bm b,\tau_{-j}^2,\eta^{-2})&\propto\left(\tau_{j}^2\right)^{-\frac12}\exp(-\frac{a_{j}^2+b_j^2}{2\tau_{j}^2}-\frac{\tau_{j}^2\eta^{-2}}{2}), \quad j = 1,\dots,k.
    \label{eq:tau_cond_signal}
\end{align}
By change of variable $\mathcal{V}_j^2= (\mathcal{T}_j^{2})^{-1}$, with realization $\nu_{j}^2=(\tau_{j}^{2})^{-1}$,  we obtain
\begin{align}
    f_{\mathcal{V}_{j}^2|\mathcal{A},\mathcal{B},\mathcal{T}_{-j}^2,\upeta^{-2}}(\nu_{j}^2|\bm a,\bm b,\tau_{-j}^2,\eta^{-2})&\propto\left(\nu_{j}^2\right)^{-\frac32}\exp(-\frac{1}{{2\nu_{j}^2}}\left(a_{j}^2+b_j^2\right)\left(\nu_{j}^2-\sqrt{\frac{\eta^{-2}}{(a_{j}^2+b_j^2)}}\right)^2).\label{eq:invgauss_signalL1}
\end{align}
Observe that comparable to \eqref{eq:inversegauss_edgeL1}, \eqref{eq:invgauss_signalL1} is the probability density function for an inverse Gaussian distribution with mean parameter $\mu'=\sqrt{\eta^{-2}/(a_j^2+b_j^2)}$ and shape parameter $\lambda'=\eta^{-2}$. Moreover, since each $(\mathcal{T}_{j}^2)^{-1}$ is mutually independent from $(\mathcal{T}_{-j}^2)^{-1}$, $j=1,\dots,k$,  we can efficiently sample the conditionals on each $(\mathcal{T}_{j}^2)^{-1}$ in parallel. 

\begin{remark}
Similar to the discussion in \cref{rem:gammaforconditional}, we observe that the density in \eqref{eq:invgauss_signalL1} is not well-defined when $a_j = b_j = 0$. Although the probability of sampling such $a_j$ and $b_j$ is zero, the probability of sampling $a_j$ and $b_j$ within machine precision is not. Thus in the case where $a_j$ and $b_j$ are sufficiently small, we instead use $\mathcal{T}_{j}^2|\mathcal{A},\mathcal{B},\mathcal{T}_{-j}^2,\upeta^{-2} \sim \Gamma(1/2,\eta^{-2}/2)$ following \eqref{eq:tau_cond_signal}.
\end{remark}

\subsubsection*{Sampling on hyperparameter $\upeta^{-2}$}
The conditional posterior of $\upeta^{-2}$ depends exclusively on $\mathcal{T}^2$, hence $f_{\upeta^{-2}|\mathcal{A},\mathcal{B},\mathcal{Y},\mathcal{T}^2}=f_{\upeta^{-2}|\mathcal{T}^2}$. This density then has the form of \eqref{eq:eta_conditional}, from which we are able to sample directly. 

\subsection*{CVBL algorithm for sparse magnitude signals} The conditional distributions in \eqref{eq:condA_signalL1}, \eqref{eq:condB_signalL1}, and \eqref{eq:invgauss_signalL1} are now combined to form the CVBL Gibbs sampling method provided in Algorithm \ref{alg:gibbs_signalL1}. Methods to efficiently sample $\bm a$ and $\bm b$ in \cref{alg:gibbs_signalL1} are discussed in \Cref{subsec:numerics_efficiency}.

\begin{algorithm}[H]
	\caption{CVBL for complex-valued signal with sparse magnitude}
	\label{alg:gibbs_signalL1}
	\begin{algorithmic}
		\item[] {\bf Input} Data vector $\bm y$, noise variance $\sigma^2$, hyperparameters $r=1$ and $\delta=10^{-3}$, chain length $N_M$, and burn-in length $B$.
        \item[] {\bf Output} samples $\bm z^{(s-B+1)}=\bm a^{(s)}+ i\bm b^{(s)}$ for $s=B,\dots,N_M$.
		\item[1] Set $\bm a^{(0)}=\Re(A^H\bm y)$, $\bm b^{(0)}=\Im(A^H\bm y)$, and $\left(\tau_{j}^2\right)^{(0)}=\left(\eta^{-2}\right)^{(0)}=1$ for $j=1,\dots,k$.
		\item[2] {\bf For} $\ell=1,\dots,N_M$ {\bf do}
		\begin{itemize}
		    \item[i.] Sample $\bm a^{(\ell)}$ from $\mathcal{A}|\mathcal{Y}=\bm y,\mathcal{B}=\bm b^{(\ell-1)},\mathcal{T}^2=\left(\bm \tau^{2}\right)^{(\ell-1)}$ \eqref{eq:condA_signalL1}.
            \item[ii.] Sample $\bm b^{(\ell)}$ from $\mathcal{B}|\mathcal{Y}=\bm y,\mathcal{A}=\bm a^{(\ell)},\mathcal{T}^2=\left(\bm \tau^{2}\right)^{(\ell-1)}$ \eqref{eq:condB_signalL1}.
            \item[iii.] Sample $\left(\tau_{j}^{-2}\right)^{(\ell)}$ from $\mathcal{T}_{j}^{-2}|\mathcal{A}=\bm a^{(\ell)},\mathcal{B}=\bm b^{(\ell)},\upeta^{-2}=\left(\eta^{-2}\right)^{(\ell-1)}$ \eqref{eq:invgauss_signalL1} for  $j=1,\dots,n$.
            \item[iv.] Sample $\left(\eta^{-2}\right)^{(\ell)}$ using 
            $\upeta^{-2}|\mathcal{T}^2=\left(\bm \tau^{2}\right)^{(\ell)}$ \eqref{eqn:hier_hypereta}.
		\end{itemize}
	\end{algorithmic}
\end{algorithm}

\subsection{Sparsity in a Transform Domain of $|\mathcal{Z}|$}\label{subsec:edgeL1}
We now turn our attention to the case where sparsity is expected in some transform domain of $|\mathcal{Z}|$.  For simplicity  we assume the signal magnitude is piecewise constant, so that there is sparsity in the corresponding gradient domain. Hence we define  $L$ in \eqref{eq:prior_edge} to be the TV operator (see \Cref{sec:sparseoperator} for more discussion).

Starting from \eqref{eq:splitmagphase} and assuming  $f_{\mathcal{Y}}(\bm y)>0$, Bayes' Theorem yields
\begin{align}
f_{\mathcal{G},\Phi,\upeta|\mathcal{Y}}(\bm g,\bm \phi,\eta|\bm y)
	&\propto\exp\left(-\frac{1}{\sigma^2}\norm{\bm y- A\left(\bm g\odot e^{i\bm \phi}\right)}^2_2-\frac{1}{\eta}\norm{L\bm g}_1\right)\bm1_{\mathbb{R}^+}(\bm g)\bm1_{[-\pi,\pi)}(\bm \phi).
 \label{eq:post1norm}
\end{align}
Analogously rewriting the prior as a scale mixture of normals as in \eqref{eqn:joint_to_marginal} provides the posterior density function
\begin{align}
    f_{\mathcal{G},\Phi,\mathcal{T}^2,\upeta^{-2}|\mathcal{Y}}(\bm g,\bm\phi,&\bm\tau^2,\eta^{-2}|\bm y)\propto f_{\mathcal{Y}|\mathcal{G},\Phi}(\bm y|\bm g,\bm\phi)f_\Phi(\bm\phi)f_{\mathcal{G}|\mathcal{T}^2}(\bm g|\bm\tau^2)f_{\upeta^{-2}}(\eta^{-2})\prod_{j=1}^kf_{\mathcal{T}_{j}^2|\upeta^{-2}}(\bm\tau_{j}^2|\eta^{-2}),
    \label{eq:posterior_edgeL1}
\end{align}
where
\begin{eqnarray*}
 f_{\mathcal{Y}|\mathcal{G},\Phi}(\bm y|\bm g,\bm\phi)  & \propto & \exp(-\frac{1}{\sigma^2}\norm{y-F(\bm g\odot e^{i\bm\phi})}_2^2)\bm1_+(\bm g),\\
 f_\Phi(\bm\phi) &\propto& \bm1_{[-\pi,\pi)}(\bm\phi),\\
f_{\mathcal{G}|\mathcal{T}^2}(\bm g|\bm\tau^2)  & \propto & \frac{1}{\sqrt{|D(\bm\tau^2)|}}\exp\left(-\frac12(L\bm g)^T[D(\bm\tau^2)]^{-1}(L\bm g)\right),\\
 f_{\upeta^{-2}}(\eta^{-2}) &\propto& (\eta^{-2})^{r-1}\exp(-\delta\eta^{-2}),\\
 f_{\mathcal{T}_{j}^2|\upeta^{-2}}(\bm\tau_{j}^2|\eta^{-2})&\propto&\eta^{-2}\exp(-\frac{\tau_{j}^2\eta^{-2}}{2}).
\end{eqnarray*}

Below we introduce a four-step Gibbs sampling process to sample from \eqref{eq:posterior_edgeL1} that combines (1) the use of Metropolis-within-Gibbs to sample from $\mathcal{G}|\mathcal{Y},\Phi,\mathcal{T}^2,\upeta^{-2}$; (2) Gibbs steps to sample from $\mathcal{T}^2|\mathcal{Y},\mathcal{G},\Phi,\upeta^{-2}$ and $\upeta^{-2}|\mathcal{Y},\mathcal{G},\Phi,\mathcal{T}^2$; and (3) rejection sampling to sample from $\Phi|\mathcal{Y},\mathcal{G},\mathcal{T}^2,\upeta^{-2}$.

\subsubsection*{Sampling on magnitude $\mathcal{G}$} The conditional distribution on $\mathcal{G}$ is a nonnegatively-constrained Gaussian density. For computational efficiency in sampling the magnitude, we use as the posterior density the untruncated density $\tilde{f}_{\mathcal{G}|\mathcal{Y},\mathcal{T}^2,\Phi}(\bm g|\bm y,\bm\tau^2,\bm \phi)$, where
\begin{align*}
    f_{\mathcal{G}|\mathcal{Y},\Phi,\mathcal{T}^2}(\bm g|\bm y,\bm \phi,\bm\tau^2)=\tilde{f}_{\mathcal{G}|\mathcal{Y},\Phi,\mathcal{T}^2}(\bm g|\bm y,\bm \phi,\bm\tau^2)\bm1_+(\bm g).
\end{align*}
Thus
\begin{align}
    \tilde{f}_{\mathcal{G}|\mathcal{Y},\Phi,\mathcal{T}^2}(\bm g|\bm y,\bm \phi,\bm\tau^2)\propto\exp\left(-\frac{\norm{\bm y-F_1\bm g}^2_2}{\sigma^2}-\frac{1}{2}(L\bm g)^T[D(\bm\tau^2)]^{-1}(L\bm g)\right),
	\label{eq:magpost}
\end{align}
where $F_1=FD(e^{i\bm \phi})$. For consistency with \eqref{eq:posterior_edgeL1}, we simply reject any proposed samples of $\mathcal{G}$ that contain a negative entry. This technique can be categorized as acceptance-rejection sampling, and to this end, \Cref{thm:magpost} shows that \eqref{eq:magpost} is a conditional Gaussian distribution, implying that we can directly draw its samples.

\begin{theorem}\label{thm:magpost}
Let  $F_1 = FD(e^{i\phi})$, $\tilde{F}_1=[\Re(F_1)^T \ \Im(F_1)^T]^T$, and $\tilde{\bm y}=[\Re(\bm y)^T \ \Im(\bm y)^T]^T$. Assume that $L$ has rank $n$. The function $\tilde{f}_{\mathcal{G}|\mathcal{Y},\Phi,\mathcal{T}^2}(\bm g|\bm y,\bm \phi,\bm\tau^2)$ in \eqref{eq:magpost} defines a Gaussian density over $\mathbb{R}^n$ with center and precision given by
\begin{subequations}
\label{eq:centerprecision}
    \begin{equation}
	\bar{\bm g} = \Gamma^{-1}\left(\frac{2}{\sigma^2}\tilde{F}_1^T\tilde{\bm y}\right),\label{eq:barg}
 \end{equation}
 \begin{equation}
	\Gamma = L^T[D(\bm\tau^2)]^{-1}L+\frac{2}{\sigma^2}\tilde{F}_1^T\tilde{F}_1.\label{eq:Gamma}
 \end{equation}
\end{subequations}
\end{theorem}
The proof of \cref{thm:magpost} is provided in \cref{sec:appendix_proofs}.

When $F$ is unitary, the resulting Gaussian density attains a simpler form provided by \Cref{thm:Funitary} which enables faster sampling.

\begin{corollary} \label{thm:Funitary}
Let $F_1 = FD(e^{i\phi})$ and suppose $F$ is unitary and $L$ has rank $n$.  For $F_1^H\bm y=\bm q\odot e^{i\bm\varphi}$ where $\bm q\in\left(\mathbb{R}^+\right)^n$ and $\bm\varphi\in[-\pi,\pi)^n$, the function $\tilde{f}_{\mathcal{G}|\mathcal{Y},\mathcal{T}^2,\Phi}(\bm g|\bm y,\bm\tau^2,\bm \phi)$ in \eqref{eq:magpost} defines a Gaussian density over $\mathbb{R}^n$ with center and precision given by
\begin{subequations}
\label{eq:precisioncenter_unity}
    \begin{equation}
	\bar{\bm g} = \Gamma^{-1}\left(\frac{2}{\sigma^2}\bm q\odot\cos(\bm\varphi)\right),\label{eq:barg_unitary}
 \end{equation}
 \begin{equation}
	\Gamma = L^T[D(\bm\tau^2)]^{-1}L+\frac{2}{\sigma^2}\mathbb{I}_n.\label{eq:Gamma_unitary}
 \end{equation}
\end{subequations}
\end{corollary}
The proof of \cref{thm:Funitary} can be found in \cref{sec:appendix_proofs}. A couple of remarks are in order.

\begin{remark}\label{rem:little_pos_mass}
If the density \eqref{eq:magpost} has little mass in the region where $\bm g\geq0$, this rejection method may become computationally infeasible. Other techniques exist for sampling from truncated multivariate normal distributions, such as the exact Hamiltonian Monte Carlo method in \cite{pakman2014exact} or a Gibbs sampling technique that updates each $g_i$ separately for $i=1,\dots,n$. These methods are often more computationally expensive and may become cost-prohibitive for high dimensional problems, however.
\end{remark}

\begin{remark}
    When $\bar{\bm g}$ in \eqref{eq:barg} and \eqref{eq:barg_unitary} has mostly nonnegative entries with a few negative elements, sampling $f_{\mathcal{G}|\mathcal{Y},\Phi,\mathcal{T}^2}$ using the rejection technique involving $\tilde{f}_{\mathcal{G}|\mathcal{Y},\Phi,\mathcal{T}^2}$ may become computationally inefficient. In this case the rejection sampling from the mode (RSM) \cite{maatouk2016new} provides a possible alternative. In short, RSM generates an exact sample of $f_{\mathcal{G}|\mathcal{Y},\Phi,\mathcal{T}^2}$ by sampling a shifted truncated Gaussian density followed by an acceptance-rejection step. In our numerical experiments RSM is implemented when the rejection technique involving \eqref{eq:magpost} fails to generate a sample after $N_s$ attempts (in our experiments  $N_s=10$).
\end{remark}

\subsubsection*{Updating the scale mixture parameter $\mathcal{T}^2$}
From \eqref{eq:posterior_edgeL1} we have the conditional posterior 
\begin{align}
f_{\mathcal{T}_{j}^2|\mathcal{Y},\mathcal{G},\Phi,\mathcal{T}_{-j}^2,\upeta^{-2}}(\tau_{j}^2|\bm y,\bm g,\bm\phi,\tau_{-j}^2,\eta^{-2})&\propto\left(\tau_{j}^2\right)^{-\frac12}\exp(-\frac{[L\bm g]_{j}^2}{2\tau_{j}^2}-\frac{\tau_{j}^2\eta^{-2}}{2}),\quad{j = 1,\dots,k}, \label{eq:tau_posterior_complex}
\end{align}
which is equivalent to the density function in the real-valued signal case
\eqref{eq:tau_posterior_real}  with $\mathcal{X} = \mathcal{G}$. Analogously to what followed there, we recognize \eqref{eq:tau_posterior_complex} as an inverse Gaussian distribution with mean parameter $\sqrt{\eta^{-2}/[L\bm g]_j^{2}}$ and shape parameter $\eta^{-2}$. The conditional posterior on $\upeta^{-2}$ is  given by \eqref{eq:eta_conditional}.

\subsubsection*{Updating the phase $\Phi$} Lastly from \eqref{eq:posterior_edgeL1} we have the phase posterior distribution
\begin{align}\label{eq:edge_postphi}
	f_{\Phi|\mathcal{Y},\mathcal{G},\mathcal{T}^2}(\bm \phi|\bm y,\bm g,\bm\tau^2)\propto\exp\left(-\frac{1}{\sigma^2}\norm{\bm y-F_2e^{i\bm \phi}}^2_2\right)\bm1_{[-\pi,\pi)}(\bm \phi),
\end{align}
where $F_2=FD(\bm g)$. The nature of \eqref{eq:edge_postphi} can be better understood by defining random variable $\Theta=e^{i\Phi}$ and corresponding realization $\bm \theta=e^{i\bm \phi}$. The posterior is then expressed as
\begin{align}\label{eq:edge_posttheta}
	f_{\Theta|\mathcal{Y},\mathcal{G}}(\bm \theta|\bm y,\bm g)\propto\exp\left(-\frac{1}{\sigma^2}\norm{\bm y-F_2\bm \theta}^2_2\right)\bm1_{CS_1}(\bm \theta),
\end{align}
where $CS_1$ is the unit circle in the complex plane. Clearly \eqref{eq:edge_posttheta}, and by extension \eqref{eq:edge_postphi}, are probability density functions of  complex Gaussian distributions restricted to the unit circle. For a general forward operator $F$, each $f_{\Phi_i|\Phi_{-i},\mathcal{Y},\mathcal{G},\mathcal{T}^2}(\phi_i|\phi_{-i},\bm y,\bm g,\bm\tau^2)$ is conditionally von Mises \eqref{eq:vonmises_def} for $i=1,\dots,n$, which we now state in \cref{thm:cond_vonmises}.

\begin{theorem}\label{thm:cond_vonmises}
    Let $F_2=FD(\bm g)$ as in \eqref{eq:edge_postphi} and \eqref{eq:edge_posttheta} and define $A=F_2^HF_2$ with elements $[A]_{j,k}=a_{j,k}e^{i\alpha_{j,k}}$, where $a_{j,k}\in\mathbb{R}$ and $\alpha_{j,k}\in[-\pi,\pi)$. Further denote $F_2^H\bm y=\bm q\odot e^{i\bm \varphi}$, where $\bm q\in\mathbb{R}_+^{n}$ and $\bm \varphi\in[-\pi,\pi)^{n}$, with
    \begin{gather*}
        u_i=\frac{2q_i}{\sigma^2}\cos\varphi_i-\sum\limits_{\substack{k=1 \\ k\neq i}}^n\frac{2a_{i,k}}{\sigma^2}\cos(\tilde{\alpha}_{i,k}+\phi_k),\quad\quad
        v_i=\frac{2q_i}{\sigma^2}\sin\varphi_i-\sum\limits_{\substack{k=1 \\ k\neq i}}^n\frac{2a_{i,k}}{\sigma^2}\sin(\tilde{\alpha}_{i,k}+\phi_k),
    \end{gather*}
   and  $\tilde{\alpha}_{i,k}=\sgn(k-i)\alpha_{i,k}$ for $k=1,\dots,n$. Then
\begin{align*}
    f_{\Phi_i|\Phi_{-i},\mathcal{Y},\mathcal{G},\mathcal{T}^2}(\phi_i|\phi_{-i},\bm y,\bm g,\bm\tau^2)\propto f_{vM}\left(\phi_{i}\Big|\mu_i,\kappa_i\right),\quad i=1,\dots,n,
\end{align*}
where $f_{vM}(x|\mu,\kappa)$ is the von Mises probability density function, given as
\begin{align}\label{eq:vonmises_def}
    f_{vM}(x|\mu,\kappa)=\frac{\exp(\kappa\cos(x-\mu))}{2\pi I_0(\kappa)},
\end{align}
with location $\mu$ and concentration $\kappa$. Here
\begin{gather}                 \kappa_i=\sqrt{u_i^2+v_i^2},\quad\quad\mu_i=\begin{cases}\arctan\left(-\frac{v_i}{u_i}\right) & \text{if } u_i>0 \\ \pi/2 & \text{if } u_i=0 \\ \arctan\left(-\frac{v_i}{u_i}\right)+\pi & \text{if } u_i<0. 
    \end{cases}
    \label{eq:kappamu}
\end{gather}
Note that in \eqref{eq:vonmises_def}, $I_0$ is the  zeroth order modified Bessel function of the first kind.
\end{theorem}
The proof to \cref{thm:cond_vonmises} is given in \cref{sec:appendix_proofs}.

Since the calculation of $\{\kappa_i\}_{i = 1}^n$ and $\{\mu_i\}_{i = 1}^n$ in \cref{eq:kappamu}  require only scalar operations, our method does not suffer from the curse of dimensionality as $n$ grows large. Furthermore, when $F$ is unitary, $f_{\Phi|\mathcal{Y},\mathcal{G},\mathcal{T}^2}(\bm \phi|\bm y,\bm g,\bm\tau^2)$ is precisely a product of von Mises density functions, as is told in \cref{thm:vonmises}.

\begin{theorem}\label{thm:vonmises}
Suppose $F_2=FD(\bm g)$, where $F$ is unitary. Let $F_2^H\bm y=\bm q\odot e^{i\bm \varphi}$ where $\bm q\in\mathbb{R}_+^{n}$ and $\bm \varphi\in[-\pi,\pi)^{n}$. Then 
\begin{align*}
    f_{\Phi|\mathcal{Y},\mathcal{G}}(\bm \phi|\bm y,\bm g)=\prod_{i=1}^nf_{vM}\left(\phi_{i}\Big|\varphi_{i},\frac{2}{\sigma^2}q_{i}\right).
\end{align*}
\end{theorem}
The proof to \cref{thm:vonmises} is provided in \cref{sec:appendix_proofs}.

When $F$ is unitary, the conditional independence shown in \cref{thm:vonmises} allows us to update each $\Phi_j$ independently of $\Phi_{-j}$, increasing the opportunity for parallelization in our technique. Even when $F$ is not unitary, \cref{thm:cond_vonmises} allows us to update each $\Phi_i$ sequentially using a Gibbs sampling scheme, although we do not benefit from the same parallelization opportunities as when $F$ is unitary.

To sample the von Mises distribution, we use the wrapped Cauchy distribution, which has probability density function
\begin{align*}
	f_{WC}(\theta;\mu,\gamma)=\sum_{k=-\infty}^\infty\frac{\gamma}{\pi(\gamma^2+(\theta-\mu+2\pi k)^2)}, \quad -\pi<\theta<\pi.
\end{align*}
Here $\gamma$ is the scale factor and $\mu$ is the mode of the unwrapped distribution. The acceptance-rejection method introduced in \cite{best1979efficient} utilizes a wrapped Cauchy density as an envelope for sampling from the von Mises distribution \eqref{eq:vonmises_def}.

\begin{algorithm}
	\caption{Sampling the von Mises density using the wrapped Cauchy}
	\label{alg:rejAccCauchy}
	\begin{algorithmic}
		\item[] {\bf Input} Location $\mu$ and concentration $\kappa$ of von Mises distribution.
        \item[] {\bf Output} Sample $\theta$.
		\item[1] Set $\tau=1+(1+4\kappa^2)^\frac12$, $\rho=(\tau-(2\tau)^\frac12)/(2\kappa)$, and $r=(1+\rho^2)/(2\rho)$.
		\item[2] Generate $u_1\sim U(0,1)$, then set $z=\cos(\pi u_1)$, $f=(1+rz)/(r+z)$, $c=\kappa(r-f)$.
		\item[3] Generate $u_2\sim U(0,1)$, then if $c(2-c)-u_2>0$, go to step $5$.
		\item[4] If $\ln(c/u_2)+1-c<0$, return to step $2$.
		\item[5] Generate $u_3\sim U(0,1)$, then set $\theta=[\text{sgn}(u_3-0.5)]\cos^{-1}(f)$.
	\end{algorithmic}
\end{algorithm}

\begin{algorithm}[H]
	\caption{CVBL for a signal with sparsity in the transformed magnitude}
	\label{alg:propMethod}
	\begin{algorithmic}
		\item[] {\bf Input} data $\bm y$, noise variance $\sigma^2$, hyperparameters $r=1$ and $\delta=10^{-3}$, sparse transform operator $L$, chain length $N_M$, and burn-in length $B$.
        \item[] {\bf Output} samples $\bm z^{(s-B+1)}=\bm g^{(s)}\odot\exp(i\bm \phi^{(s)})$ for $s=B,\dots,N_M$.
		\item[1] Set $\bm g^{(0)}=|A^H\bm y|$, $(\bm\tau^2)^{(0)}=1$, and $\bm \phi^{(0)}=\arg(A^H\bm y)$
		\item[2] {\bf For} $\ell=1,\dots,N_M$ {\bf do}
		\begin{itemize}
		    \item[i.] Draw a sample $\bm g^*$ from $\tilde{f}_{\mathcal{G}|\mathcal{Y},\Phi,\mathcal{T}^2}\left(\bm g\big|\bar{\bm y},\bm \phi^{(\ell-1)},\left(\bm\tau^2\right)^{(\ell-1)}\right)$ using Theorem \ref{thm:magpost}.
		    \item[ii.] If $\bm g^*$ contains a negative element, return to (i). Otherwise, set $\bm g^{(\ell)}=\bm g^*$.
                \item[iii.] Sample each $\left(\tau_{j}^{-2}\right)^{(\ell)}$ from $\mathcal{T}_{j}^{-2}|\mathcal{G}=\bm g^{(\ell)},\upeta^{-2}=\left(\eta^{-2}\right)^{(\ell-1)}$ \eqref{eq:tau_posterior_complex} for  $j=1,\dots,k$.
                \item[iv.] Sample $\left(\eta^{-2}\right)^{(\ell)}$ from $\mathcal{T}_{j}^{-2}=\left(\tau^2\right)^{(\ell)}$ \eqref{eqn:hier_hypereta}.
		    \item[v.] For each $i=1,\dots,n$, draw a sample $\phi_i^*$ from $\Phi_i|\mathcal{Y}=\bm y, \mathcal{G}=\bm g^{(\ell)}$ \eqref{eq:edge_postphi} using Algorithm \ref{alg:rejAccCauchy} and set $\phi_i^{(\ell)}=\phi_i^*$.
		\end{itemize}
	\end{algorithmic}
\end{algorithm}

We are now ready to sample from the joint distribution $f_{\mathcal{G},\Phi|\mathcal{Y},\bm\eta}(\bm g,\bm \phi|\bm y,\eta)$ in \eqref{eq:posterior_edgeL1}. After initializing our chain, the three-stage Gibbs sampler is implemented, where $\bm g$ is updated first, followed by the $\bm\tau^2$ update, and concluded by the $\bm \phi$ update. This is done for some predetermined number of iterations $N_M$, after which the output chain is formed using all the samples generated after the burn-in period $B$. This method is summarized in \cref{alg:propMethod}.\footnote{\
Both \cref{alg:gibbs_signalL1} and \cref{alg:propMethod} may be easily modified to use hyperprior \eqref{eq:delta_hyper} instead of \eqref{eq:hyperprior} by fixing the value of $\left(\eta^{-2}\right)^{(\ell)}=\hat{\eta}^{-2}$ for all $\ell=1,\dots,N_M$.}

\section{Numerical Results} \label{sec:numerics} 

We now demonstrate the efficacy of the CVBL algorithm given in \cref{alg:gibbs_signalL1} for the sparse magnitude case and \cref{alg:propMethod} for the sparse transform of the magnitude case by performing experiments in 1D and 2D using three different forward operators for $F$ in \eqref{eqn:observation_model}:
\begin{itemize}
    \item $F_F\in\mathbb{C}^{n\cross n}$: The discrete Fourier transform operator  with entries (in 1D) given by 
    \begin{align}
\label{eq:Fouriertransform}
    [F_F]_{j,k} = \frac{1}{\sqrt{n}}\exp(-2\pi i\frac{jk}{n}),\quad\quad j,k=0,\dots,n-1.
\end{align}
Observe that $F_F$ is unitary, so \Cref{thm:vonmises} applies.
\item $F_B\in\mathbb{R}^{n\cross n}$: A blurring operator  with entries (in 1D) given by 
\begin{align}
    [F_B]_{j,k} = \begin{cases}
        \frac{1}{\sqrt{26}}(2^{2-|j-k|}) & \text{ if } |j-k| \leq 2\\
       0 & \text{ else}
    \end{cases}, \quad\quad j,k=0,\dots,n-1. \label{eq:blurtransform}
\end{align}
Observe that $F_B$ is a banded Toeplitz matrix, and although non-singular, it is notoriously ill-conditioned.
\item $F_U\in\mathbb{C}^{m\cross n}$: A random under-sampled Fourier transform matrix which has entries (in 1D) given by \eqref{eq:Fouriertransform} but with randomly zeroed-out rows so that $m = \lceil \nu n\rceil$, where $0< \nu \leq 1$. For $\mathcal{M}_\nu\subseteq\{2,3,\dots,n\}$ such that $|\mathcal{M}_\nu|=m$, $F_U$ is defined as
\begin{align}
    [F_U]_{j,k} = \begin{cases}
        [F_F]_{j,k} & \text{ if } j \in \mathcal{M}_\nu\\
       0 & \text{ else}
    \end{cases}, \quad\quad j,k=1,\dots,n. \label{eq:undertransform}
\end{align}
Note that the zeroth frequency term is never zeroed-out in the construction of $\mathcal{M}_\nu$.
\end{itemize} 

\begin{remark}\label{rem:blur}
The extension to 2D for $F_F$ and $F_U$ is straightforward. For $F_B$, the 2D operator is defined such that $F_B\bm z$ convolves the image $\bm z$ with the kernel $K$ given by

\begin{align*}
    K = \frac{1}{2\sqrt{70}}\begin{bmatrix}
        0 & 0 & 1 & 0 & 0\\
        0 & 1 & 2 & 1 & 0\\
        1 & 2 & 16 & 2 & 1\\
        0 & 1 & 2 & 1 & 0\\
        0 & 0 & 1 & 0 & 0
    \end{bmatrix},
\end{align*}
causing the blurring effect to occur in both dimensions.\footnote{The 1D and 2D forward operators used in our experiments are explicitly stated for reproducibility purposes.}
\end{remark}

In the numerical examples that follow, $L=\mathbb{I}_n$ in the sparse signal case or the first order differencing operator in the case of the magnitude having a sparse gradient, where we enforce zero boundary conditions to ensure that it is of rank $n$. For 1D signals this amounts to
\begin{align}
\label{eq:Ldiff}
    [L]_{i,j}=\begin{cases}
        1&\text{if } i=j,\\
        -1&\text{if } i=j+1,\\
        0&\text{else}.
    \end{cases}
\end{align}

Finally we assume that the variance $\sigma^2$ in our observation model \eqref{eqn:observation_model} is known, with corresponding signal-to-noise (SNR) given by 
\begin{align}
    \text{SNR}=10\log_{10}\left(\frac{\norm{F\bm z^{exact}}_2^2}{m\sigma^2}\right), \label{eq:SNR}
\end{align}
where $\bm z^{exact} \in \mathbb{C}^m$ is the exact solution to \eqref{eqn:observation_model}.  The SNR values chosen in our experiments highlight the effectiveness of the CVBL method in noisy environments. 

In all experiments, a total of $5000$ samples are drawn with a burn-in period of $B = 200$. To encourage the sampler to move towards high-mass areas of the probability density, for $\ell < B$ we set $\bm g^{(\ell)}=|\bm g^*|$, regardless of whether or not $\bm g^*$ contains negative elements. Our 1D analysis includes figures showing the means and $90\%$ credibility intervals (CI) for the marginal magnitude of the posterior distribution. We also approximate the marginal density function of the phase of random individual pixels using a kernel density estimation technique \cite{botev2010kernel}. In 2D we provide the magnitude means and the size of the $90\%$ CIs. 
Finally, for each choice of $F$ we also compare the MAP estimate of our CVBL posterior to the corresponding classical LASSO solution $\bm z^{LASSO}$: 
\begin{itemize}
    \item {\bf sparse signal case:}
We compute the solution to the objective function
\begin{align}\label{eq:sparsesig_classical_lasso}
    \bm z^{LASSO}=\argmin_{\bm z}\norm{F\bm z - \bm y}_2^2 + \lambda\norm{\bm z}_1
\end{align}
by considering the real and imaginary parts of $\bm z$, namely $\bm a$ and $\bm b$, respectively, giving the objective function as
\begin{align}\label{eq:sparsesig_classical_lasso_ab}
    \begin{bmatrix}
        \bm a^{LASSO} \\ \bm b^{LASSO}
    \end{bmatrix}=\argmin_{\bm a,\bm b}\left(\norm{\begin{bmatrix}
        F_R & -F_I \\ F_I & F_R
    \end{bmatrix}\begin{bmatrix}
        \bm a \\ \bm b
    \end{bmatrix} - \begin{bmatrix}
        \bm y_R \\ \bm y_I
    \end{bmatrix}}_2^2 + \lambda\sum_{j=1}^n\sqrt{a_j+b_j}\right),
\end{align}
which we solve using the alternating direction method of multipliers (ADMM) \cite{boyd2011distributed}.

\item {\bf  sparse transform case:}
We again employ ADMM  to compute the solution in the sparse transform case to the generalized LASSO problem. Here, however, since $|\bm z|$ is not differentiable, in order to solve
\begin{align}\label{eq:sparseedge_classical_lasso}
    \bm z^{LASSO}=\argmin_{\bm z}\norm{F\bm z - \bm y}_2^2 + \lambda\norm{L|\bm z|}_1,
\end{align}
we follow what was done in \cite{churchill2023sub,sanders2017composite} and instead use the diagonal matrix $\Theta^{(k)}$ with non-zero entries $\Theta^{(k)}_{jj}=z^{(k)}_j/|z^{(k)}_j|$
at each iteration $k$ of the algorithm. The objective function then becomes
\begin{align*}
    \bm z^{LASSO}=\argmin_{\bm z}\norm{F\bm z - \bm y}_2^2 + \lambda\norm{L\left(\Theta^{(k)}\right)^H\bm z}_1.
\end{align*}

\end{itemize}
In both cases we test $\lambda = \alpha \sigma^2/\bar{\eta}$, where $\bar{\eta}$ is the mean of the samples of $\upeta$ generated by the CVBL method. In some sense  $\alpha = 1$ represents the ``best case'' scenario for selecting suitable parameters, while other values of $\alpha$ allow us to test for robustness. We use $\alpha = 1$ in each experiment unless otherwise specified. 

\subsection{Numerical Efficiency}
\label{subsec:numerics_efficiency}
Sampling Gaussians using  matrix (e.g.~Cholesky) factorization  can be prohibitively expensive in high dimensions. This is because new matrices must be factorized for each re-sampling of the hyper-parameters $\tau^2$, yielding a general cost $\mathcal{O}(n^3)$ flops per sample. The approach detailed in \Cref{alg:pert_bef_opt} provides an efficient way to generate samples of multivariate Gaussian distributions \cite{papandreou2010gaussian,vono2022high}, which we use to sample $\mathcal{A}|\mathcal{Y},\mathcal{B},\mathcal{T}^2$ and $\mathcal{B}|\mathcal{Y},\mathcal{A},\mathcal{T}^2$ in \cref{alg:gibbs_signalL1} and $\mathcal{G}|\mathcal{Y},\Phi,\mathcal{T}^2$ in \cref{alg:propMethod}. Moreover, we can efficiently solve the system in Step $3$ of \cref{alg:pert_bef_opt} using the conjugate gradient method \cite{hestenes1952methods}.

\begin{algorithm}[h!]
	\caption{Gaussian sampling by perturbation-optimization }
	\label{alg:pert_bef_opt}
	\begin{algorithmic}
		\item[] {\bf Input} Mean $\bm\mu_p\in\mathbb{R}^n$, measurement $\bm y\in\mathbb{R}^m$, variances $\Sigma_p$, $\Sigma_\ell$, sparsifying transform $L\in\mathbb{R}^{k\cross n}$, and forward model $F$ corresponding to Gaussian prior $f_\mathcal{X}=\mathcal{N}(L\bm\mu_p,\Sigma_p)$ and likelihood $f_{\mathcal{Y}|\mathcal{X}}=\mathcal{N}(F\bm x,\Sigma_\ell)$.
        \item[] {\bf Output} Sample $\bm x_s$ of the posterior $f_{\mathcal{X}|\mathcal{Y}}$.
		\item[1] Perturb the prior mean $\tilde{\bm\mu}_p\sim\mathcal{N}(\bm\mu_p,\Sigma_p)$.
		\item[2] Perturb the data $\tilde{\bm y}\sim\mathcal{N}(\bm y,\Sigma_\ell)$.
		\item[3] Solve $(L^T\Sigma_p^{-1}L+F^T\Sigma_\ell^{-1} F)\bm x_s=(L^T\Sigma_p^{-1}\tilde{\bm\mu}_p+F^T\Sigma_\ell^{-1}\tilde{\bm y})$ with respect to $\bm x_s$.
	\end{algorithmic}
\end{algorithm}

Finally we point out that when the forward operator is given by $F_B$ \eqref{eq:blurtransform} and the magnitude sparsity is in the gradient domain, we use a block sampling approach to generate samples of the phase \eqref{eq:edge_postphi},  allowing us to  takes advantage of the inherent sparsity in $F_B$ and to update portions of the phase in parallel.

\subsection{Sparsity in Magnitude}
\label{subsec:numerics_sparsemag}
\begin{figure}[h!]

    \begin{subcaptionblock}{.28\textwidth}

    \includegraphics[width=.85\linewidth]{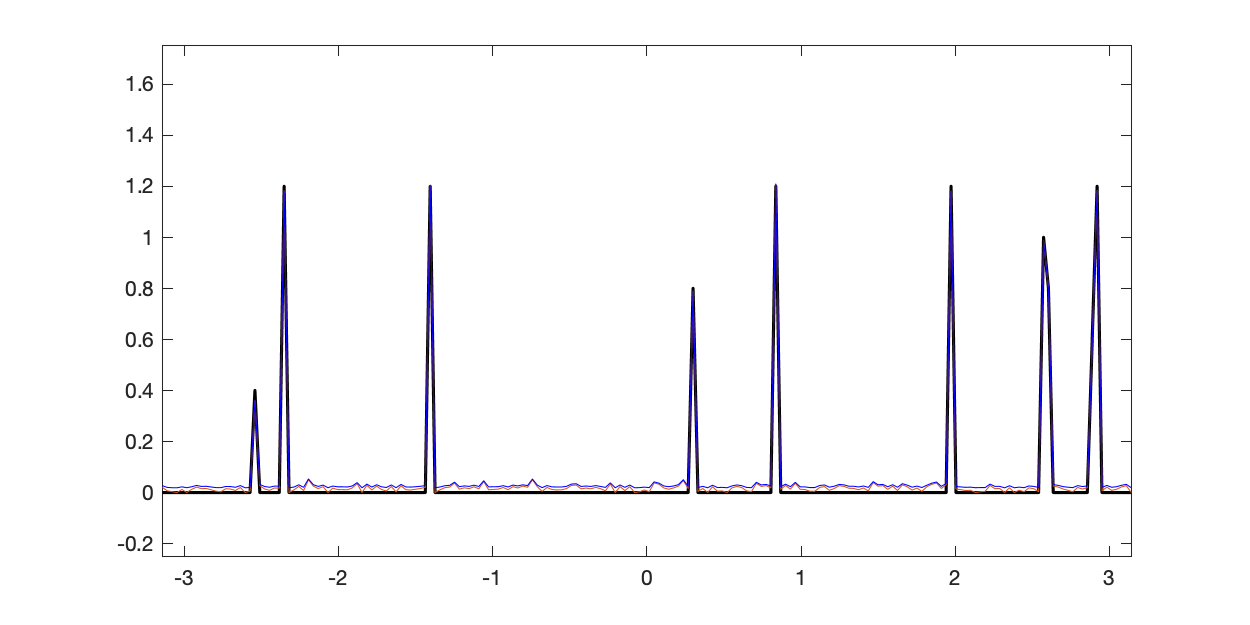}
    \caption{$F_F$ in \eqref{eq:Fouriertransform}}
    \end{subcaptionblock}
    \begin{subcaptionblock}{.28\textwidth}

    \includegraphics[width=.85\linewidth]{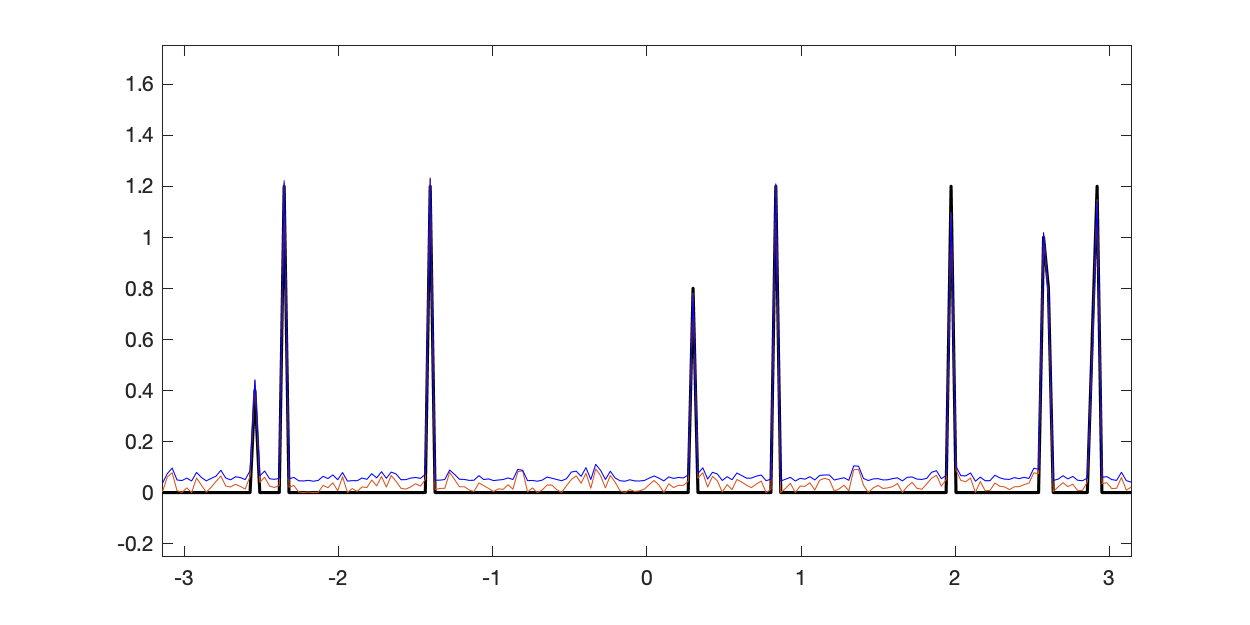}
    \caption{$F_B$ in \eqref{eq:blurtransform}}
    \end{subcaptionblock}
    \begin{subcaptionblock}{.28\textwidth}

    \includegraphics[width=.85\linewidth]{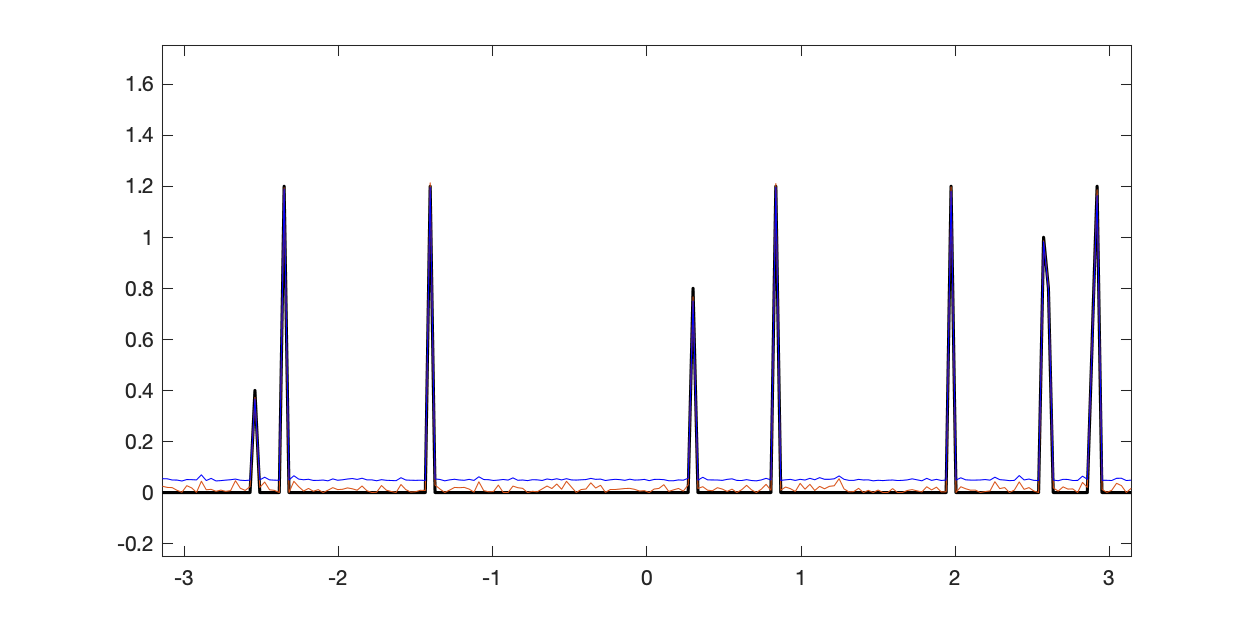}
    \caption{$F_U$ in \eqref{eq:undertransform}; $\nu=0.8$}
    \end{subcaptionblock}
    \begin{subcaptionblock}[b][2.3cm][t]{.13\textwidth}

    \includegraphics[width=.85\linewidth]{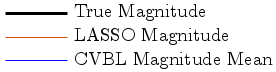}
    \end{subcaptionblock}
    \begin{subcaptionblock}{.28\textwidth}

    \includegraphics[width=.85\linewidth]{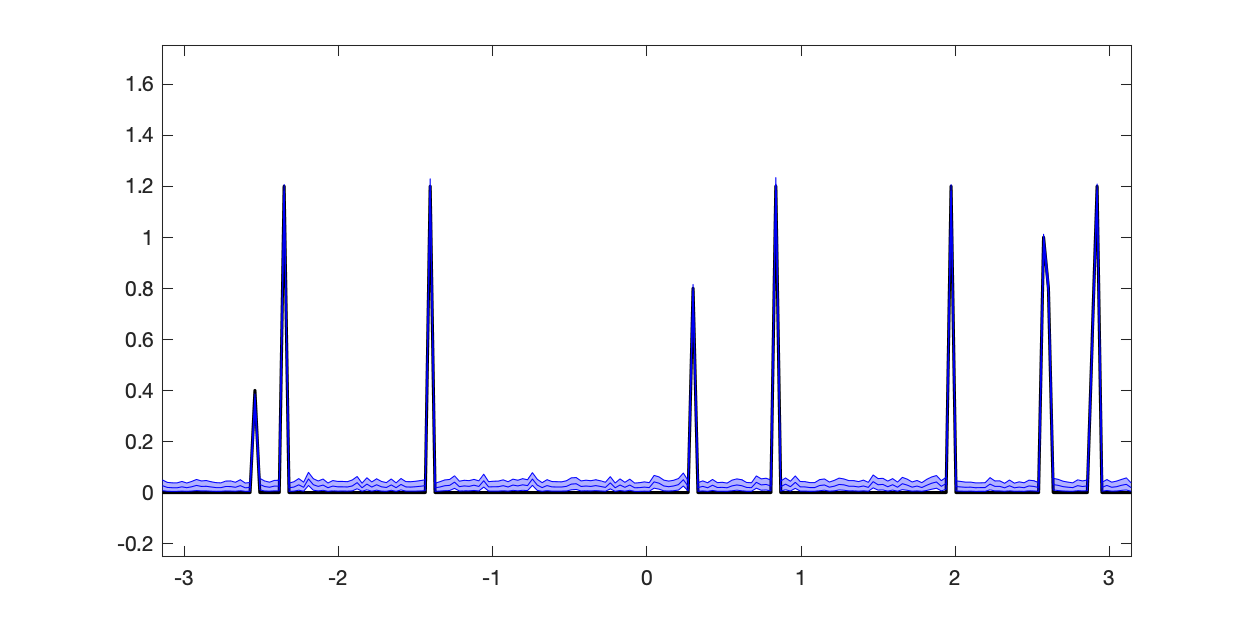}
    \caption{$F_F$ in \eqref{eq:Fouriertransform}}
     \end{subcaptionblock}
    \begin{subcaptionblock}{.28\textwidth}

    \includegraphics[width=.85\linewidth]{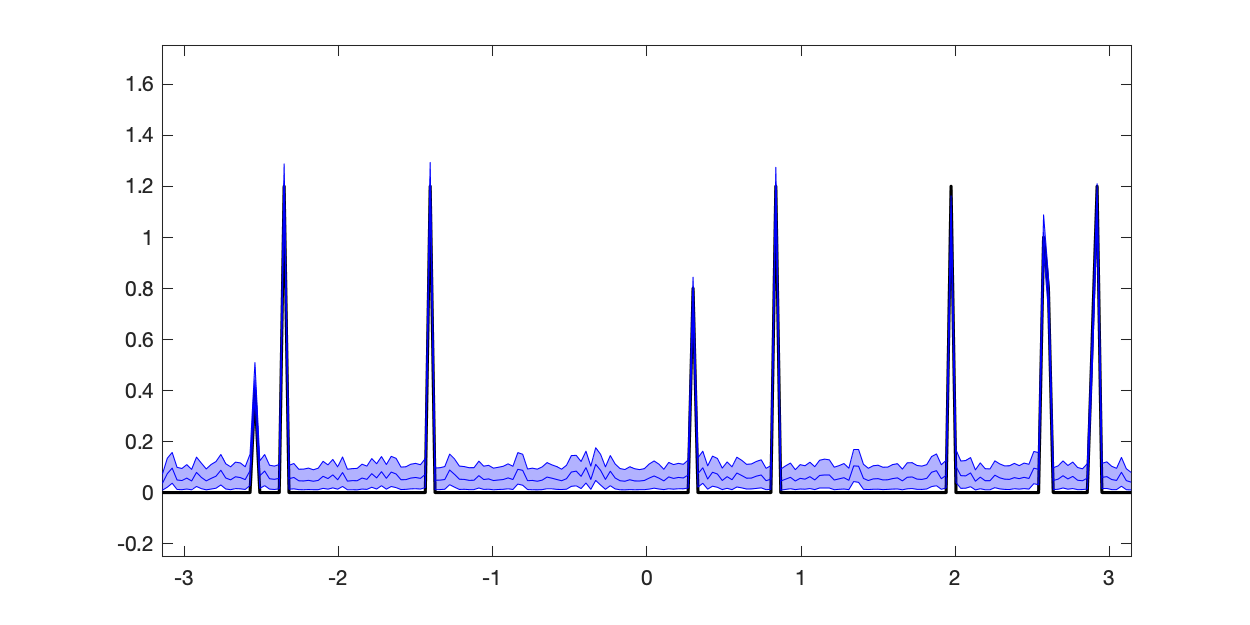}
    \caption{$F_B$ in \eqref{eq:blurtransform}}
     \end{subcaptionblock}
    \begin{subcaptionblock}{.28\textwidth}

    \includegraphics[width=.85\linewidth]{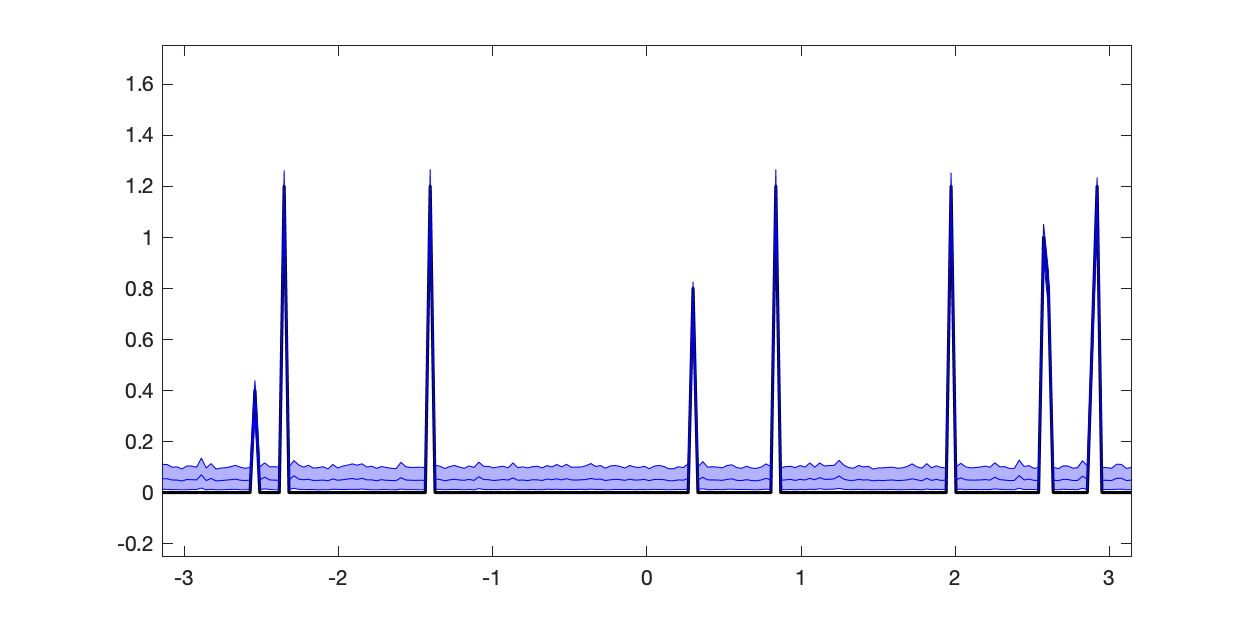}
    \caption{$F_U$ in \eqref{eq:undertransform}; $\nu=0.8$}
     \end{subcaptionblock}
     \begin{subcaptionblock}[b][2.1cm][t]{.13\textwidth}

    \includegraphics[width=.85\linewidth]{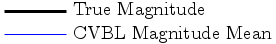}
     \end{subcaptionblock}
    \caption{Recovered magnitude values for complex-valued signals with sparse magnitude.  (a)-(c) The magnitude of the LASSO solution and the mean of the CVBL output; (d)-(f) $90\%$ CI for the magnitude reconstruction. Here SNR $= 20\dB$.}
    \label{fig:1D_sparsemag_mag}
\end{figure}

For the sparse magnitude case, we perform experiments with \cref{alg:gibbs_signalL1} using $5000$ samples and the forward operators in \eqref{eq:Fouriertransform}, \eqref{eq:blurtransform}, and \eqref{eq:undertransform}.  \Cref{fig:1D_sparsemag_mag} shows the magnitude means and $90\%$ CIs for the 1D noisy experiments for $\text{SNR}  = 20\dB$, \eqref{eq:SNR}. 

\begin{figure}[h!]

    \begin{subcaptionblock}{.28\textwidth}

    \includegraphics[width=.85\linewidth]{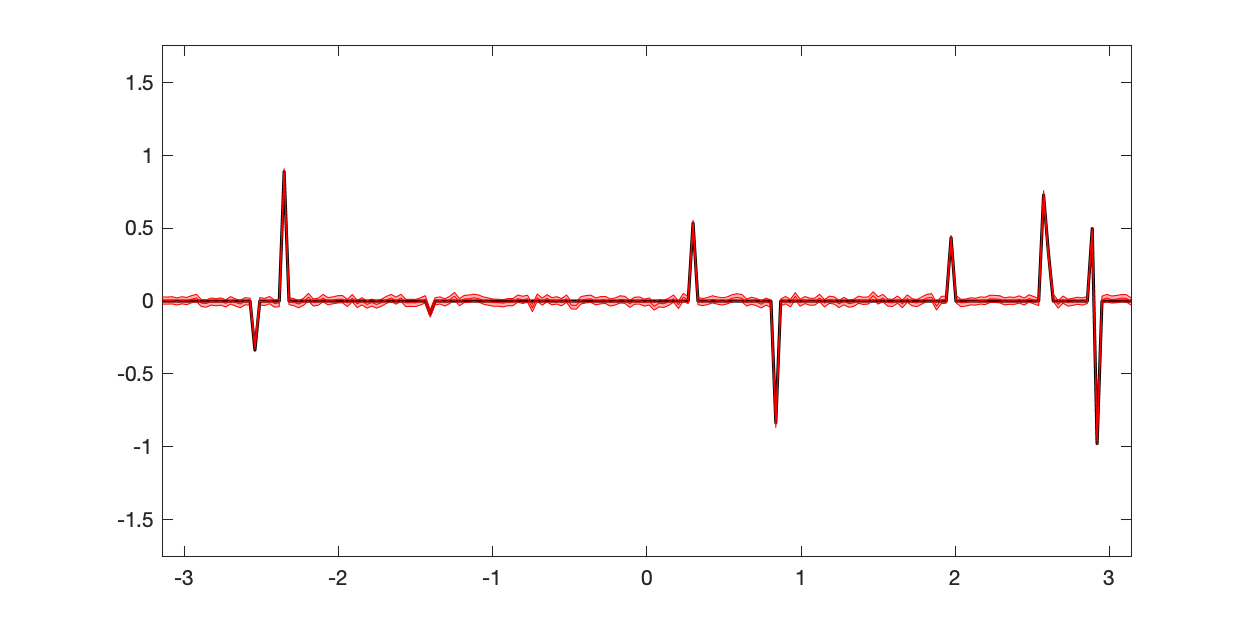}
    \caption{$F_F$ in \eqref{eq:Fouriertransform}}
    \end{subcaptionblock}
    \begin{subcaptionblock}{.28\textwidth}

    \includegraphics[width=.85\linewidth]{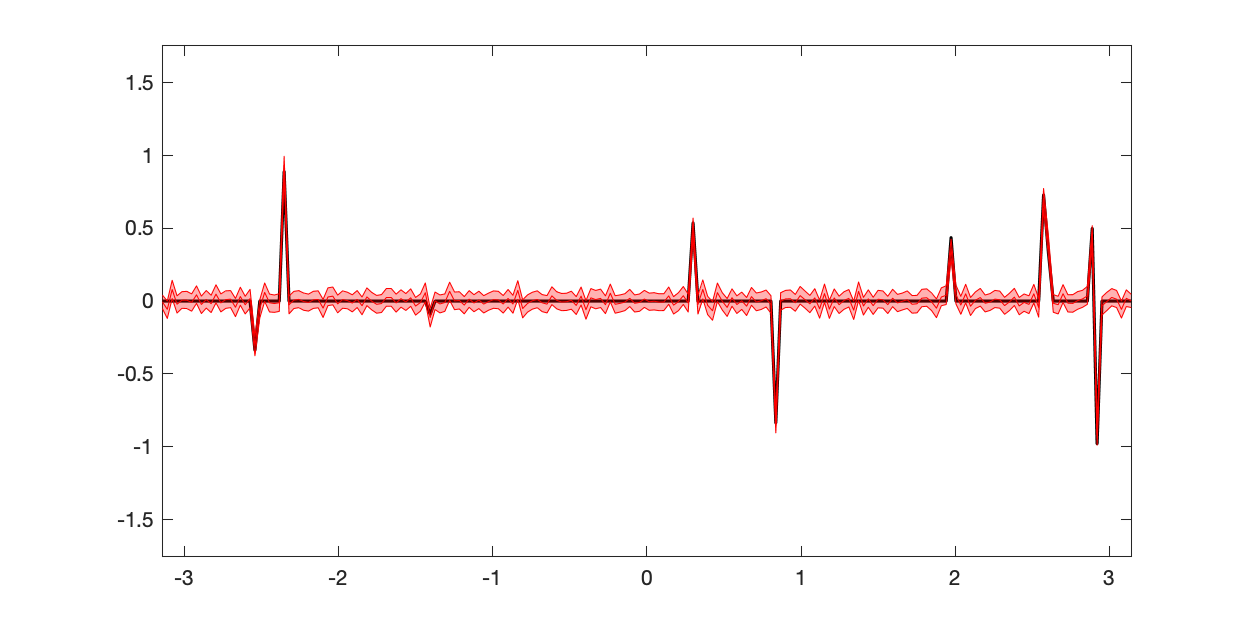}
    \caption{$F_B$ in \eqref{eq:blurtransform}}
    \end{subcaptionblock}
    \begin{subcaptionblock}{.28\textwidth}

    \includegraphics[width=.85\linewidth]{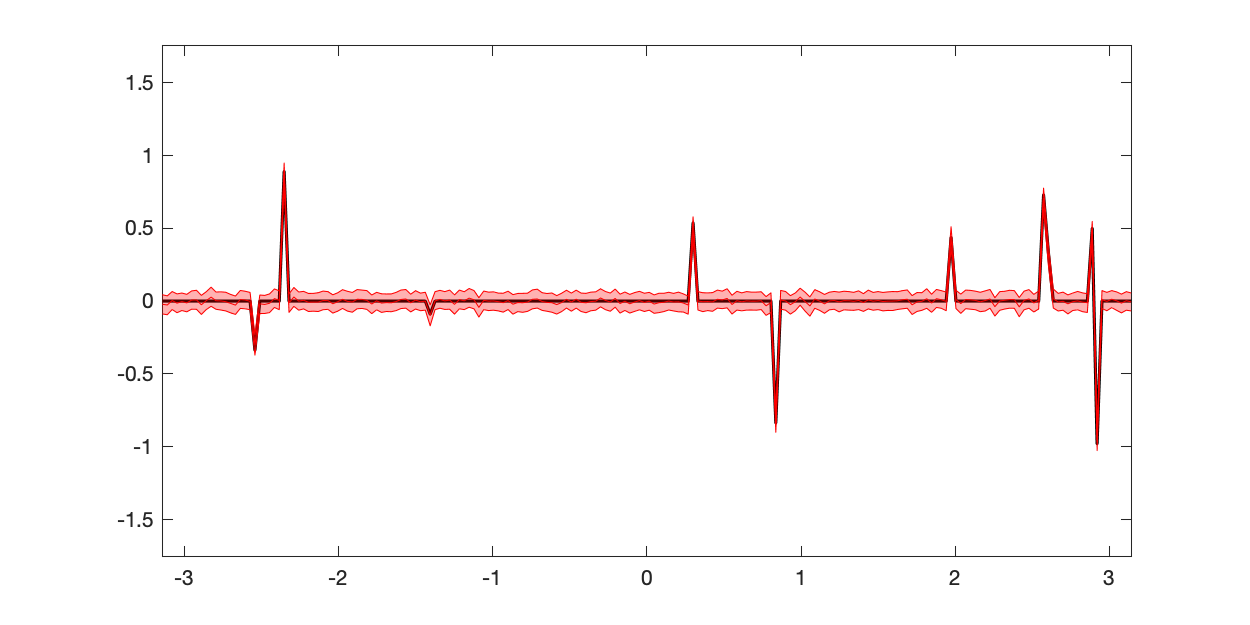}
    \caption{$F_U$ in \eqref{eq:undertransform}; $\nu=0.8$}
    \end{subcaptionblock}
    \begin{subcaptionblock}[b][2.3cm][t]{.11\textwidth}

    \includegraphics[width=.85\linewidth]{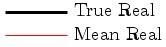}
     \end{subcaptionblock}
     \begin{subcaptionblock}{.28\textwidth}

    \includegraphics[width=.85\linewidth]{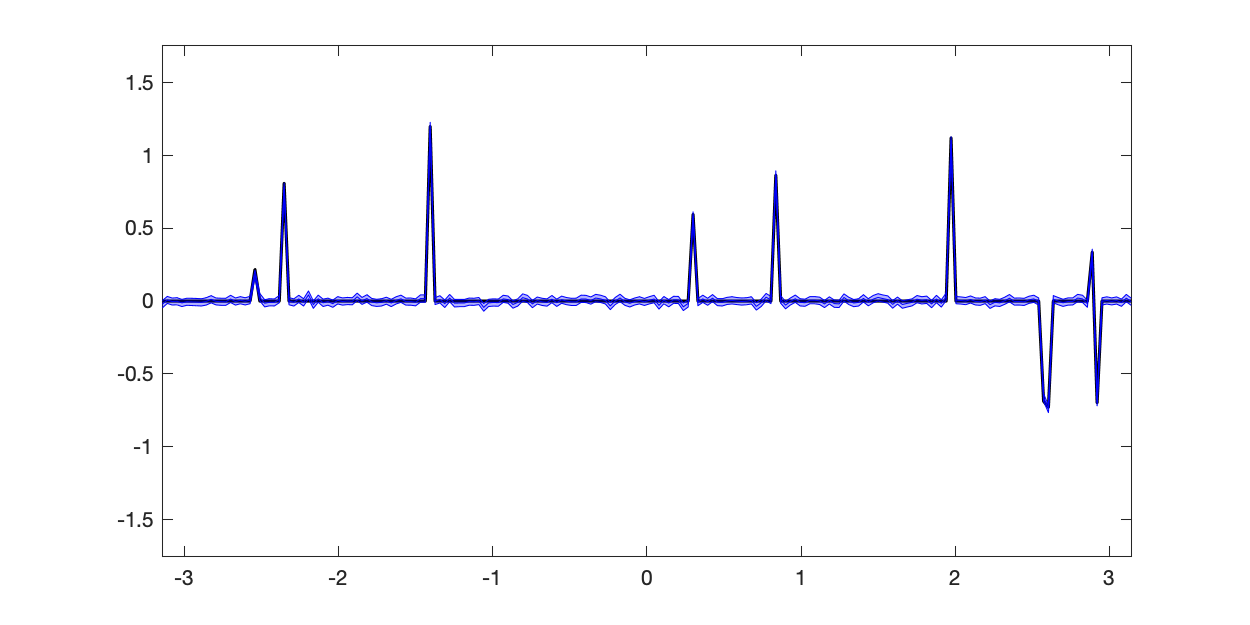}
    \caption{$F_F$ in \eqref{eq:Fouriertransform}}
    \end{subcaptionblock}
    \begin{subcaptionblock}{.28\textwidth}

    \includegraphics[width=.85\linewidth]{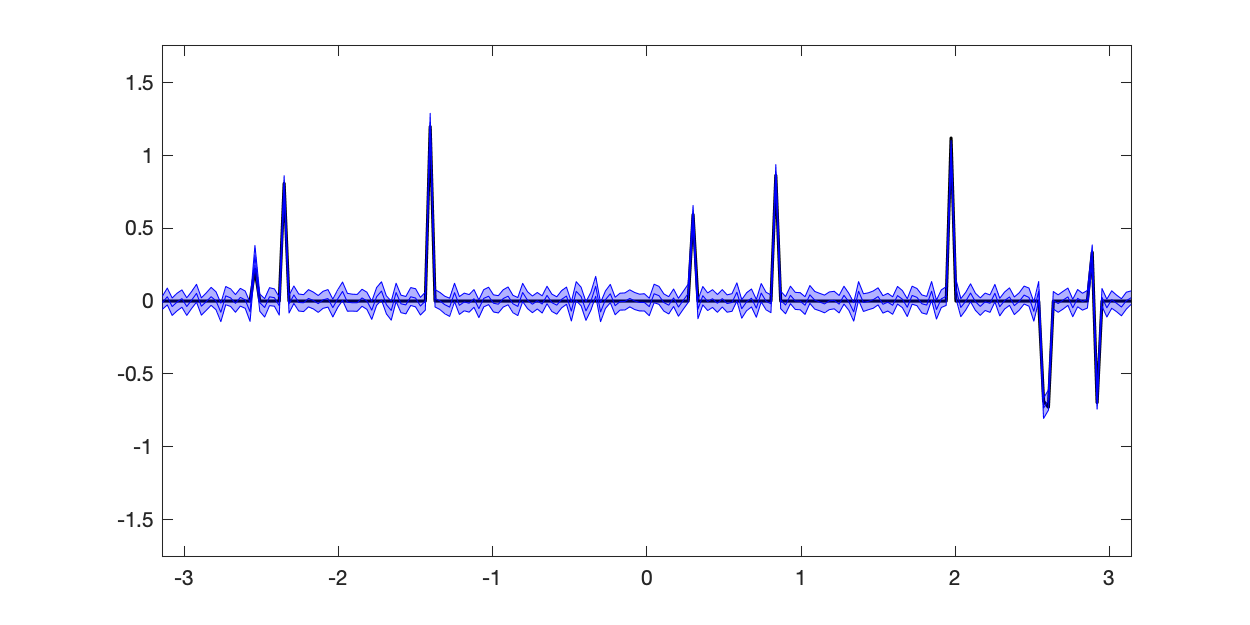}
    \caption{$F_B$ in \eqref{eq:blurtransform}}
    \end{subcaptionblock}
    \begin{subcaptionblock}{.28\textwidth}

    \includegraphics[width=.85\linewidth]{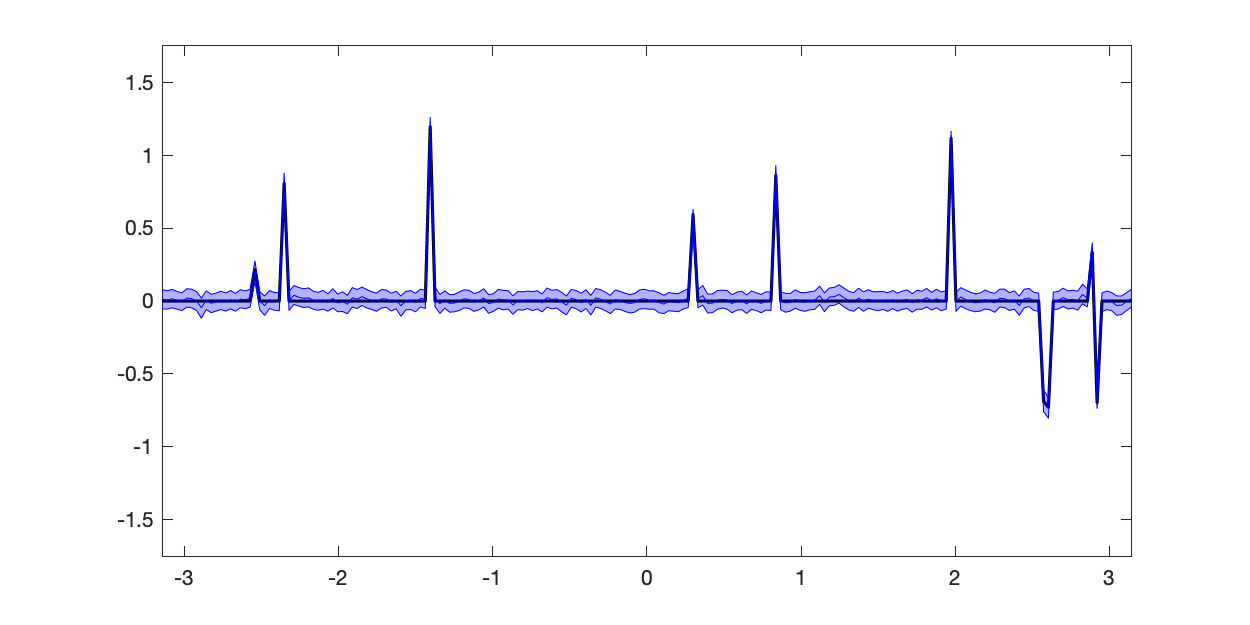}
    \caption{$F_U$ in \eqref{eq:undertransform}; $\nu=0.8$}
    \end{subcaptionblock}
    \begin{subcaptionblock}[b][2.3cm][t]{.13\textwidth}

    \includegraphics[width=.85\linewidth]{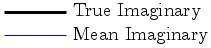}
     \end{subcaptionblock}
    \caption{Means of (top) real and (bottom) imaginary parts of the signal in the sparse magnitude case. Shaded regions indicate $90\%$ CIs for the CVBL recovery. SNR $= 20\dB$.}
    \label{fig:1D_sparsemag_AandB}
\end{figure}

Although the recovered magnitude means are not sparse, we see in \cref{fig:1D_sparsemag_AandB} that the real and imaginary components of the signal are close to zero outside of the signal support. Combined with \cref{fig:1D_sparsemag_mag} our results show that, as expected, the CVBL provides mean information similar to that of the classical LASSO point-estimate technique \eqref{eq:sparsesig_classical_lasso}. 

\begin{figure}
    \centering
    \begin{subcaptionblock}{.21\textwidth}
    \centering
    \includegraphics[width=1\linewidth]{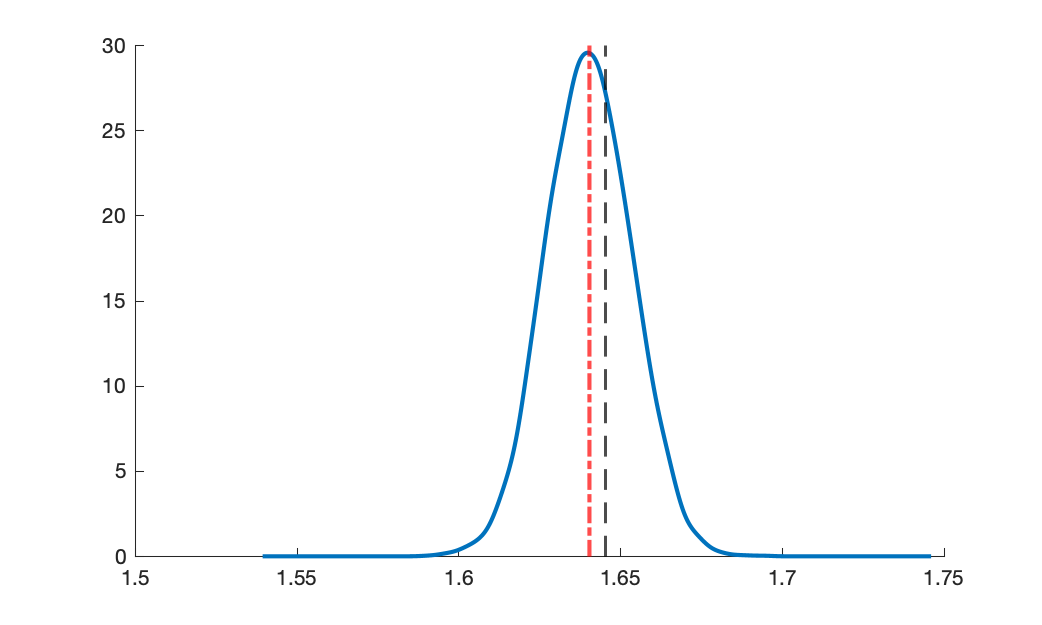}
    \caption{$F_F$ in \eqref{eq:Fouriertransform}}
    \end{subcaptionblock}
    \begin{subcaptionblock}{.21\textwidth}
    \centering
    \includegraphics[width=1\linewidth]{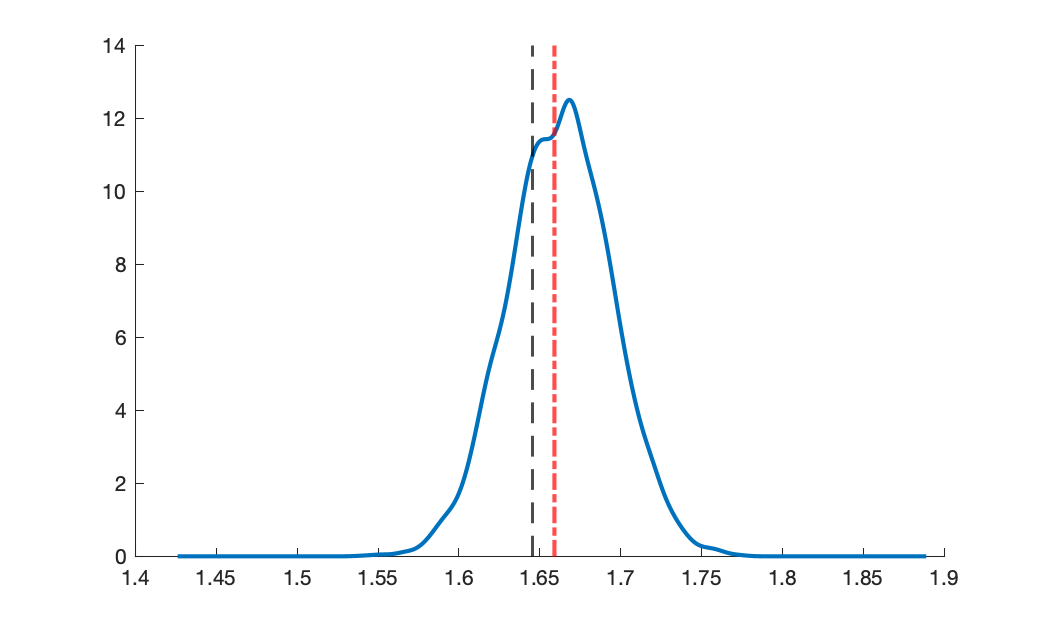}
    \caption{$F_B$ in \eqref{eq:blurtransform}}
    \end{subcaptionblock}
    \begin{subcaptionblock}{.21\textwidth}
    \centering
    \includegraphics[width=1\linewidth]{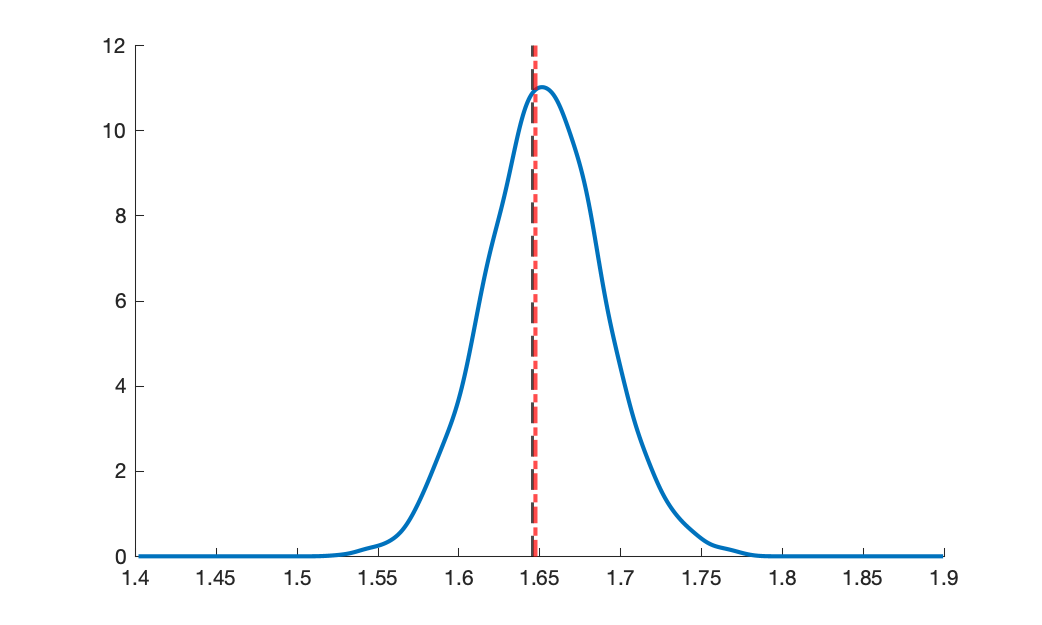}
    \caption{$F_U$ in \eqref{eq:undertransform}; $\nu=0.8$}
    \end{subcaptionblock}
    \begin{subcaptionblock}[b][2.4cm][t]{.13\textwidth}
    \centering
    \includegraphics[width=1\linewidth]{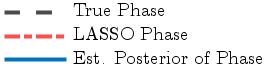}
     \end{subcaptionblock}
    \caption{Phase estimate in the sparse magnitude case at a randomly selected ``on'' pixel (${\bm g}_j = 1$). True phase, LASSO phase point estimate, and CVBL marginal posterior for the phase. SNR = $20 \dB$.}
    \label{fig:1D_sparsemag_phase}
\end{figure}

\Cref{fig:1D_sparsemag_phase} displays the phase results for each forward operator at a randomly selected ``on" pixel (where ${\bm g}_j >0$).\footnote{The probability density plots in \cref{fig:1D_sparsemag_phase}, \cref{fig:1D_sparseedge_phase}, and \cref{fig:2D_transform_phase} are formed using kernel density estimation techniques \cite{botev2010kernel}.} In all three experiments, the CVBL method recovers the support of the signal as well as uncertainty information for both the magnitude and the phase.
 
\subsection{Sparsity in Transform of Magnitude}
\label{subsec:numerics_sarsemagtransform}
We now consider the case where sparsity is expected in some linear transform domain of the magnitude. Here the magnitude $\bm g = \{g_j\}_{j= 1}^{200}$ is given by $\{f(t_j)\}_{j = 1}^{200}$  for the function $f:[-\pi,\pi)\to\mathbb{R}_+$ defined as
\begin{align}\label{eq:1D_sparseedge_testfunct}
    f(t)=\begin{cases}
        2, & -2.8 \leq t \leq -2.1\\
        1.5 & -1.6 \leq t \leq -1.3\\
        1 + \frac32\exp(\left(\frac{t-\pi/2}{2/3}\right)^2) & t>0\\
        1 & \text{else}
    \end{cases}
\end{align}
and $t_j = -\pi + \frac{j\pi}{100}$. The corresponding phase $\Phi = \{\phi_j\}_{j = 1}^{200}$  is randomly chosen uniformly from $[-\pi,\pi)$. Observe that while  $|f|(t)$ is sparse in the gradient domain when $t\leq0$, this is not the case when $t>0$. Since our prior is based on the assumption that the gradient domain is sparse, using \cref{eq:1D_sparseedge_testfunct} allows us to test the effectiveness of the CVBL method  in regions where the gradient is nonzero. As already noted in \Cref{sec:sparseoperator}, other choices of sparsifying transform operators such as HOTV  may be more suitable and our method is not inherently limited to using the identity or TV operators.  For simplicity, as well as to emphasize robustness of our approach, we use the differencing operator in \eqref{eq:Ldiff} and leave other operators for future investigations.

\begin{figure}[h!]

    \begin{subcaptionblock}{.28\textwidth}

    \includegraphics[width=.9\linewidth]{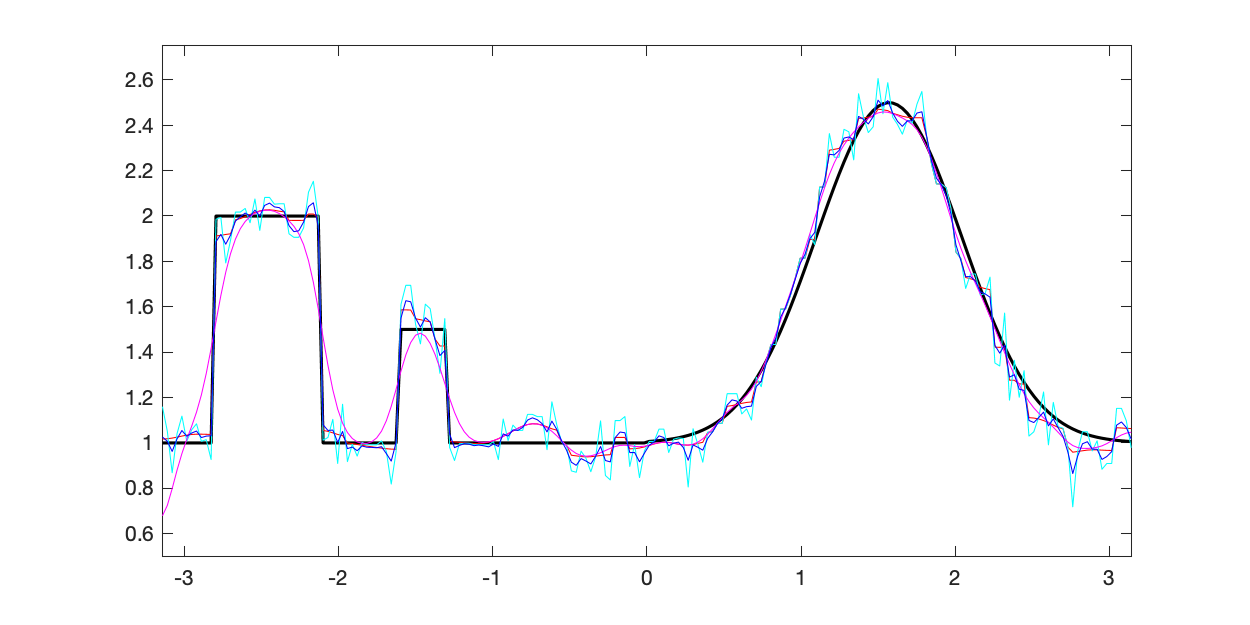}
    \caption{$F_F$ in \eqref{eq:Fouriertransform}}
    \end{subcaptionblock}
    \begin{subcaptionblock}{.28\textwidth}

    \includegraphics[width=.9\linewidth]{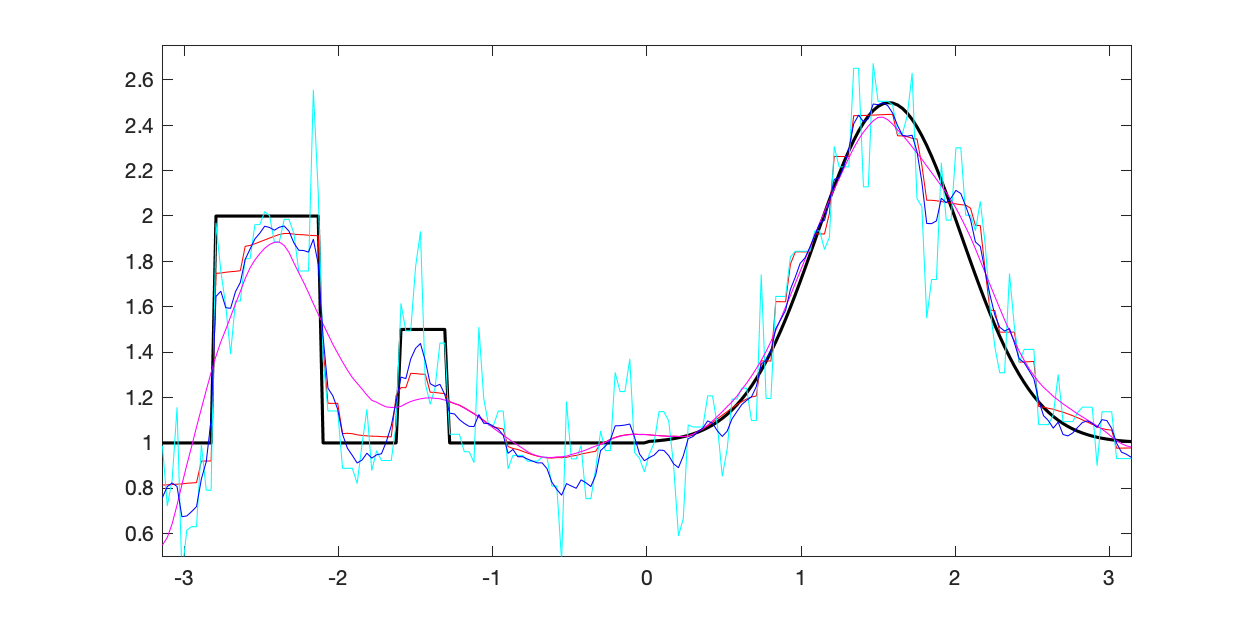}
    \caption{$F_B$ in \eqref{eq:blurtransform}}
    \end{subcaptionblock}
    \begin{subcaptionblock}{.28\textwidth}

    \includegraphics[width=.9\linewidth]{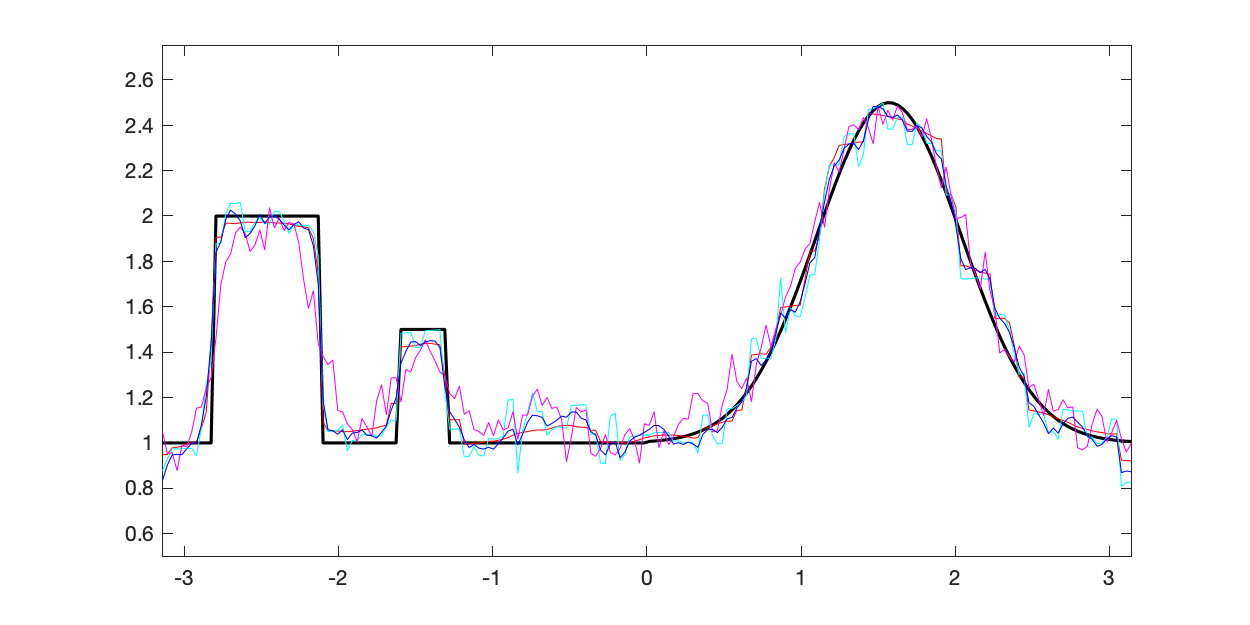}
    \caption{$F_U$ in \eqref{eq:undertransform}; $\nu=0.8$}
    \end{subcaptionblock}
    \begin{subcaptionblock}[b][2.3cm][t]{.13\textwidth}

    \includegraphics[width=.9\linewidth]{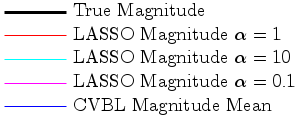}
    \end{subcaptionblock}
    \begin{subcaptionblock}{.28\textwidth}

    \includegraphics[width=.9\linewidth]{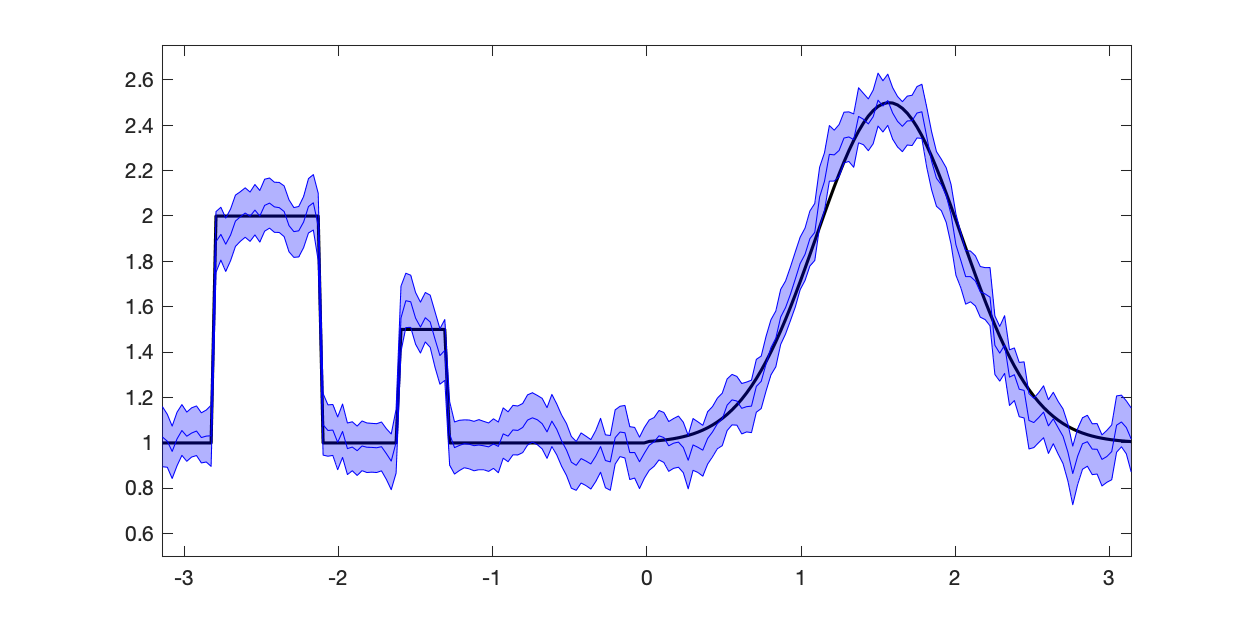}
    \caption{$F_F$ in \eqref{eq:Fouriertransform}}
     \end{subcaptionblock}
    \begin{subcaptionblock}{.28\textwidth}

    \includegraphics[width=.9\linewidth]{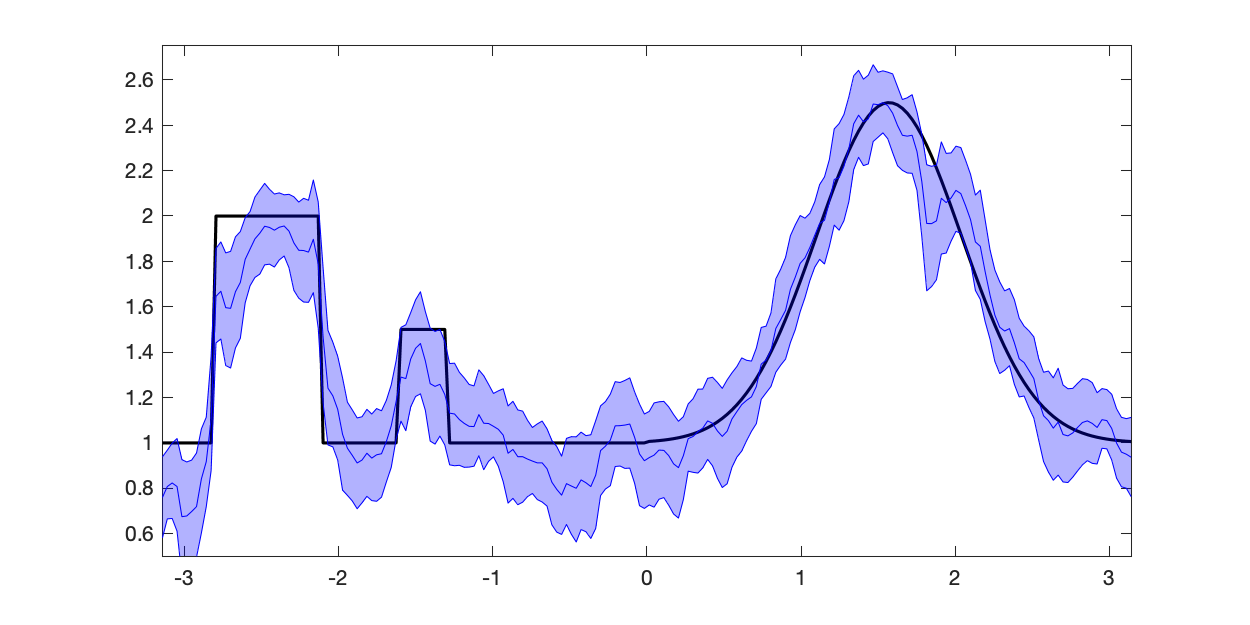}
    \caption{$F_B$ in \eqref{eq:blurtransform}}
     \end{subcaptionblock}
    \begin{subcaptionblock}{.28\textwidth}

    \includegraphics[width=.9\linewidth]{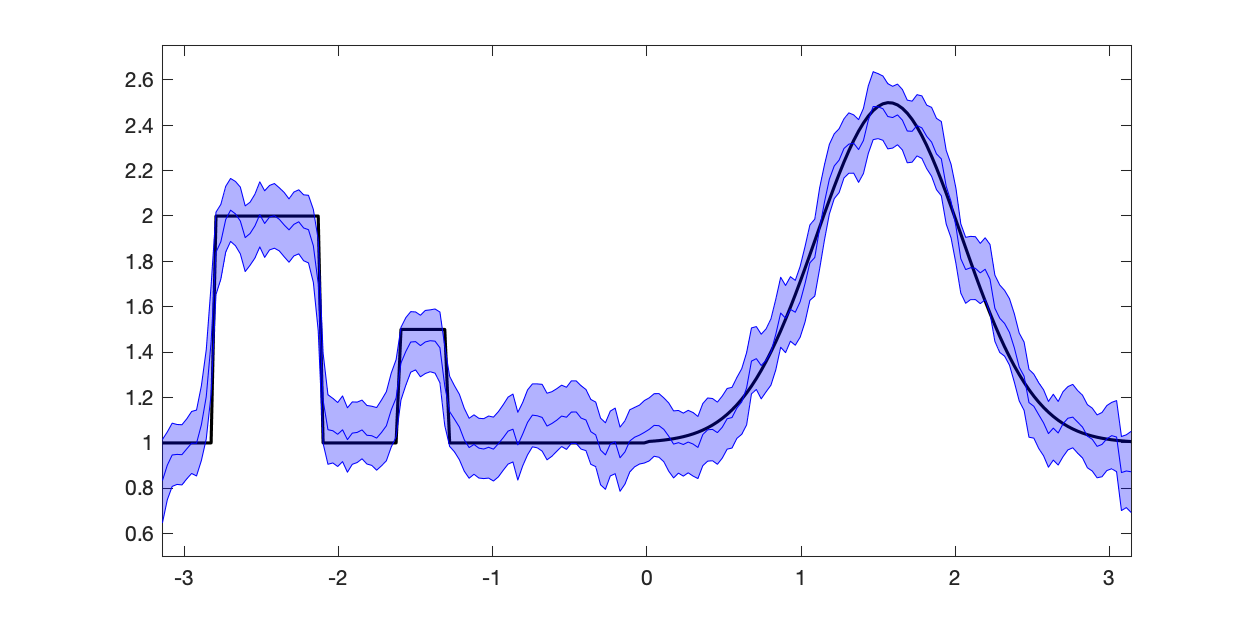}
    \caption{$F_U$ in \eqref{eq:undertransform}; $\nu=0.8$}
     \end{subcaptionblock}
     \begin{subcaptionblock}[b][2.1cm][t]{.13\textwidth}

    \includegraphics[width=.9\linewidth]{legendMag.png}
     \end{subcaptionblock}
    \caption{Magnitude recovery for ${\bm g}$ given by \eqref{eq:1D_sparseedge_testfunct}. (top) The CVBL mean along with the magnitude of LASSO solutions using different regularization parameters; (bottom) $90\%$ CI for the magnitude reconstruction. Here SNR $=20\dB$.}
    \label{fig:1D_sparseedge_mag}
\end{figure}

\Cref{fig:1D_sparseedge_mag} compares results for recovering \eqref{eq:1D_sparseedge_testfunct} from \cref{alg:propMethod} where $5000$ samples were generated to the LASSO solution in \eqref{eq:sparseedge_classical_lasso} with $\lambda = \alpha \sigma^2 /\bar{\eta}$, where $\bar{\eta}$ is again the mean of the smaples of $\upeta$ generated by the CVBL method and $\alpha = .1,1,10$. Each of the three transforms, $F_F$ \eqref{eq:Fouriertransform}, $F_B$ \eqref{eq:blurtransform}, and $F_U$ \eqref{eq:undertransform} with $\nu=0.8$ were considered for $SNR = 20\dB$. The LASSO method is clearly sensitive to the choice of regularization parameter, and is most accurate for $\alpha = 1$. 

\begin{figure}
    \centering
    \begin{subcaptionblock}{.21\textwidth}
    \centering
    \includegraphics[width=.95\linewidth]{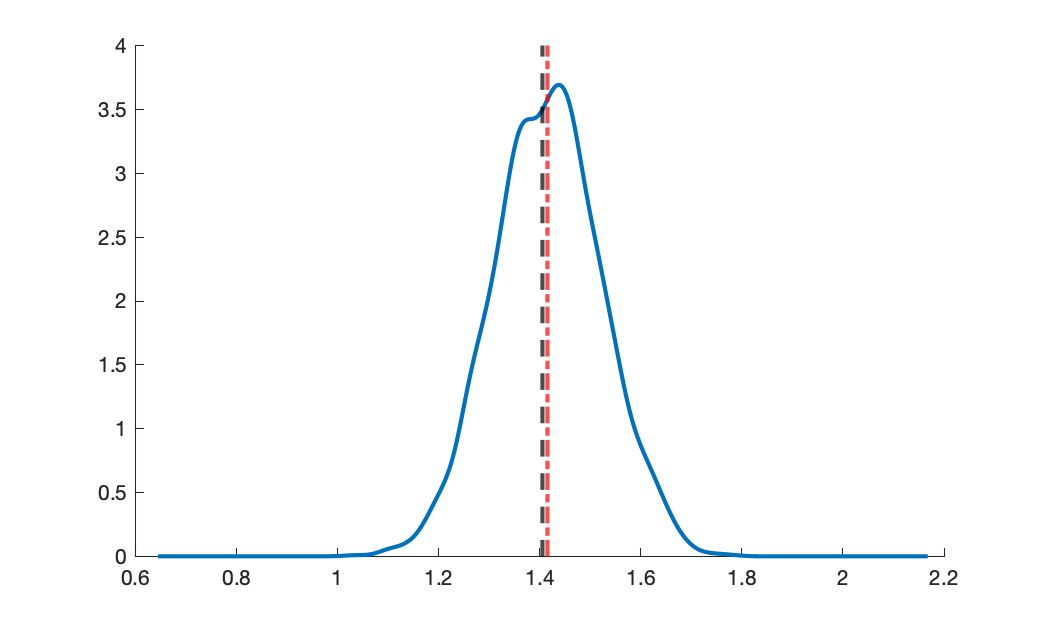}
    \caption{$F_F$ in \eqref{eq:Fouriertransform}}
    \end{subcaptionblock}
    \begin{subcaptionblock}{.21\textwidth}
    \centering
    \includegraphics[width=.95\linewidth]{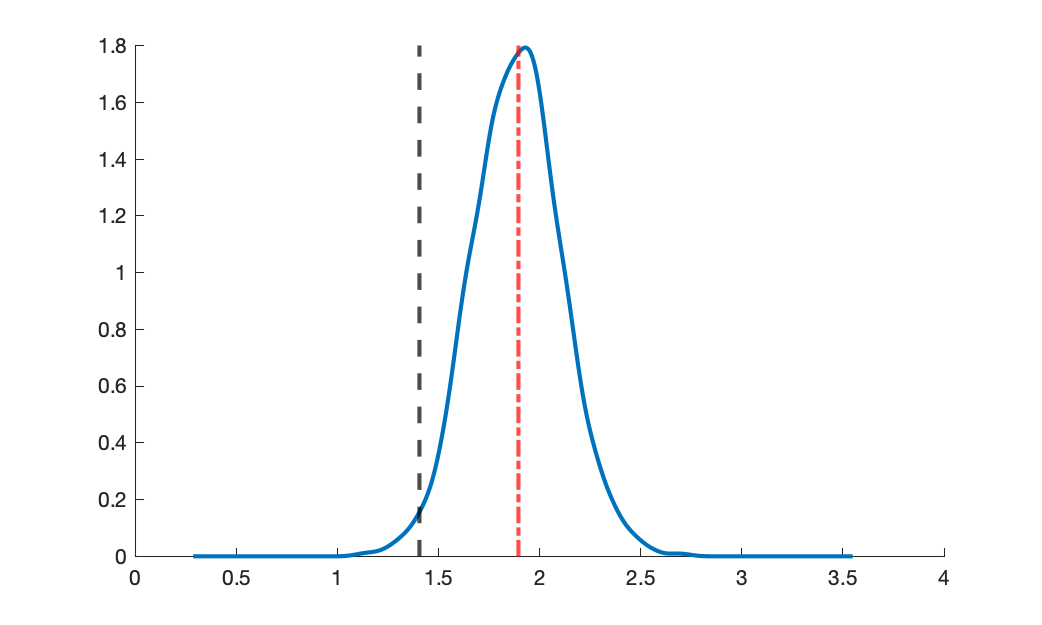}
    \caption{$F_B$ in \eqref{eq:blurtransform}}
    \end{subcaptionblock}
    \begin{subcaptionblock}{.21\textwidth}
    \centering
    \includegraphics[width=.95\linewidth]{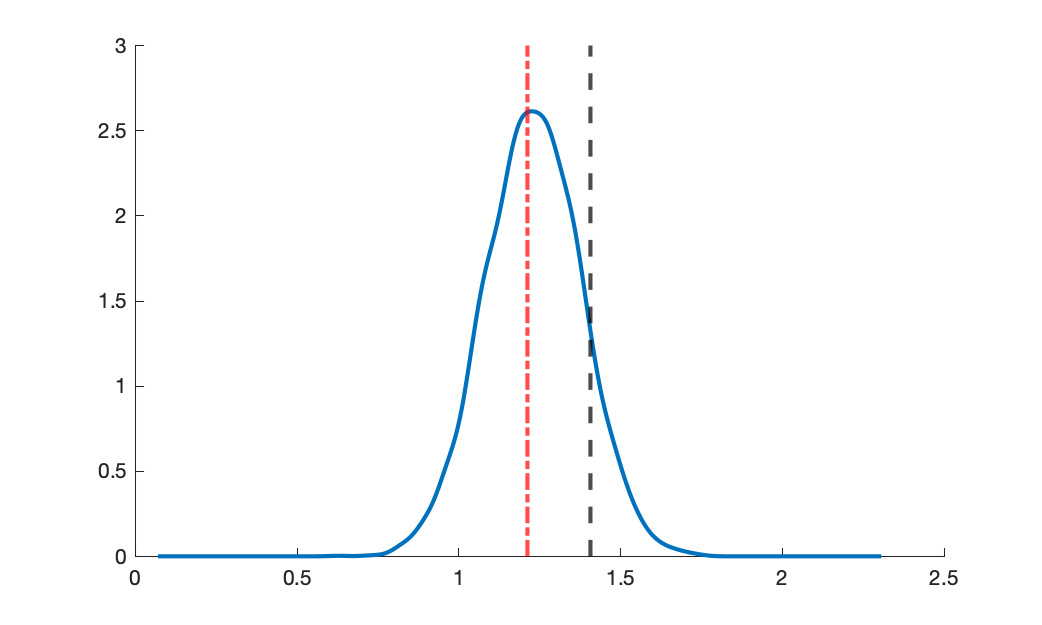}
    \caption{$F_U$ in \eqref{eq:undertransform}; $\nu=0.8$}
    \end{subcaptionblock}
    \begin{subcaptionblock}[b][2cm][t]{.13\textwidth}
    \centering
    \includegraphics[width=.95\linewidth]{legendPhase.png}
     \end{subcaptionblock}
    \caption{Phase estimate in the sparse transform case of a random pixel. True phase (black dash), LASSO phase point estimate (red dash-dot), and CVBL marginal posterior for the phase. SNR $=20\dB$.}
    \label{fig:1D_sparseedge_phase}
\end{figure}

\begin{figure}[h!]
    \centering
    \begin{subcaptionblock}{.21\textwidth}
    \centering
    \includegraphics[width=.95\linewidth]{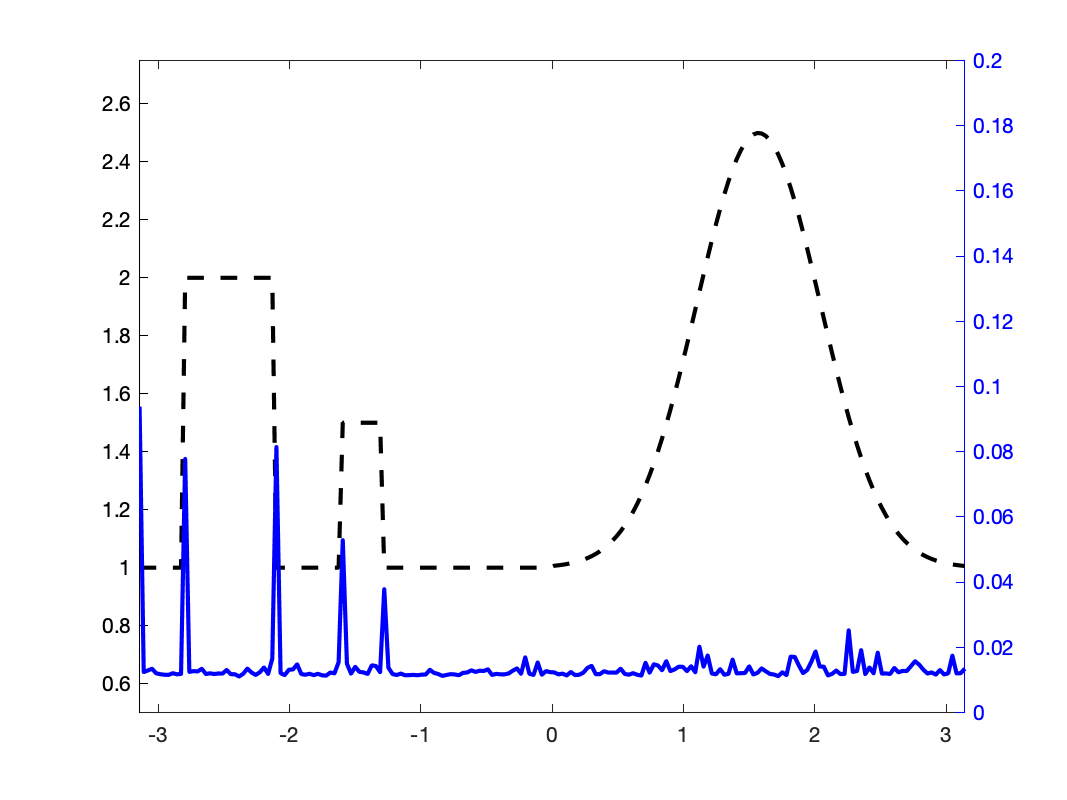}
    \caption{$F_F$ in \eqref{eq:Fouriertransform}}
    \end{subcaptionblock}
    \begin{subcaptionblock}{.21\textwidth}
    \centering
    \includegraphics[width=.95\linewidth]{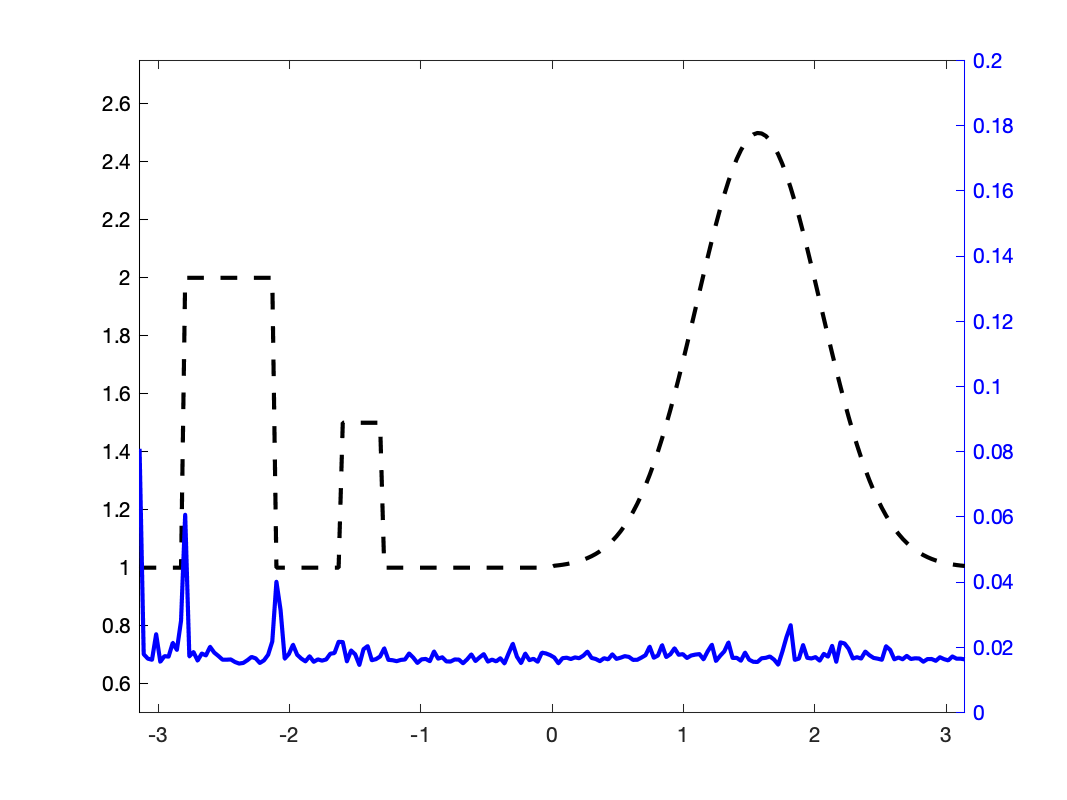}
    \caption{$F_B$ in \eqref{eq:blurtransform}}
    \end{subcaptionblock}
    \begin{subcaptionblock}{.21\textwidth}
    \centering
    \includegraphics[width=.95\linewidth]{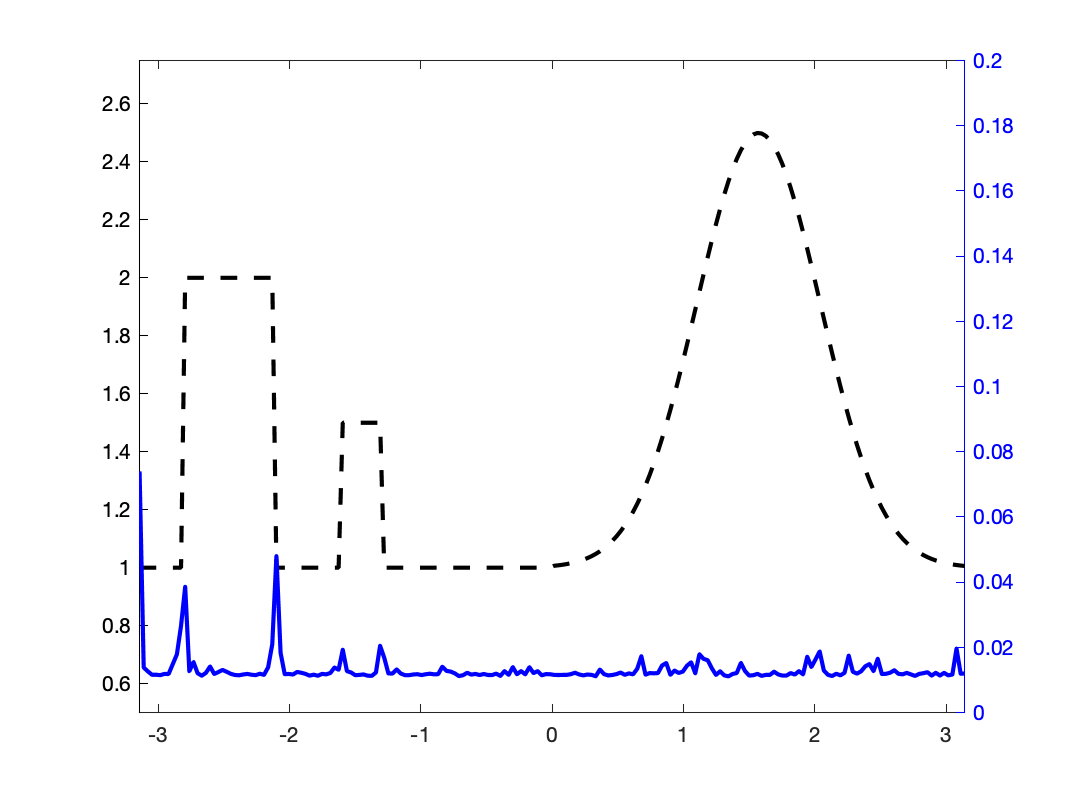}
    \caption{$F_U$ in \eqref{eq:undertransform}; $\nu=0.8$}
    \end{subcaptionblock}
    \begin{subcaptionblock}[b][2.2cm][t]{.13\textwidth}
    \centering
    \includegraphics[width=.95\linewidth]{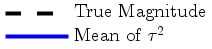}
     \end{subcaptionblock}
    \caption{Mean of $\tau^2$ (right axis) compared with the true magnitude (left axis).  SNR $=20\dB$.}
    \label{fig:1D_sparseedge_tausq}
\end{figure}

\Cref{fig:1D_sparseedge_phase} demonstrates that the CVBL method recovers marginal posterior density functions for the phase that are consistent with the estimates calculated by  \eqref{eq:sparseedge_classical_lasso}. The true phase is also located in the regions of large mass generated by the corresponding kernel density estimation. \Cref{fig:1D_sparseedge_tausq} shows the mean of the $\tau^2$ samples \eqref{eq:tau_posterior_complex} for each experiment, which as expected is  largest in support regions of the sparse domain.

\subsection{Noise Study}\label{subsec:numerics_noisestudy}
We now analyze the effect of noise on \Cref{alg:propMethod}. To this end for both the sparse magnitude case and the signal corresponding to \eqref{eq:1D_sparseedge_testfunct} we consider $10\dB \le $ SNR $ \le 30\dB$.

\begin{figure}[h!]
    \centering
    \begin{subcaptionblock}{.4\textwidth}
    \centering
    \includegraphics[width=.65\linewidth]{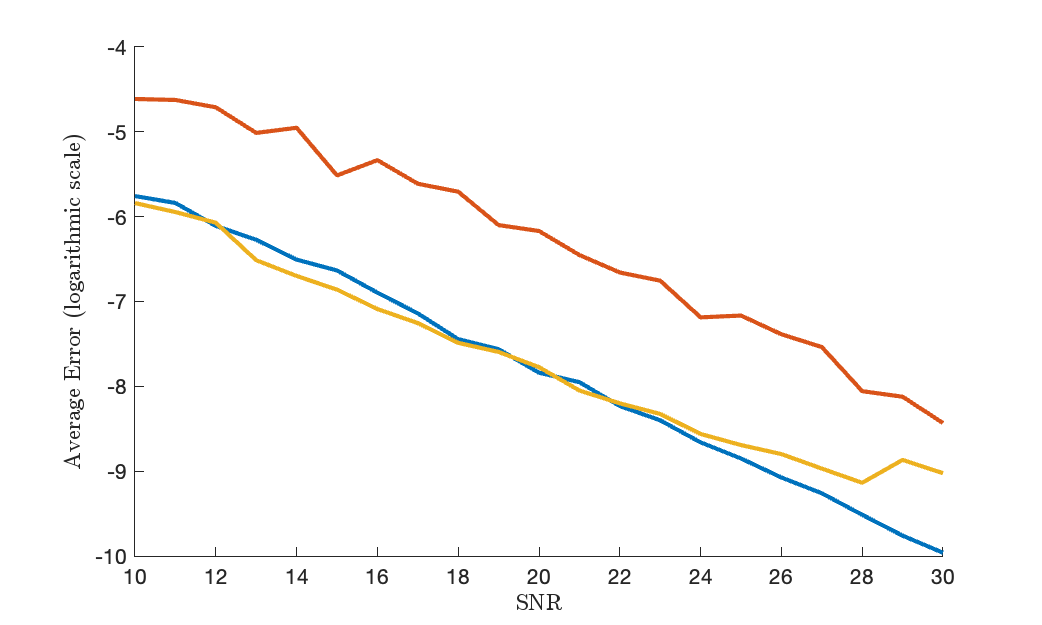}
    \caption{Sparse signal}
    \end{subcaptionblock}
    \begin{subcaptionblock}{.4\textwidth}
    \centering
    \includegraphics[width=.65\linewidth]{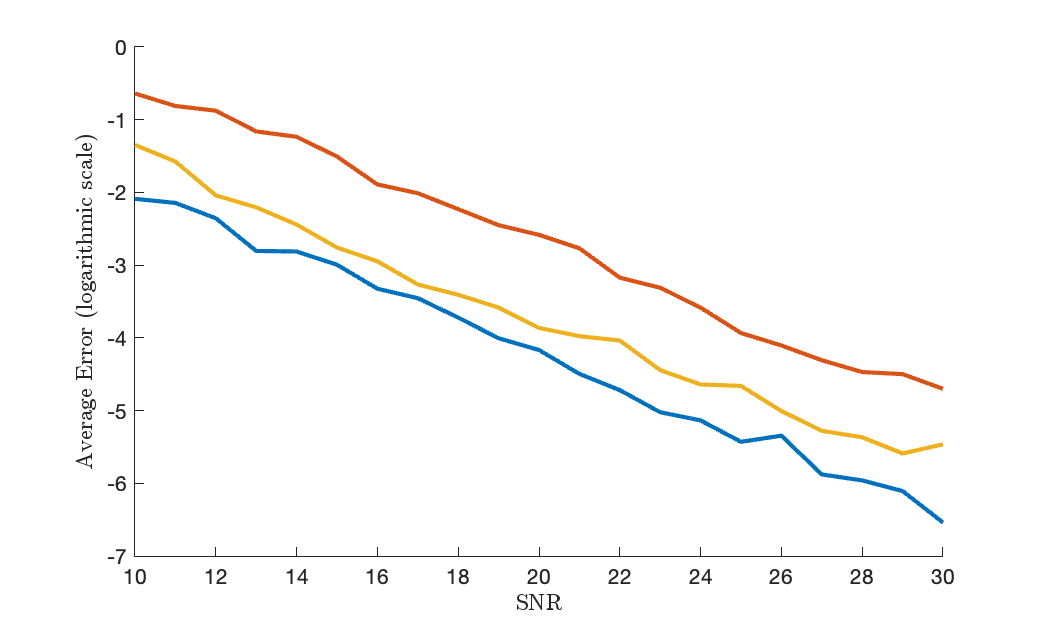}
    \caption{Sparse transform of magnitude, \eqref{eq:1D_sparseedge_testfunct}}
    \end{subcaptionblock}
    \begin{subcaptionblock}[b][2.4cm][t]{.08\textwidth}
    \centering
    \includegraphics[width=.8\linewidth]{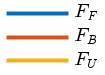}
     \end{subcaptionblock}
    \caption{Average error in the sample mean generated using \Cref{alg:propMethod} for increasing SNR. The vertical axes are in logarithmic scale.}
    \label{fig:1D_SNR_study}
\end{figure}

\Cref{fig:1D_SNR_study} demonstrates that the average error is similar for all three forward operators across the range of SNR values, with a greater overall range of errors with respect to the forward operator in the sparse transform case. These results are consistent with what is observed in \cref{fig:1D_sparseedge_mag} as well as what is  apparent in \cref{fig:1D_sparseedge_mag} (b).  In particular the approximation in the region  $[-1.7,-1]$  is not well resolved. More insight is provided in \cref{fig:1D_sparseedge_tausq} (b), where it is evident that the support in the sparse domain in that range is not clearly identified. 

\begin{figure}[h!]
    \centering
    \begin{subcaptionblock}{.21\textwidth}
    \centering
    \includegraphics[width=\linewidth]{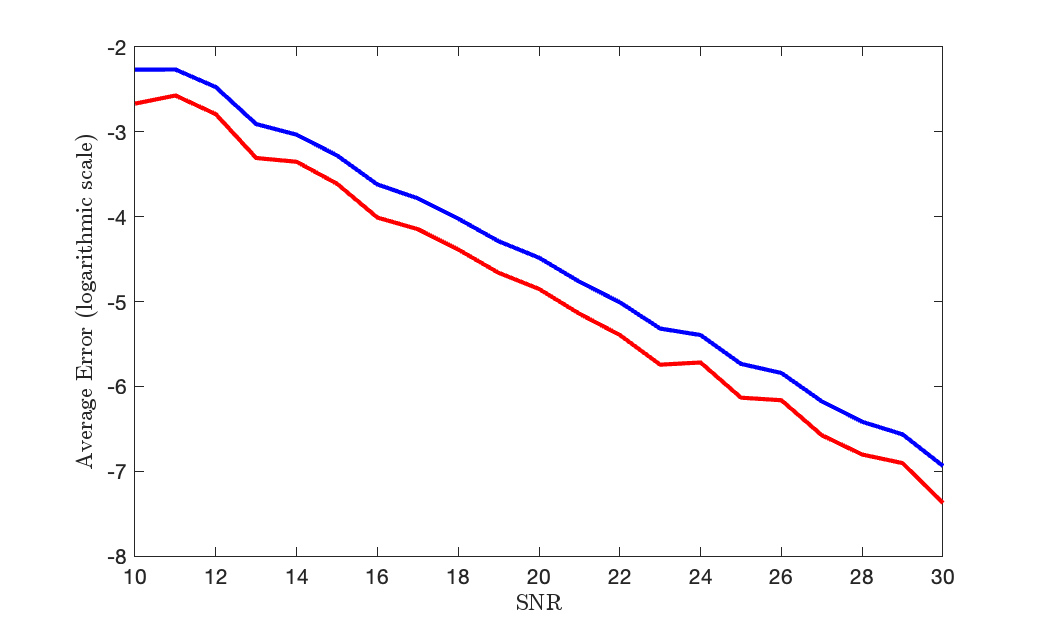}
    \caption{$F_F$ in \eqref{eq:Fouriertransform}}
    \end{subcaptionblock}
    \begin{subcaptionblock}{.21\textwidth}
    \centering
    \includegraphics[width=\linewidth]{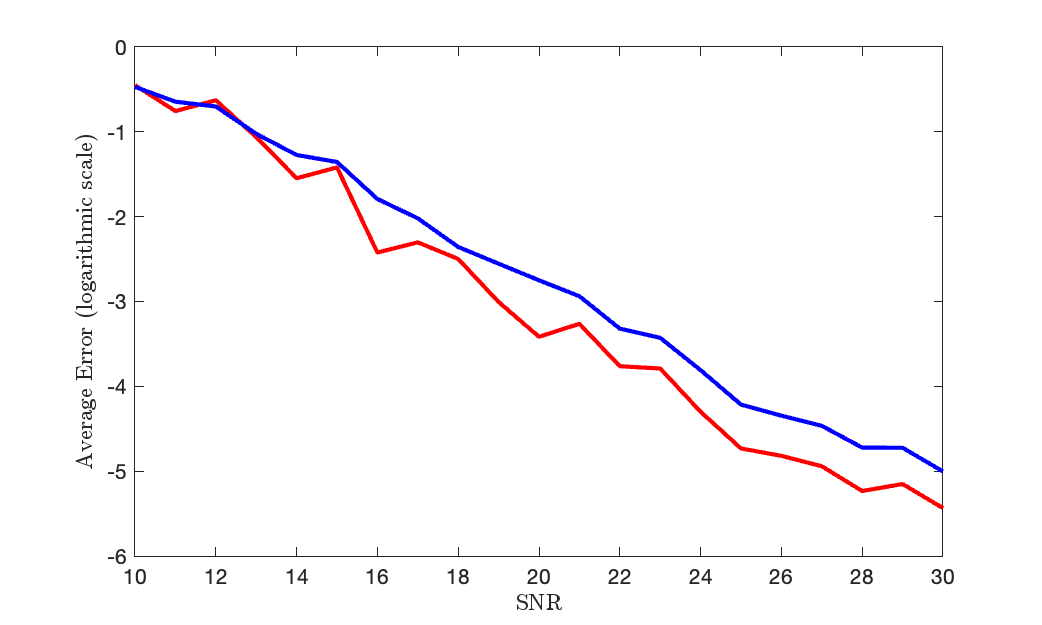}
    \caption{$F_B$ in \eqref{eq:blurtransform}}
    \end{subcaptionblock}
    \begin{subcaptionblock}{.21\textwidth}
    \centering
    \includegraphics[width=\linewidth]{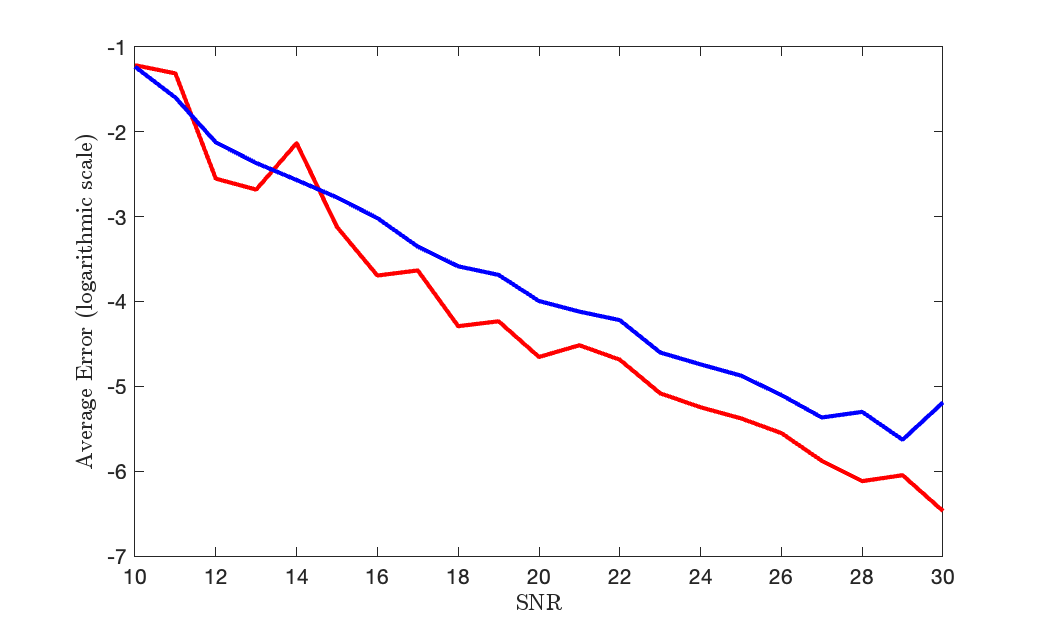}
    \caption{$F_U$ in \eqref{eq:undertransform}; $\nu=0.8$}
    \end{subcaptionblock}
    \begin{subcaptionblock}{.21\textwidth}
    \centering
    \includegraphics[width=\linewidth]{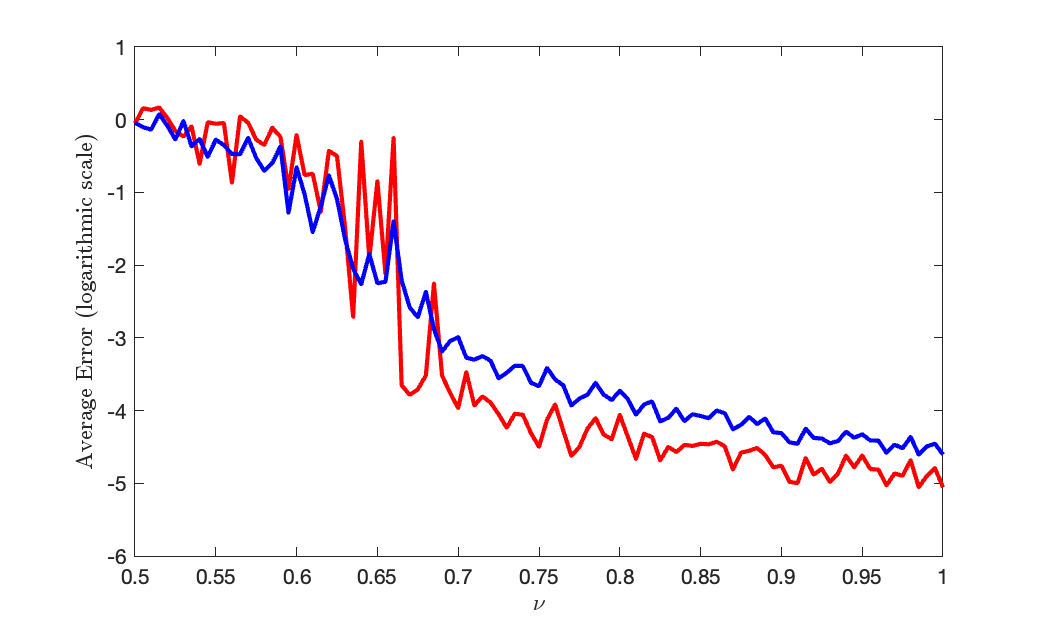}
    \caption{$F_U$; $0.5\leq\nu\leq1$}
    \end{subcaptionblock}
    \begin{subcaptionblock}[b][2.4cm][t]{.12\textwidth}
    \centering
    \includegraphics[width=\linewidth]{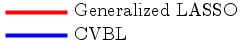}
     \end{subcaptionblock}
    \caption{Comparison of the average phase error in mean of samples generated using CVBL with the LASSO solution \eqref{eq:sparseedge_classical_lasso} for the test function \eqref{eq:1D_sparseedge_testfunct} at (a)-(c) increasing SNR and (d) increasing sampling rate $0.5\leq\nu\leq1$ with SNR $=20\dB$. The vertical axes are in logarithmic scale.}
    \label{fig:1D_error_study}
\end{figure}

\Cref{fig:1D_error_study} (a)-(c) compares the average phase error over a range of SNR values for the CVBL and generalized LASSO point estimates.  While the LASSO method outperforms \cref{alg:propMethod} for each choice of $F$ with increasing SNR, the error difference is neglible for low SNR values for $F_B$ and $F_U$.   Hence we see that the optimal {\em point estimate} recovery algorithm essentially depends on the SNR and sampling rate of the observable data. {\em Uncertainty information}, however, is only acquired when using CVBL, as the generalized LASSO technique does not infer the phase information.

Finally, \cref{fig:1D_error_study} (d) compares the average phase error in the sample mean for the function with magnitude given by \eqref{eq:1D_sparseedge_testfunct} using CVBL with the LASSO solution for the undersampled observable data case.  Specifically, the data are obtained using  $F_U$ in \eqref{eq:undertransform} for a range of sampling rates $\nu$ while  the SNR is held constant at $20\dB$. While consistent with \cref{fig:1D_error_study} (c), this result also suggests that for smaller values of $\nu$ (more undersampling) the CVBL method provides on average a lower average phase error than the generalized LASSO.

\subsection{2D Experiments}\label{subsec:2Dexperiments}
Our 2D experiments consider sparsity in the signal magnitude with SNR $= 0\dB$ and sparsity in the signal magnitude gradient  with  SNR $=15\dB$ and SNR $= 25\dB$.

\begin{figure}[h!]
    \centering
    \begin{subcaptionblock}{.24\textwidth}
    \centering
    \includegraphics[width=1\linewidth]{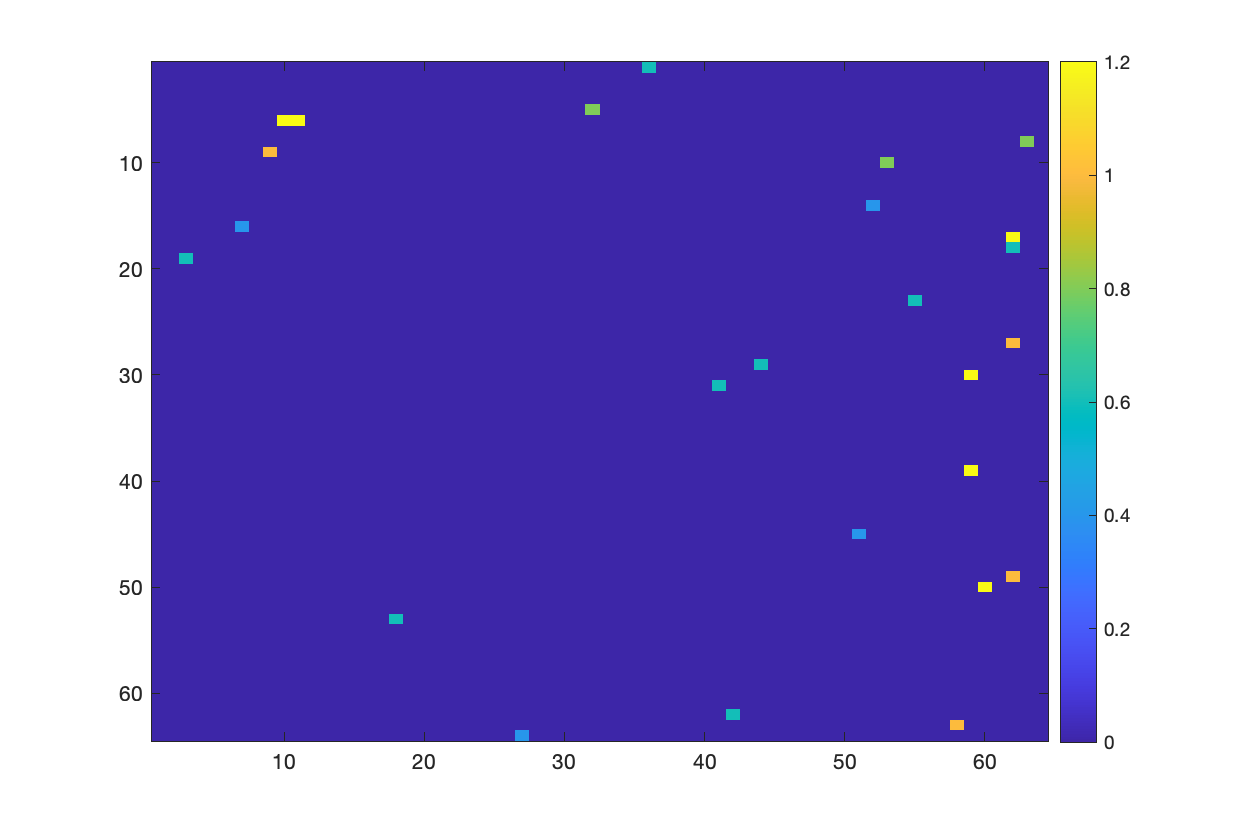}
    \caption{True magnitude}
    \end{subcaptionblock}
    \begin{subcaptionblock}{.24\textwidth}
    \centering
    \includegraphics[width=1\linewidth]{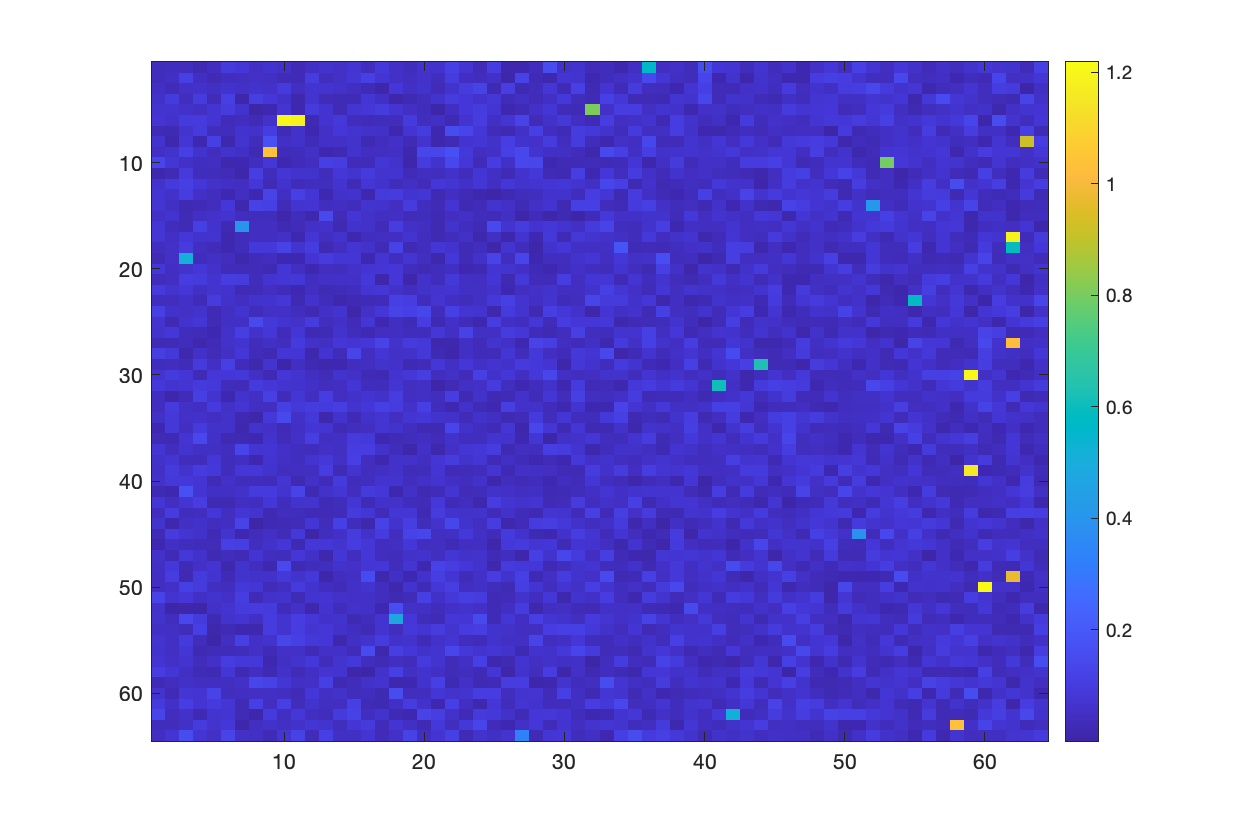}
    \caption{MLE}
    \end{subcaptionblock}
    \begin{subcaptionblock}{.24\textwidth}
    \centering
    \includegraphics[width=1\linewidth]{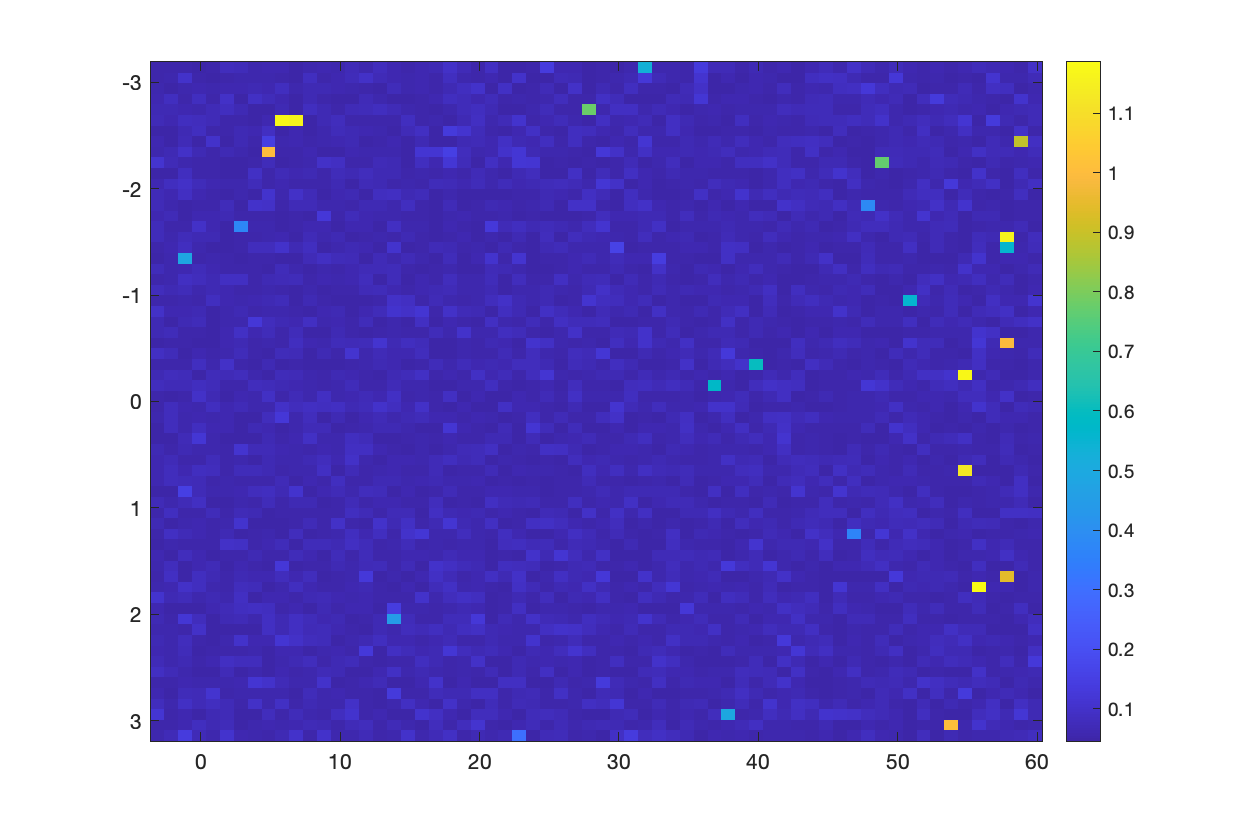}
    \caption{Mean posterior}
    \end{subcaptionblock}
    \caption{Maximum likelihood and mean posterior point estimates of a sparse magnitude complex-valued signal for forward operator $F_F$ in \eqref{eq:Fouriertransform}. SNR $= 0\dB$.}
    \label{fig:2D_signal_FF}
\end{figure}

\begin{figure}
    \centering
    \begin{subcaptionblock}{.24\textwidth}
    \centering
    \includegraphics[width=1\linewidth]{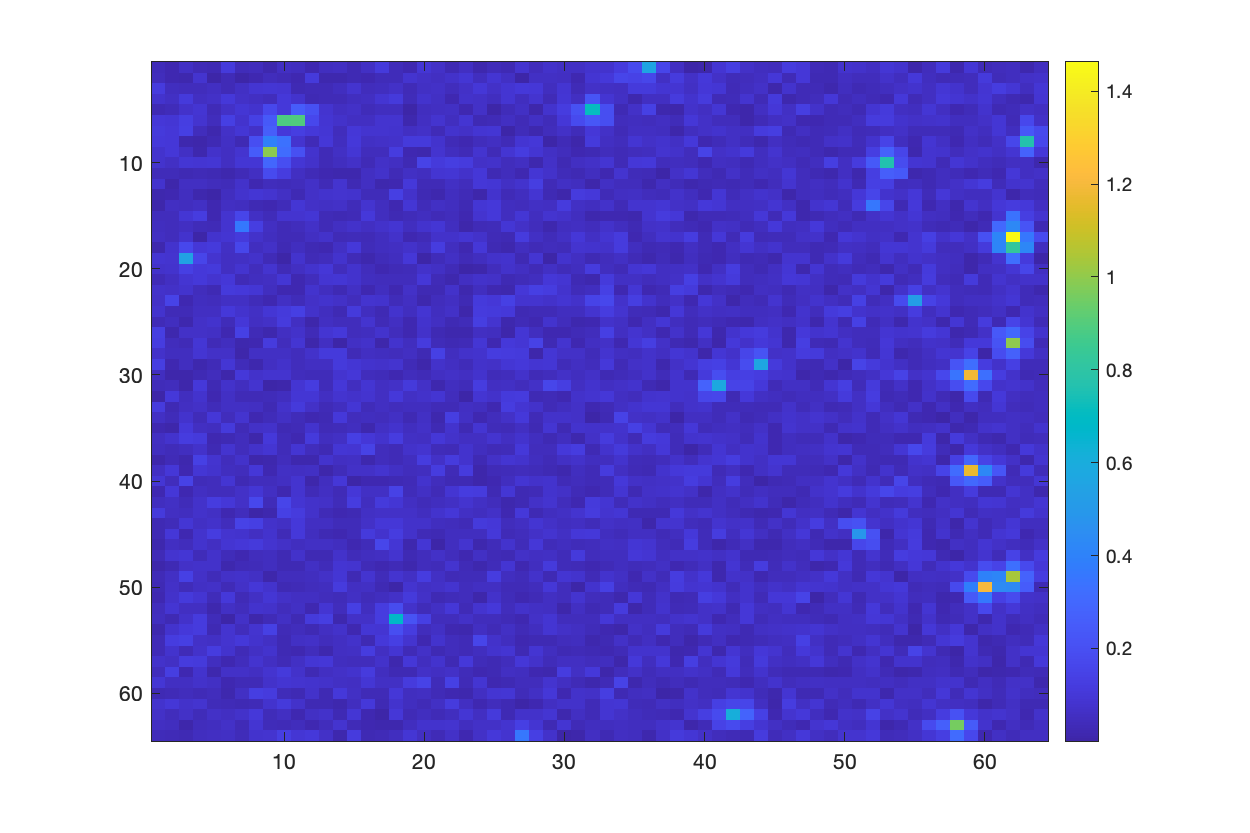}
    \caption{MLE;  $F_B$}
    \end{subcaptionblock}
    \begin{subcaptionblock}{.24\textwidth}
    \centering
    \includegraphics[width=1\linewidth]{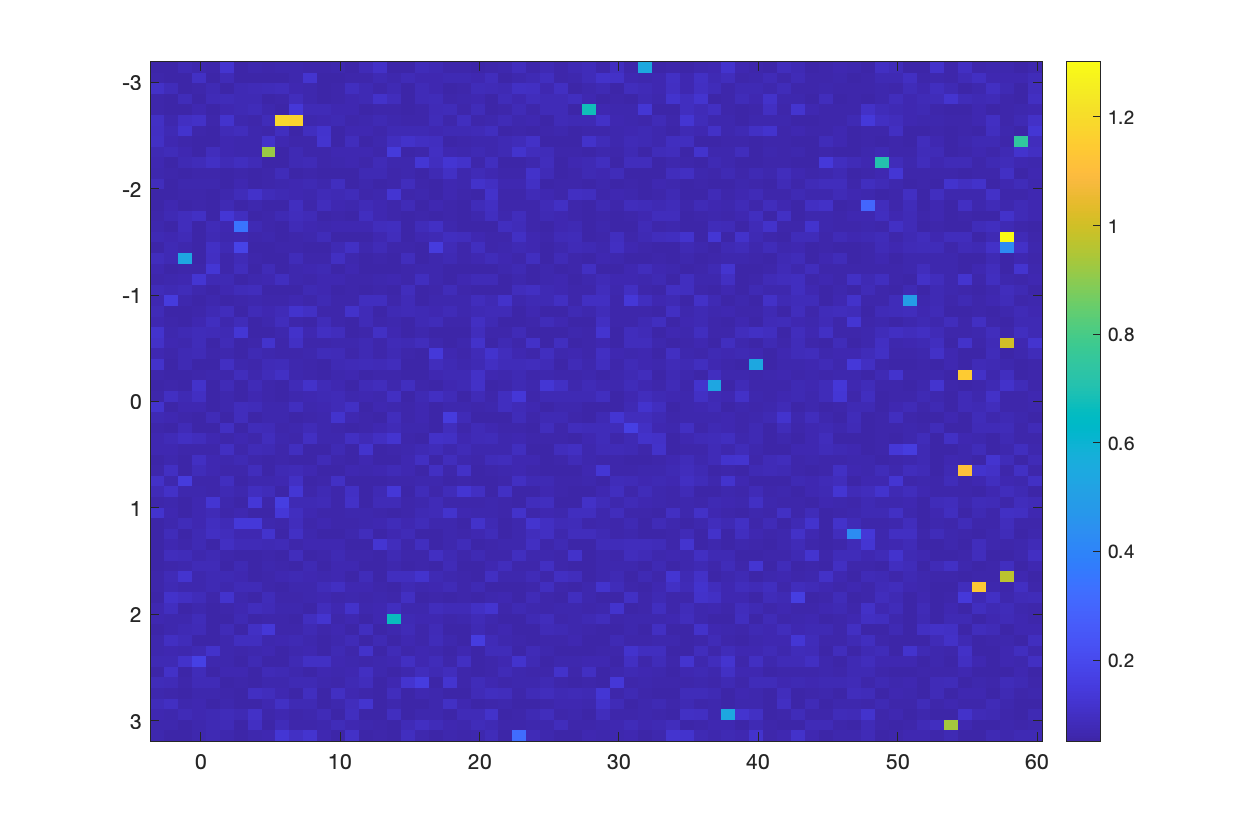}
    \caption{Mean posterior; $F_B$}
    \end{subcaptionblock}
    \begin{subcaptionblock}{.24\textwidth}
    \centering
    \includegraphics[width=1\linewidth]{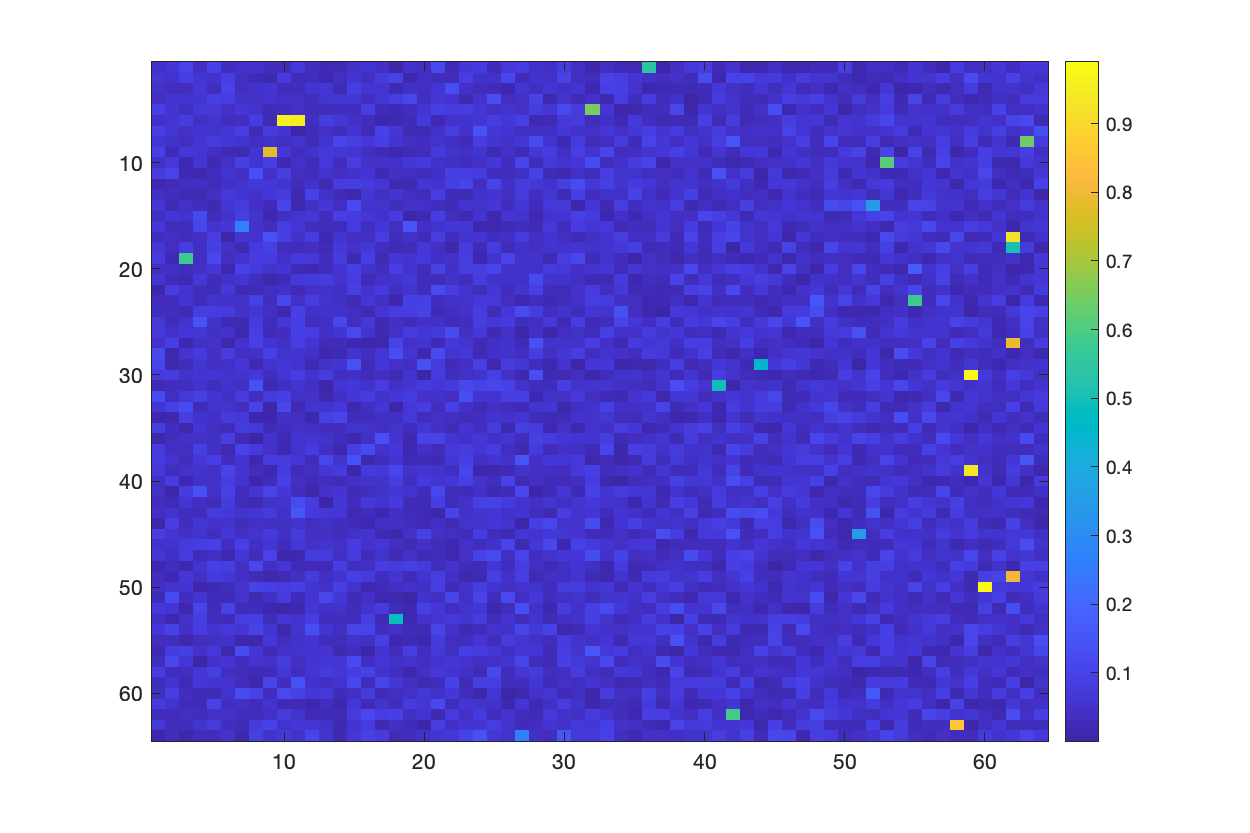}
    \caption{MLE; $F_U$}
    \end{subcaptionblock}
    \begin{subcaptionblock}{.24\textwidth}
    \centering
    \includegraphics[width=1\linewidth]{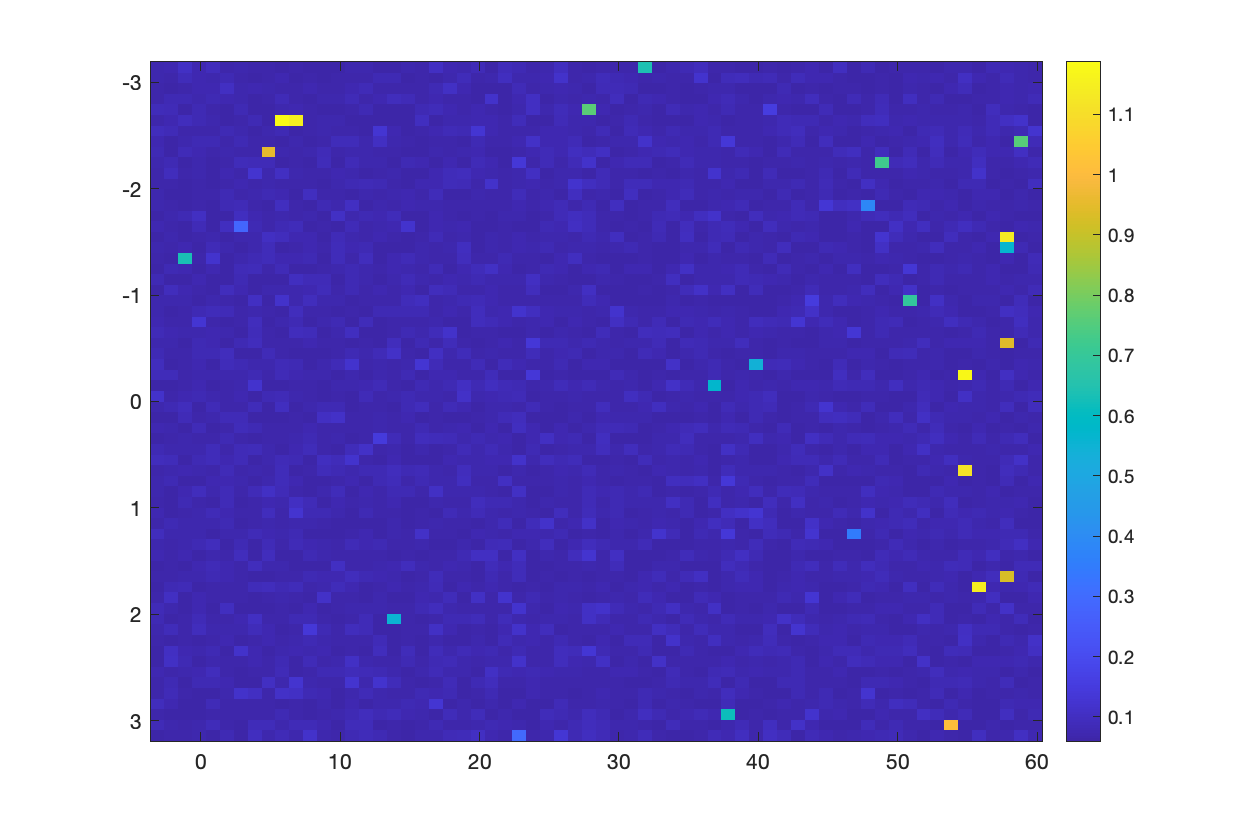}
    \caption{Mean posterior; $F_U$}
     \end{subcaptionblock}
    \caption{Maximum likelihood and mean posterior point estimates recoveries of a sparse magnitude complex-valued signal. Forward operator  (a-b) $F_B$ in \eqref{eq:blurtransform} (c-d) $F_U$ $\nu = 0.8$ in \eqref{eq:undertransform}. SNR = $0\dB$.}
    \label{fig:2D_signal_FB_FU}
\end{figure}

\Cref{fig:2D_signal_FF} displays point estimates for the magnitude using $F_F$ while \Cref{fig:2D_signal_FB_FU} shows the corresponding results for $F_B$ and $F_U$. In each case the same $64\cross64$ sparse images is used as input and we simulate $5000$ samples of the posterior. We observe that much of the background noise is suppressed while the fidelity of the support of the signal is maintained when compared to the maximum likelihood point estimate. 

\begin{figure}
    \begin{subcaptionblock}{.23\textwidth}
    \includegraphics[width=1\linewidth]{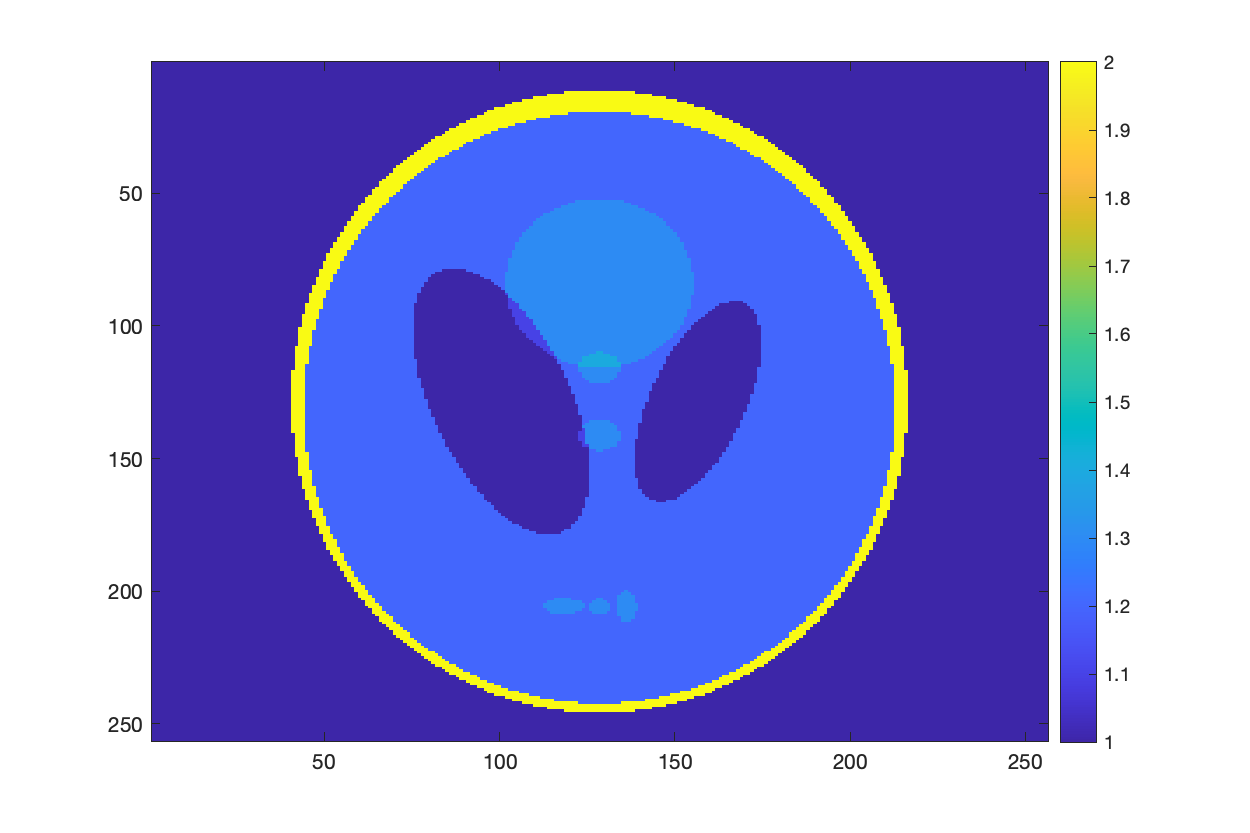}
     \caption{Ground truth}
    \end{subcaptionblock}
    \centering
    \begin{subcaptionblock}{.23\textwidth}
    \centering
    \includegraphics[width=1\linewidth]{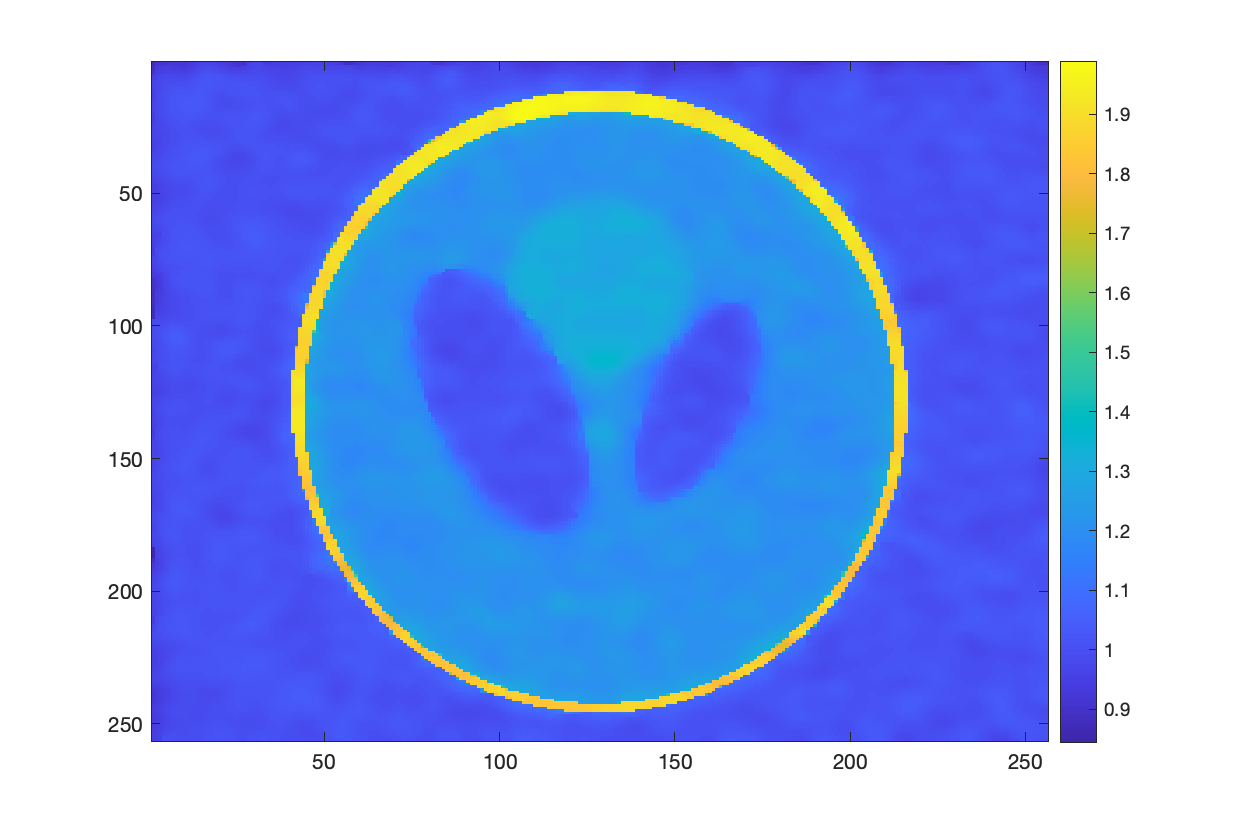}
    \caption{Generalized LASSO}
    \end{subcaptionblock}
    \begin{subcaptionblock}{.23\textwidth}
    \centering
    \includegraphics[width=1\linewidth]{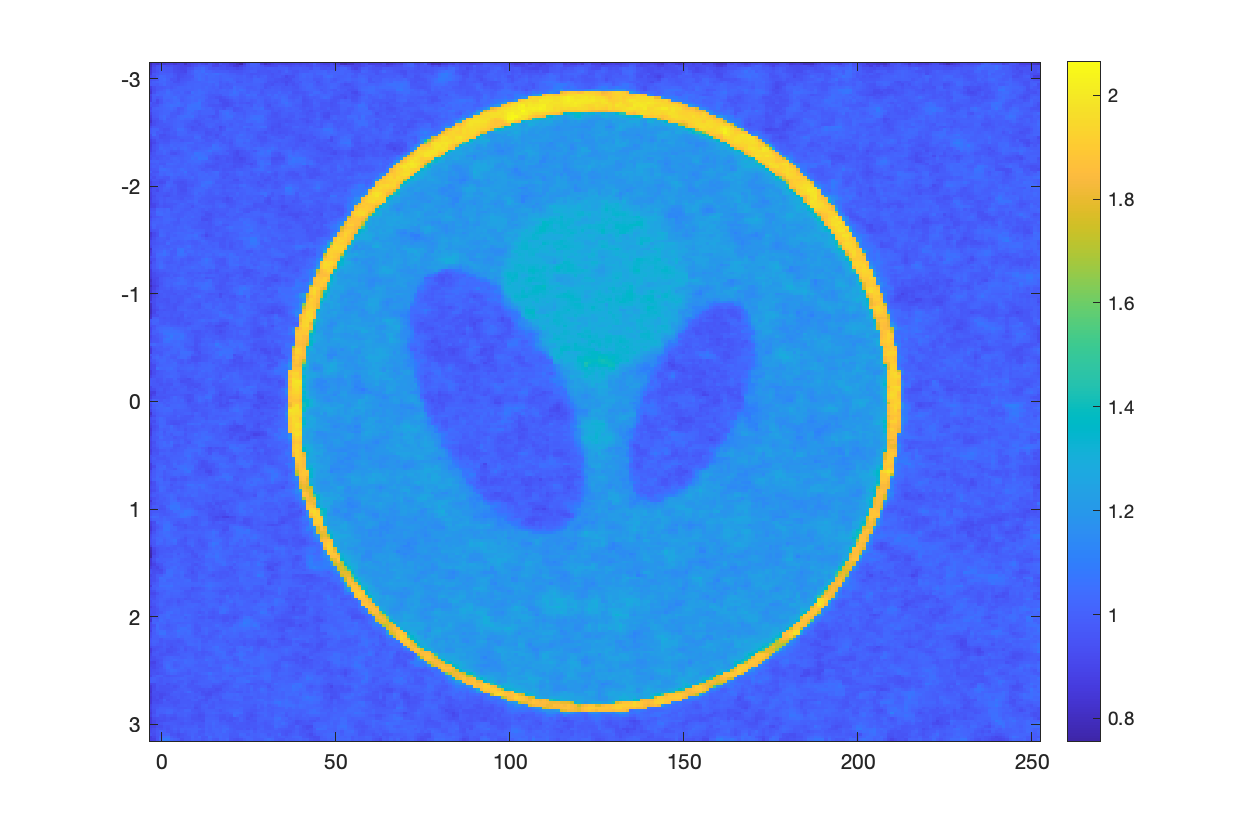}
    \caption{Posterior mean}
    \end{subcaptionblock}
    \begin{subcaptionblock}{.23\textwidth}
    \centering
    \includegraphics[width=1\linewidth]{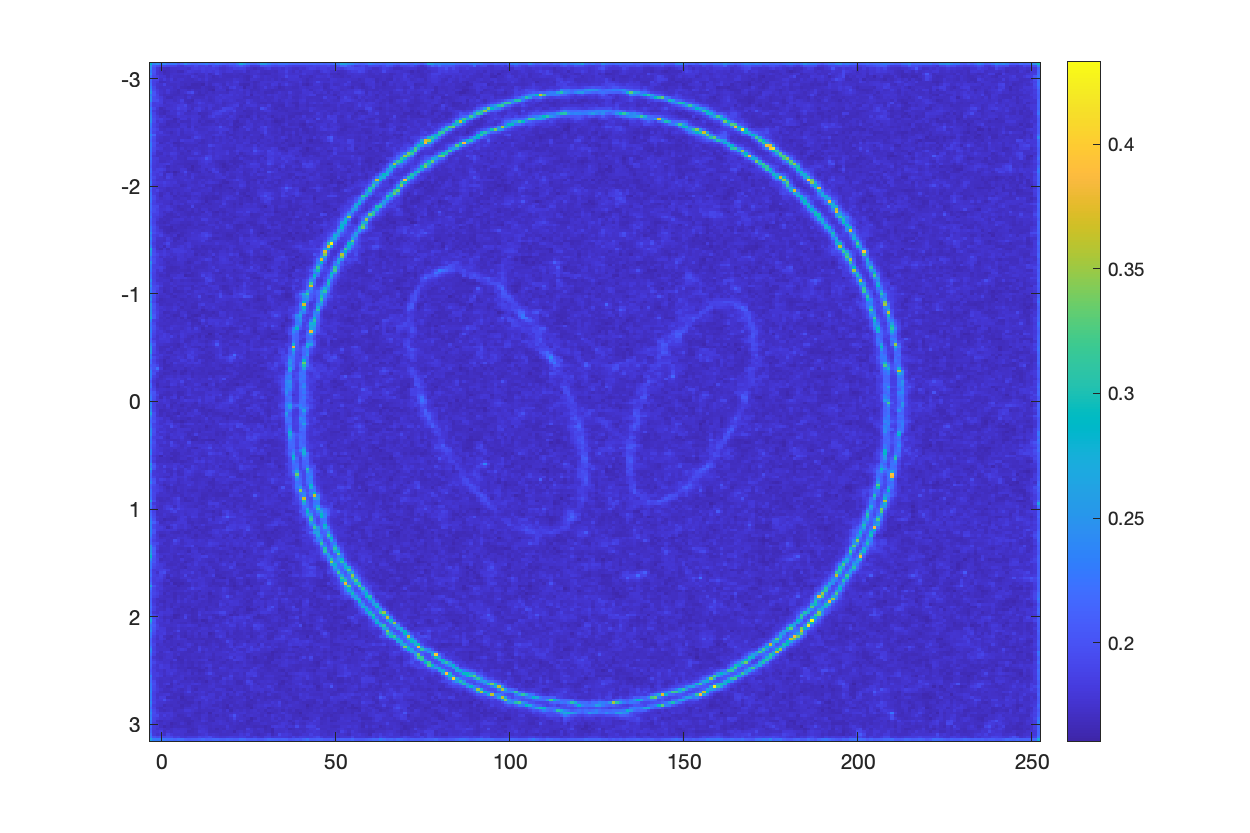}
    \caption{$90\%$ CI}
    \end{subcaptionblock}
    \begin{subcaptionblock}{.23\textwidth}
    \includegraphics[width=1\linewidth]{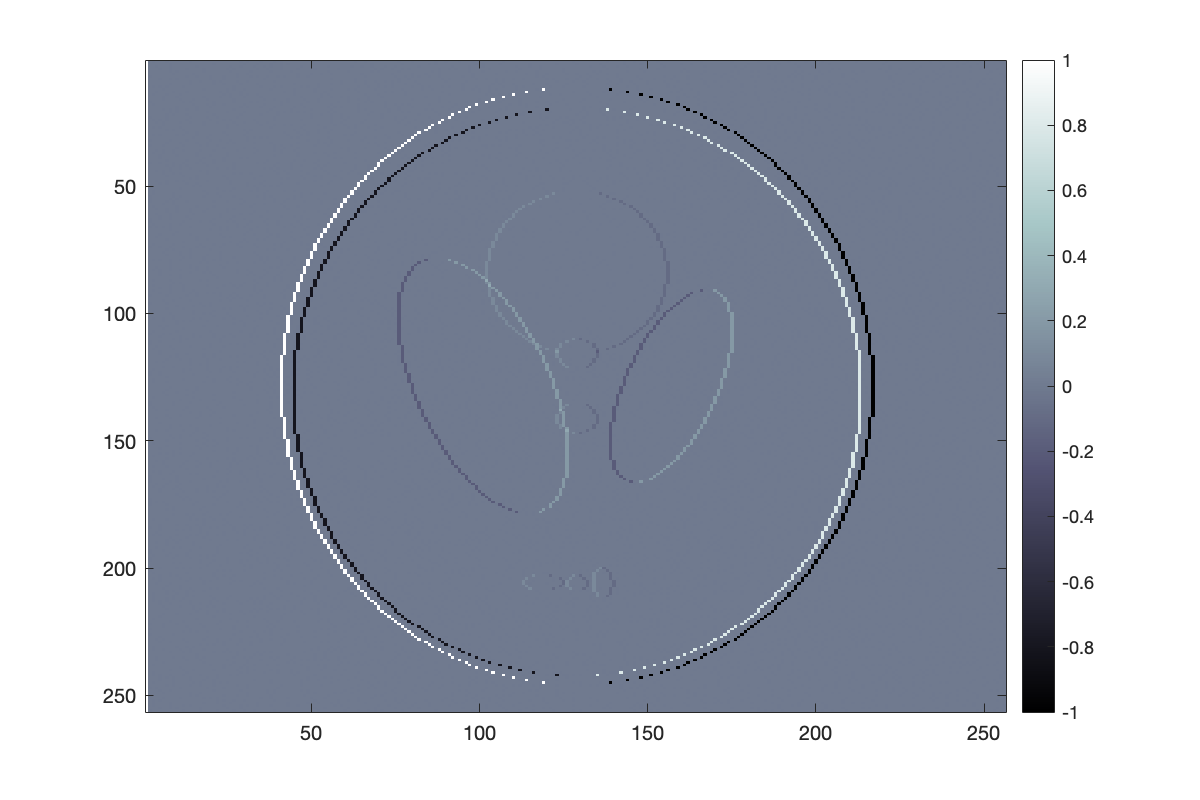}
     \caption{Horizontal differences}
    \end{subcaptionblock}
    \begin{subcaptionblock}{.23\textwidth}
    \centering
    \includegraphics[width=1\linewidth]{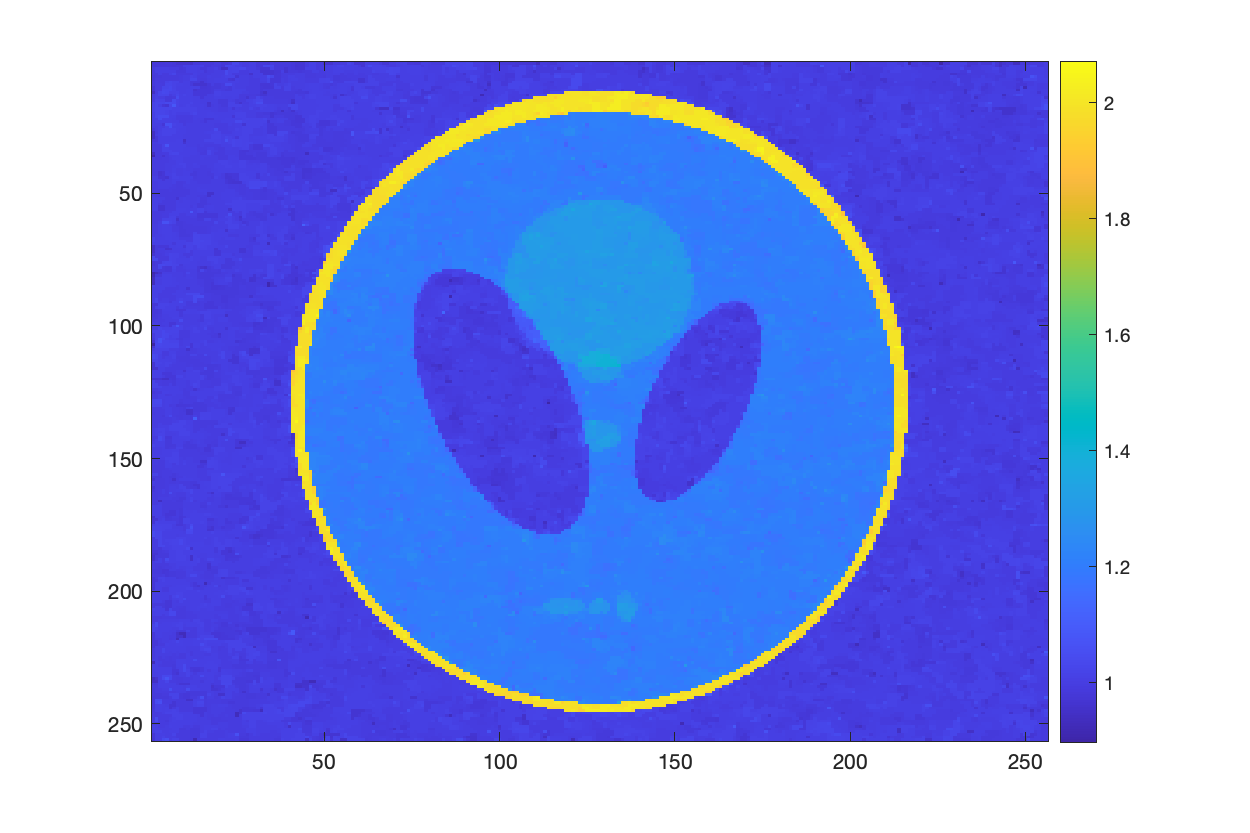}
    \caption{Generalized LASSO}
    \end{subcaptionblock}
    \begin{subcaptionblock}{.23\textwidth}
    \centering
    \includegraphics[width=1\linewidth]{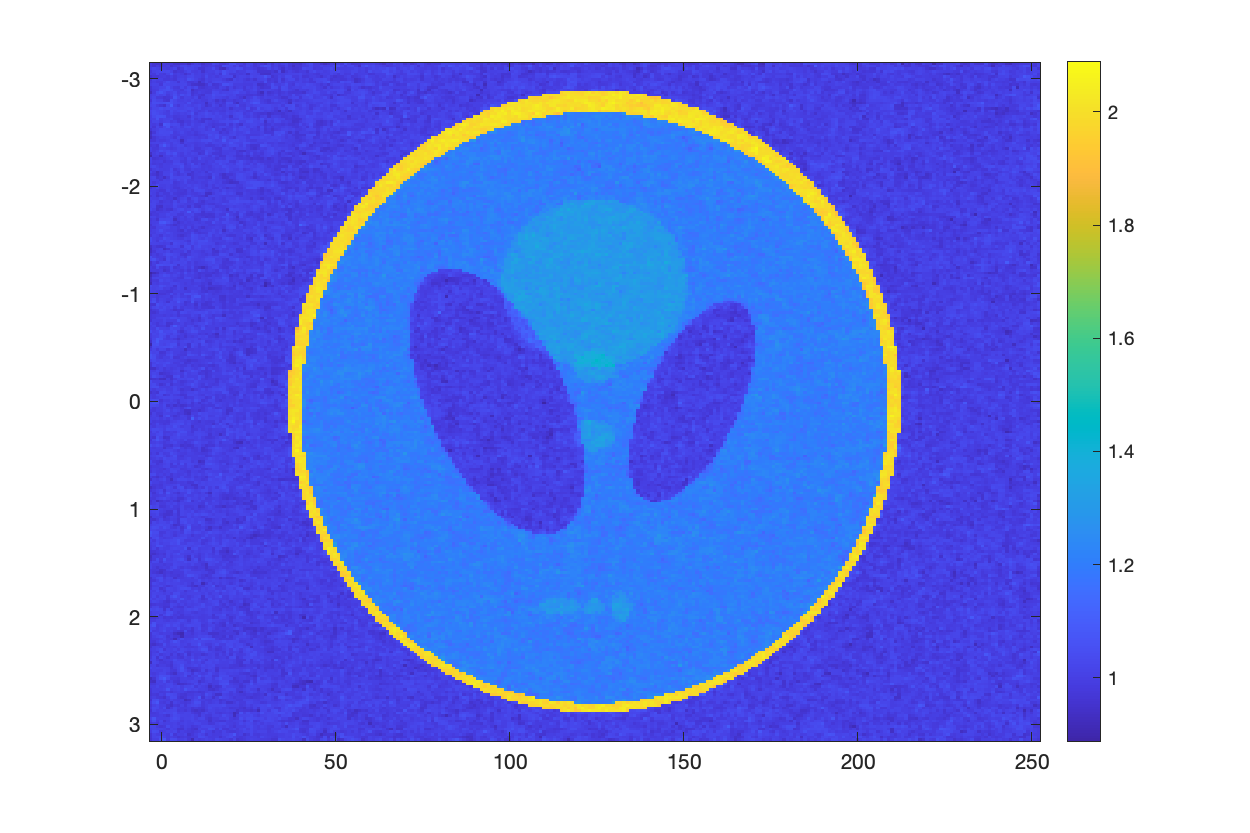}
    \caption{Posterior mean}
    \end{subcaptionblock}
    \begin{subcaptionblock}{.23\textwidth}
    \centering
    \includegraphics[width=1\linewidth]{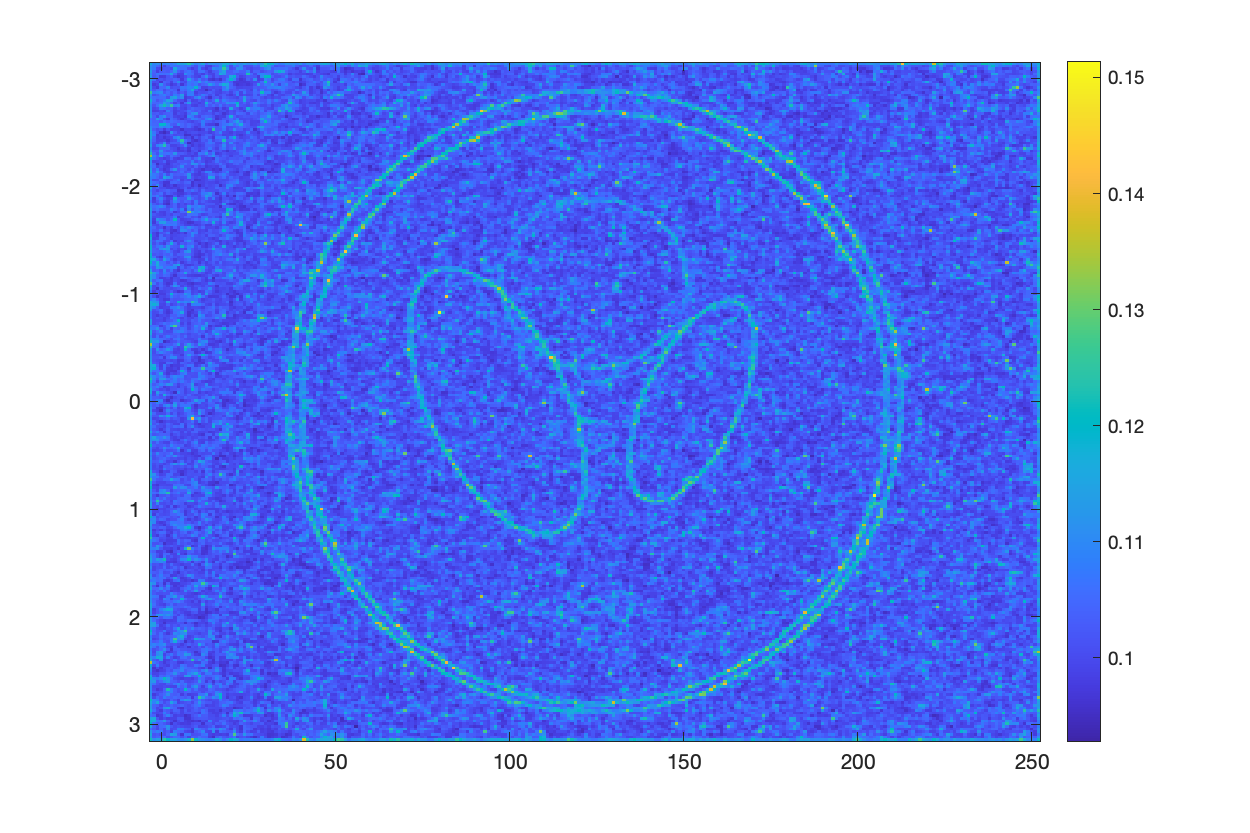}
    \caption{$90\%$ CI}
    \end{subcaptionblock}
    \caption{Magnitude recovery of the complex-valued Shepp Logan phantom  using $F_F$ in \eqref{eq:Fouriertransform}. (a) Shepp-Logan phantom; (e) horizontal differences of the ground truth. (b)-(d) and (f)-(h) correspond to SNR values of $15\dB$ and $25\dB$, respectively. $\hat{\eta}=10^{-2}$.}
    \label{fig:2D_transform_FF}
\end{figure}

\begin{figure}[h!]
    \centering
    \begin{subcaptionblock}{.23\textwidth}
    \centering
    \includegraphics[width=1\linewidth]{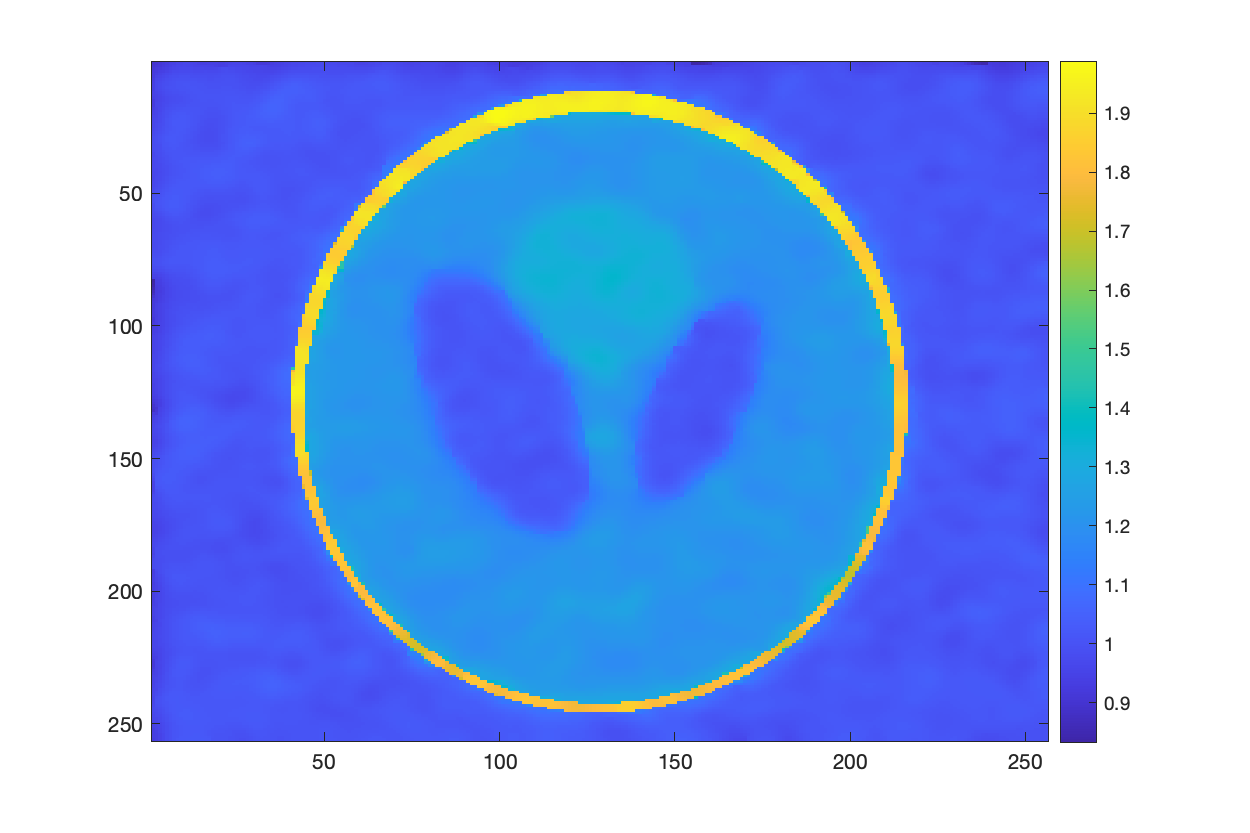}
    \caption{Generalized LASSO}
    \end{subcaptionblock}
    \begin{subcaptionblock}{.23\textwidth}
    \centering
    \includegraphics[width=1\linewidth]{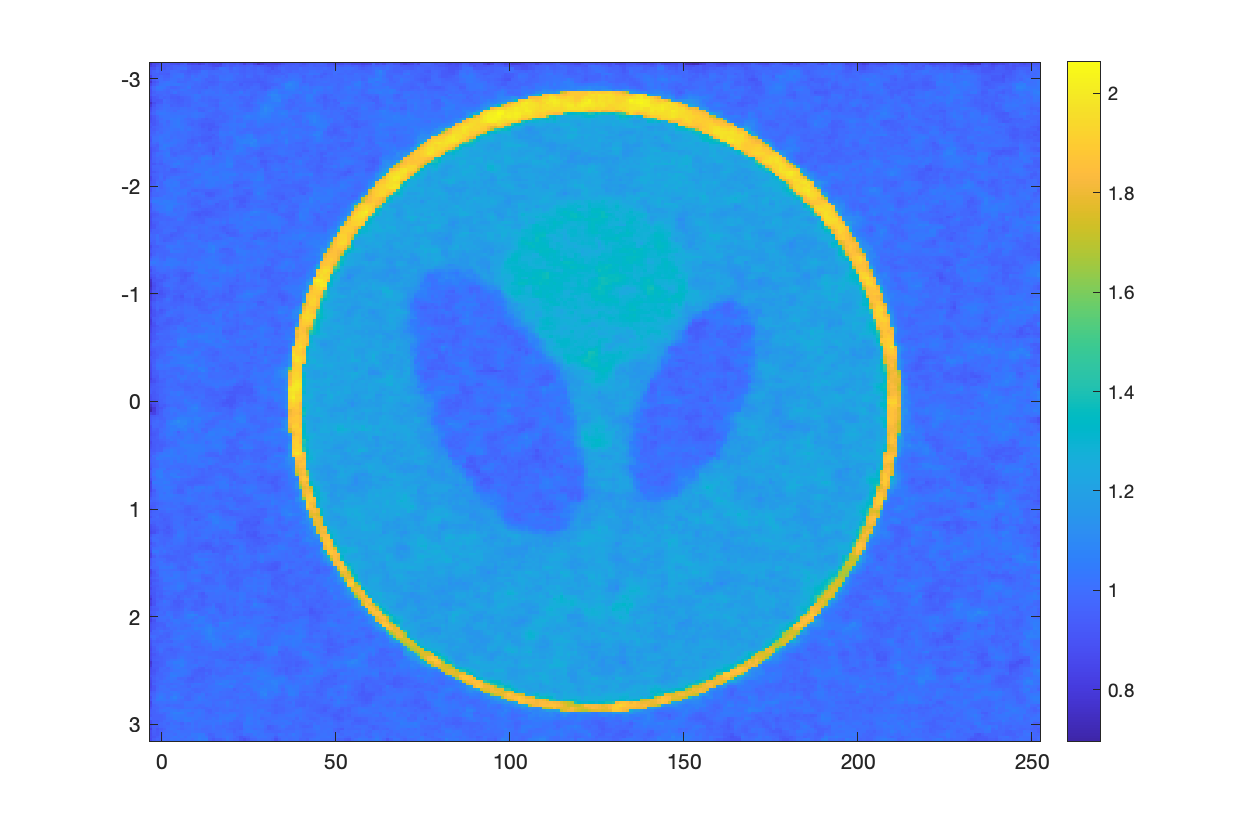}
    \caption{Posterior mean}
    \end{subcaptionblock}
    \begin{subcaptionblock}{.23\textwidth}
    \centering
    \includegraphics[width=1\linewidth]{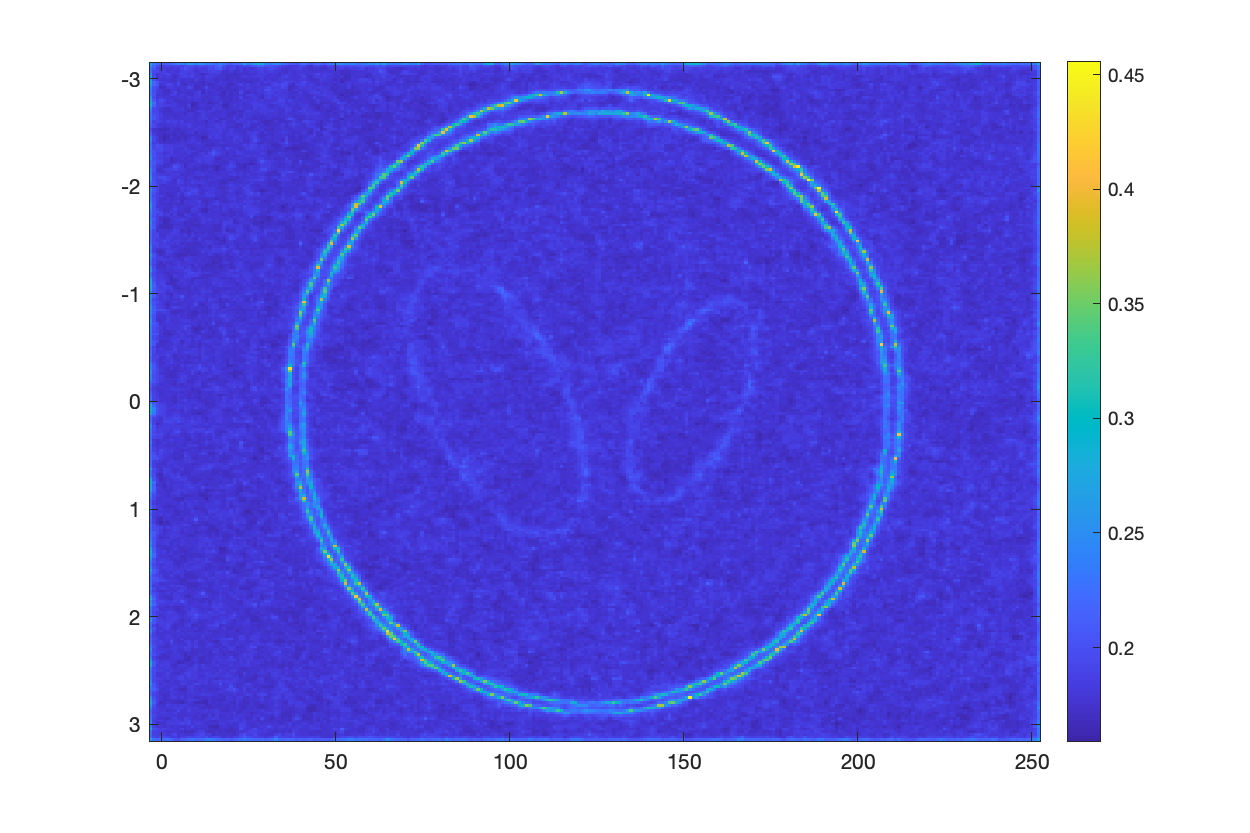}
    \caption{$90\%$ CI}
    \end{subcaptionblock}\\
    \begin{subcaptionblock}{.23\textwidth}
    \centering
    \includegraphics[width=1\linewidth]{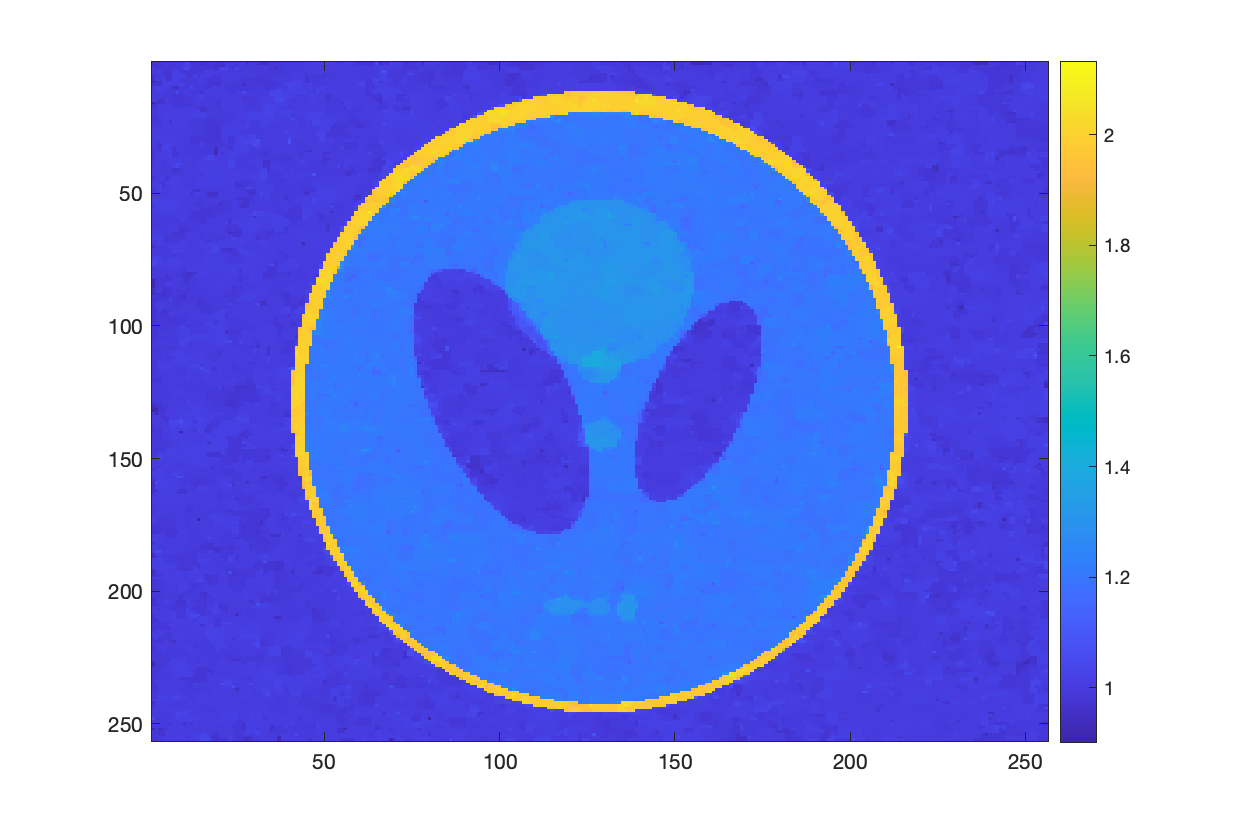}
    \caption{Generalized LASSO}
    \end{subcaptionblock}
    \begin{subcaptionblock}{.23\textwidth}
    \centering
    \includegraphics[width=1\linewidth]{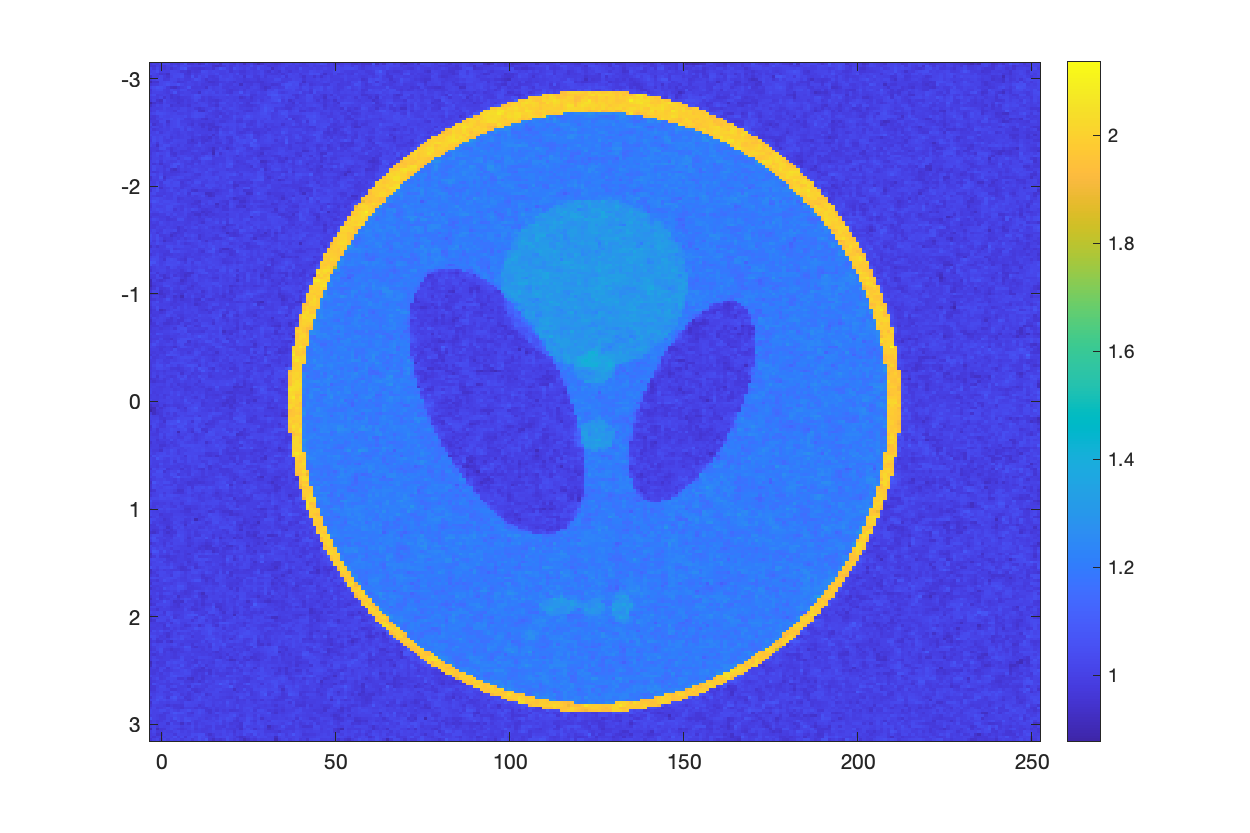}
    \caption{Posterior mean}
    \end{subcaptionblock}
    \begin{subcaptionblock}{.23\textwidth}
    \centering
    \includegraphics[width=1\linewidth]{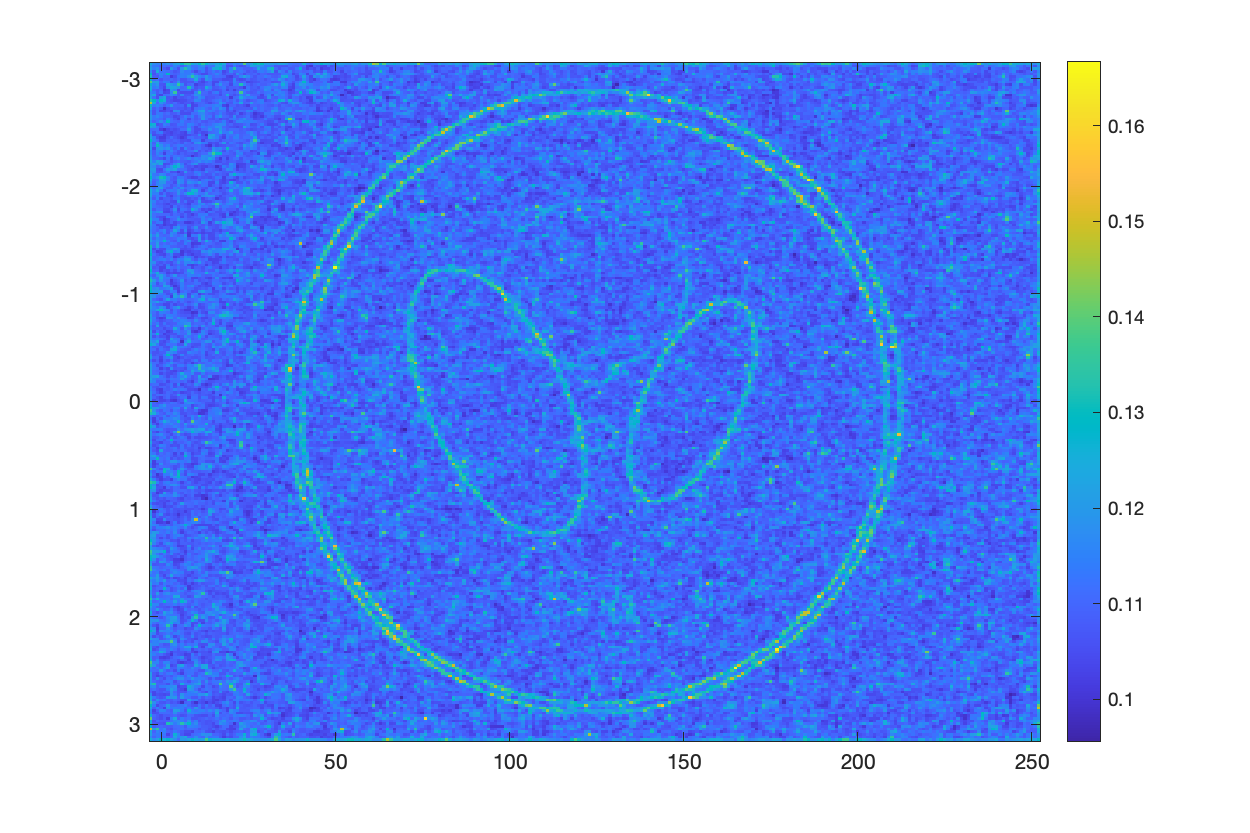}
    \caption{$90\%$ CI}
    \end{subcaptionblock}
     \caption{Magnitude recovery of the complex-valued Shepp Logan phantom  using $F_B$ in \eqref{eq:blurtransform}. (a)-(c) and (d)-(f) correspond to SNR values of $15\dB$ and $25\dB$, respectively. $\hat{\eta}=10^{-2}$.}
    \label{fig:2D_transform_FB}
\end{figure}

We  now consider the case of promoting sparsity in the  magnitude gradient for the Shepp-Logan phantom image \cite{shepp1974fourier} depicted  in \cref{fig:2D_transform_FF} (a)  on a $256\cross256$ grid. A random phase is then added to each pixel, so that the resulting image  $\bm z$ is modeled by  \eqref{eqn:observation_model}. \Cref{fig:2D_transform_FF} and \cref{fig:2D_transform_FB} compare the posterior means of the  CVBL to the generalized LASSO estimates for $F_F$ in \cref{eq:Fouriertransform} and $F_B$ in \eqref{eq:blurtransform}.\footnote{The dense nature of $F_U$ in \eqref{eq:undertransform} makes exact sampling of the phase $\bm\phi$ impractical for large signals. Approximate or stochastic techniques may expedite sampling, but is beyond the scope of the current investigation.} To avoid the issue discussed in \cref{rem:flat_post}, in the Shepp Logan experiments we fix $\eta$ to be $\hat{\eta}=10^{-2}$, which was chosen heuristically, with $\lambda = \sigma^2/ \hat{\eta}$ used for the generalized LASSO regularization parameter. The hyperprior in turn is given by \eqref{eq:delta_hyper} rather than \eqref{eq:hyperprior}.

\cref{fig:2D_transform_FF} and \cref{fig:2D_transform_FB} indicate that the generalized LASSO estimate  \eqref{eq:sparseedge_classical_lasso} is smoother when compared to the CVBL posterior mean. In particular we observe that the three lower circles of the generalized LASSO estimate of the phantom are ``blurred'' together but remain  distinct for the CVBL. These results seem to indicate that when using the CVBL method, a smaller value of $\hat{\eta}$, that is, increasing the reliance on the prior, is needed to match the amount of regularization apparent the generalized LASSO approach.

\begin{figure}[h!]
    \centering
    \begin{subcaptionblock}{.21\textwidth}
    \centering
    \includegraphics[width=1\linewidth]{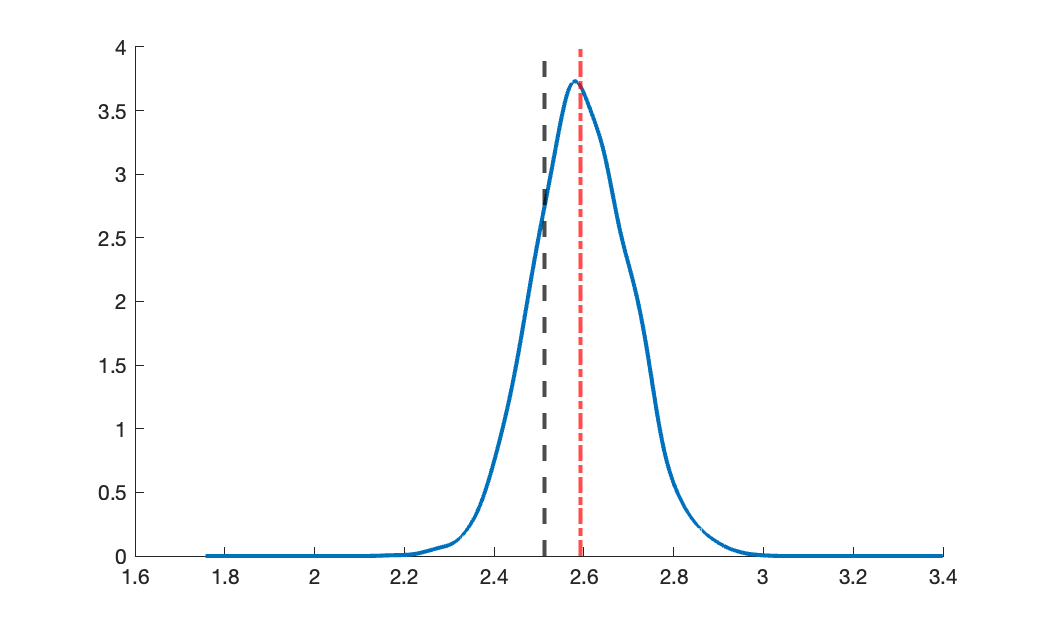}
    \caption{$F_F$, $15\dB$}
    \end{subcaptionblock}
    \begin{subcaptionblock}{.21\textwidth}
    \centering
    \includegraphics[width=1\linewidth]{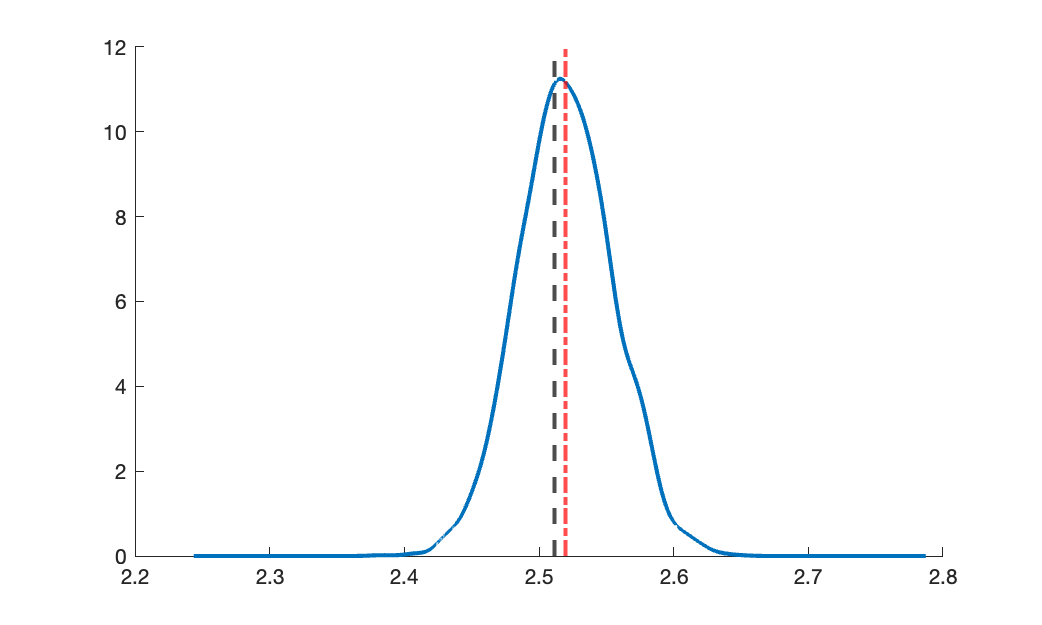}
    \caption{$F_F$, $25\dB$}
    \end{subcaptionblock}
    \begin{subcaptionblock}{.21\textwidth}
    \centering
    \includegraphics[width=1\linewidth]{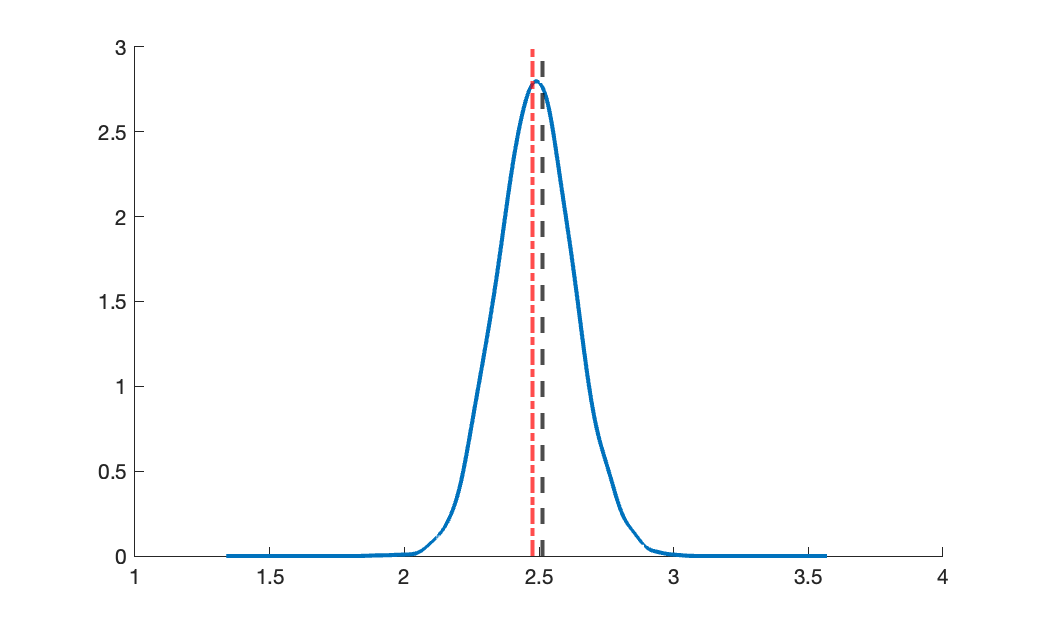}
    \caption{$F_B$, $15\dB$}
    \end{subcaptionblock}
    \begin{subcaptionblock}{.21\textwidth}
    \centering
    \includegraphics[width=1\linewidth]{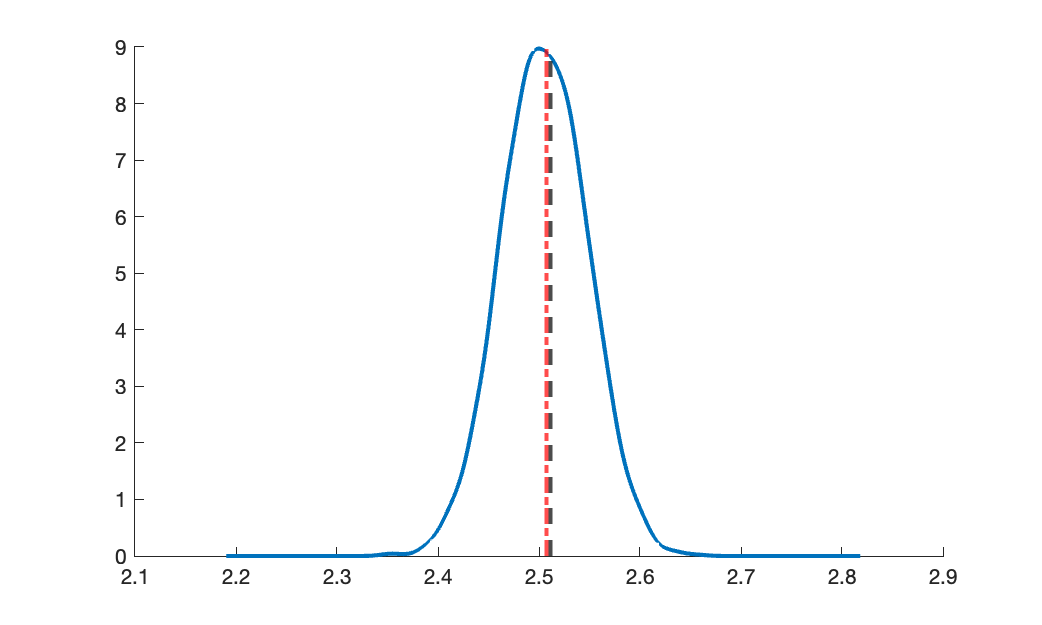}
    \caption{$F_B$, $25\dB$}
    \end{subcaptionblock}
    \begin{subcaptionblock}[b][2.4cm][t]{.12\textwidth}
    \centering
    \includegraphics[width=1\linewidth]{legendPhase.png}
    \end{subcaptionblock}
    \caption{Recovered phase at a randomly chosen pixel of the complex-valued Shepp Logan image using $F_F$ in \eqref{eq:Fouriertransform} and $F_B$ in \eqref{eq:blurtransform} with SNR of $15\dB$ and $25\dB$.  $\hat{\eta}=10^{-2}$.}
    \label{fig:2D_transform_phase}
\end{figure}

Finally, \cref{fig:2D_transform_phase} demonstrates that the true phase of a randomly selected pixel is within an area of high density of the kernel density estimation of the marginal posterior density of the phase, and that the kernel density estimation mode is consistent with those obtained using the generalized LASSO technique.
\section{Concluding remarks and future work}\label{sec:conclusion} 
This investigation extends the real-valued Bayesian LASSO  (RVBL), which was originally designed to promote sparsity in a sparse signal, in two ways.  We first show that it can be modified to promote sparsity in a chosen transform (here the gradient) domain.   We then demonstrate that the RVBL  can be further expanded to include {\em complex}-valued signals and images.  We call our method the complex-valued Bayesian LASSO (CVBL).  Our numerical experiments show that the CVBL can efficiently recover samples from the entire complex-valued posterior density function, enabling uncertainty quantification of both the magnitude and  phase of the true signal.  

The CVBL is practical for coherent imaging problems
with unitary or sparse forward operators, since it is easily parallelizable.  Developing surrogates for dense forward operators will be necessary to efficiently sample large problems, and will be the focus of future investigations. We will also consider different sparse transform operators (such at HOTV) along with adaptive empirical hyperparameters to avoid the pitfulls of over-regularization. Recent work in \cite{zhang2022empirical} may be useful in this regard.

\appendix
\section{Proofs}\label{sec:appendix_proofs}

\begin{manualtheorem}{4.1}
Let  $F_1 = FD(e^{i\phi})$, $\tilde{F}_1=[\Re(F_1)^T \ \Im(F_1)^T]^T$, and $\tilde{\bm y}=[\Re(\bm y)^T \ \Im(\bm y)^T]^T$. Assume that $L$ is rank $n$. The function $\tilde{f}_{\mathcal{G}|\mathcal{Y},\Phi,\mathcal{T}^2}(\bm g|\bm y,\bm \phi,\bm\tau^2)$ in \eqref{eq:magpost} defines a Gaussian density over $\mathbb{R}^n$ with center and precision given by
\begin{subequations}
\label{eq:centerprecision_app}
    \begin{equation}
	\bar{\bm g} = \Gamma^{-1}\left(\frac{2}{\sigma^2}\tilde{F}_1^T\tilde{\bm y}\right),
 \end{equation}
 \begin{equation}
	\Gamma = L^T[D(\bm\tau^2)]^{-1}L+\frac{2}{\sigma^2}\tilde{F}_1^T\tilde{F}_1.
 \end{equation}
\end{subequations}
\end{manualtheorem}

\begin{proof}
Let $\Gamma=L^T[D(\bm\tau^2)]^{-1}L+\frac{2}{\sigma^2}\tilde{F}_1^T\tilde{F}_1$, and suppose $\bm x\in\Ker(\Gamma)$ so that $\bm x^T\Gamma \bm x = 0$. By \eqref{eq:Gamma} we then have
\begin{align*}
		\bm x^T\Gamma \bm x=\bm x^TL^T[D(\bm\tau^2)]^{-1}L\bm x+\frac{2}{\sigma^2}\bm x^T\tilde{F}_1^T\tilde{F}_1\bm x=0,
	\end{align*}
that is $\bm x\in\Ker(L)\cap\Ker(\tilde{F}_1)$. Since $\Ker(L)=\{0\}$, we have immediately that $\Ker(\Gamma)=\{0\}$ so that $\Gamma$ is invertible. We furthermore observe that $\Gamma$ is also symmetric positive definite. To explicitly determine the center and precision as ${\bar{\bm g}}$ and $\Gamma$ in \eqref{eq:barg} and \eqref{eq:Gamma}, consider the following from \eqref{eq:magpost}:
	\begin{align*}
		-\frac{\norm{\bm y-F_1\bm g}^2_2}{\sigma^2}-\frac{1}{2}(L\bm g)^T[D(\bm\tau^2)]^{-1}(L\bm g)&=-\frac{\norm{\tilde{\bm y}-\tilde{F}_1\bm g}^2_2}{\sigma^2}-\frac{1}{2}(L\bm g)^T[D(\bm\tau^2)]^{-1}(L\bm g)\\
  &=-\frac12(L\bm g)^T[D(\bm\tau^2)]^{-1}(L\bm g)-\frac{1}{\sigma^2}(\tilde{\bm y}-\tilde{F}_1\bm g)^T(\tilde{\bm y}-\tilde{F}_1\bm g)\\
		&=-\frac12\bm g^TL^T[D(\bm\tau^2)]^{-1}L\bm g-\frac{1}{\sigma^2}\tilde{\bm y}^T\tilde{\bm y}+\frac{1}{\sigma^2}\bm g^T\tilde{F}_1^T\tilde{\bm y}\\
  &\quad+\frac{1}{\sigma^2}\tilde{\bm y}^T\tilde{F}_1\bm g-\frac{1}{\sigma^2}\bm g^T\tilde{F}_1^T\tilde{F}_1\bm g\\
		&=-\frac12\bm g^T\Gamma \bm g+\frac{1}{\sigma^2}\bm g^T\tilde{F}_1^T\tilde{\bm y}+\frac{1}{\sigma^2}\tilde{\bm y}^T\tilde{F}_1\bm g-\frac{1}{\sigma^2}\tilde{\bm y}^T\tilde{\bm y}\\
            &=-\frac12\bm g^T\Gamma \bm g+\frac{2}{\sigma^2}\bm g^T\tilde{F}_1^T\tilde{\bm y}-\frac{1}{\sigma^2}\tilde{\bm y}^T\tilde{\bm y}\\
            &=-\frac12\bm g^T\Gamma \bm g+\bm g^T\Gamma\bar{\bm g}-\frac12\bar{\bm g}^T\Gamma\bar{\bm g} +\frac12\bar{\bm g}^T\Gamma\bar{\bm g}-\frac{1}{\sigma^2}\tilde{\bm y}^T\tilde{\bm y}\\
		&=-\frac12(\bm g-\bar{\bm g})^T\Gamma(\bm g-\bar{\bm g})+c(\tilde{\bm y}),
	\end{align*}
where $c(\tilde{\bm y})=\frac12\bar{\bm g}^T\Gamma\bar{\bm g}-\frac{1}{\sigma^2}\tilde{\bm y}^T\tilde{\bm y}$ is mutually independent of $\bm g$. Thus we have mean $\bar{\bm g}$ and precision $\Gamma$, yielding the desired result.
\end{proof}

\begin{manualcorollary}{4.2}
Let $F_1 = FD(e^{i\phi})$ and suppose $F$ is unitary and $L$ is rank $n$. For $F_1^H\bm y=\bm q\odot e^{i\bm\varphi}$ where $\bm q\in\left(\mathbb{R}^+\right)^n$ and $\bm\varphi\in[-\pi,\pi)^n$, the function $\tilde{f}_{\mathcal{G}|\mathcal{Y},\mathcal{T}^2,\Phi}(\bm g|\bm y,\bm\tau^2,\bm \phi)$ in \eqref{eq:magpost} defines a Gaussian density over $\mathbb{R}^n$ with center and precision given by
\begin{subequations}
\label{eq:precisioncenter_unity_app}
    \begin{equation}
	\bar{\bm g} = \Gamma^{-1}\left(\frac{2}{\sigma^2}\bm q\odot\cos(\bm\varphi)\right),
 \end{equation}
 \begin{equation}
	\Gamma = L^T[D(\bm\tau^2)]^{-1}L+\frac{2}{\sigma^2}\mathbb{I}_n.
 \end{equation}
\end{subequations}
\end{manualcorollary}

\begin{proof}
    Since $F$ is unitary, we have in \eqref{eq:magpost}
    \begin{align*}
        -\frac{\norm{\bm y-F_1\bm g}^2_2}{\sigma^2}-\frac{1}{2}(L\bm g)^T[D(\bm\tau^2)]^{-1}(L\bm g)&=-\frac12\bm g^T\Gamma \bm g+\frac{1}{\sigma^2}\bm g^T\left(\bm q\odot e^{i\bm \varphi}\right)\\&\quad+\frac{1}{\sigma^2}\left(\bm q\odot e^{-i\bm \varphi}\right)^T\bm g-\frac{1}{\sigma^2}\bm y^H\bm y\\
		&=-\frac12\bm g^T\Gamma \bm g+\frac{1}{\sigma^2}\sum_{j=1}^ng_jq_j\left(e^{i\varphi_j}+e^{-i\varphi_j}\right)-\frac{1}{\sigma^2}\bm y^H\bm y\\
		&=-\frac12\bm g^T\Gamma \bm g+\frac{2}{\sigma^2}\sum_{j=1}^ng_jq_j\cos(\varphi_j)-\frac{1}{\sigma^2}\bm y^H\bm y\\
		&=-\frac12\bm g^T\Gamma \bm g+\bm g^T\Gamma\bar{\bm g}-\frac12\bar{\bm g}^T\Gamma\bar{\bm g} + \frac12\bar{\bm g}^T\Gamma\bar{\bm g}-\frac{1}{\sigma^2}\bm y^H\bm y\\
		&=-\frac12(\bm g-\bar{\bm g})^T\Gamma(\bm g-\bar{\bm g})+c(\bm y)
    \end{align*}
    where $c(\bm y)=\frac12\bar{\bm g}^T\Gamma\bar{\bm g}-\frac{1}{\sigma^2}\bm y^H\bm y$ is mutually independent of $\bm g$. Thus we have mean $\bar{\bm g}$ and precision $\Gamma$, yielding the desired result.
\end{proof}

\begin{manualtheorem}{4.3}
    Let $F_2=FD(\bm g)$ as in \eqref{eq:edge_postphi} and \eqref{eq:edge_posttheta} and define $A=F_2^HF_2$ with elements $[A]_{j,k}=a_{j,k}e^{i\alpha_{j,k}}$, where $a_{j,k}\in\mathbb{R}$ and $\alpha_{j,k}\in[-\pi,\pi)$. Further denote $F_2^H\bm y=\bm q\odot e^{i\bm \varphi}$, where $\bm q\in\mathbb{R}_+^{n}$ and $\bm \varphi\in[-\pi,\pi)^{n}$, with
    \begin{gather*}
        u_i=\frac{2q_i}{\sigma^2}\cos\varphi_i-\sum\limits_{\substack{k=1 \\ k\neq i}}^n\frac{2a_{i,k}}{\sigma^2}\cos(\tilde{\alpha}_{i,k}+\phi_k),\quad\quad
        v_i=\frac{2q_i}{\sigma^2}\sin\varphi_i-\sum\limits_{\substack{k=1 \\ k\neq i}}^n\frac{2a_{i,k}}{\sigma^2}\sin(\tilde{\alpha}_{i,k}+\phi_k),
    \end{gather*}
    where $\tilde{\alpha}_{i,k}=\sgn(k-i)\alpha_{i,k}$ for $k=1,\dots,n$. Then
\begin{align*}
    f_{\Phi_i|\Phi_{-i},\mathcal{Y},\mathcal{G},\mathcal{T}^2}(\phi_i|\phi_{-i},\bm y,\bm g,\bm\tau^2)\propto\pi_{vM}\left(\phi_{i}\Big|\mu_i,\kappa_i\right),\quad i=1,\dots,n,
\end{align*}
where
\begin{gather*}                 \kappa_i=\sqrt{u_i^2+v_i^2},\quad\quad\mu_i=\begin{cases}\arctan\left(-\frac{v_i}{u_i}\right) & \text{if } u_i>0 \\ \pi/2 & \text{if } u_i=0 \\ \arctan\left(-\frac{v_i}{u_i}\right)+\pi & \text{if } u_i<0 
    \end{cases}
\end{gather*}
and $\pi_{vM}(x|\mu,\kappa)$ is the von Mises probability density function, given as
\begin{align*}
    \pi_{vM}(x|\mu,\kappa)=\frac{\exp(\kappa\cos(x-\mu))}{2\pi I_0(\kappa)},
\end{align*}
with location $\mu$ and concentration $\kappa$, and where $I_0$ is the modified Bessel function of the first kind of order zero.
\end{manualtheorem}

\begin{proof}
   By direct calculation using \eqref{eq:edge_postphi}, we have
\begin{align*}
	f_{\Phi_i|\Phi_{-i},\mathcal{Y},\mathcal{G},\mathcal{T}^2}&(\phi_i|\phi_{-i},\bm y,\bm g,\bm\tau^2)\propto\exp(-\frac{1}{\sigma^2}\norm{\bm y-F_2e^{i\bm \phi}}^2_2)\\
	&\propto\exp(-\frac{1}{\sigma^2}\left(e^{-i\bm\phi^T}Ae^{i\bm\phi}-\bm y^HF_2e^{i\bm\phi} - e^{-i\bm\phi^T}F_2^H\bm y\right))\\
        &=\exp\Bigg(-\frac{1}{\sigma^2}\bigg(\sum\limits_{\substack{j,k=1 \\ k>j}}^n \left(a_{j,k}e^{i(\alpha_{j,k}-\phi_j+\phi_k)}+a_{j,k}e^{-i(\alpha_{j,k}-\phi_j+\phi_k)}\right)-2\bm q^T\cos(\bm\phi-\bm\varphi)\bigg)\Bigg)\\
        &=\exp\Bigg(-\frac{1}{\sigma^2}\bigg(\sum\limits_{\substack{j,k=1 \\ k>j}}^n 2a_{j,k}\cos(\alpha_{j,k}-\phi_j+\phi_k)-2\bm q^T\cos(\bm\phi-\bm\varphi)\bigg)\Bigg)\\
        &\propto\exp\Bigg(-\frac{1}{\sigma^2}\bigg(\sum\limits_{\substack{k=1 \\ k\neq i}}^n 2a_{i,k}\left(\cos\phi_i\cos(\phi_k+\tilde{\alpha}_{i,k})+\sin\phi_i\sin(\phi_k+\tilde{\alpha}_{i,k})\right)-2 q_i\cos(\phi_i-\varphi_i)\bigg)\Bigg)\\
        &=\exp\Bigg(\bigg(\frac{2q_i}{\sigma^2}\cos\varphi_i-\sum\limits_{\substack{k=1 \\ k\neq i}}^n \frac{2a_{i,k}}{\sigma^2}\cos(\phi_k+\tilde{\alpha}_{j,k})\bigg)\cos\phi_i\\
        &\quad\quad\quad\quad+\bigg(\frac{2q_i}{\sigma^2}\sin\varphi_i-\sum\limits_{\substack{k=1 \\ k\neq i}}^n \frac{2a_{i,k}}{\sigma^2}\sin(\phi_k+\tilde{\alpha}_{i,k})\bigg)\sin\phi_i\Bigg)\\
        &=\exp(u_i\cos\phi_i+v_i\sin\phi_i)\\
        &=\exp(\kappa_i\cos(\phi_i-\mu_i)).
\end{align*}
\end{proof}

\begin{manualtheorem}{4.4}
Suppose $F_2=FD(\bm g)$, where $F$ is unitary. Let $F_2^H\bm y=\bm q\odot e^{i\bm \varphi}$ where $\bm q\in\mathbb{R}_+^{n}$ and $\bm \varphi\in[-\pi,\pi)^{n}$. Then 
\begin{align*}
    f_{\Phi|\mathcal{Y},\mathcal{G}}(\bm \phi|\bm y,\bm g)=\prod_{i=1}^n\pi_{vM}\left(\phi_{i}\Big|\varphi_{i},\frac{2}{\sigma^2}q_{i}\right).
\end{align*}
\end{manualtheorem}

\begin{proof}
Consider the following:
\begin{align*}
	f_{\Phi|\mathcal{Y},\mathcal{G}}(\bm \phi|\bm y,\bm g)&=C\exp(-\frac{1}{\sigma^2}\norm{\bm y-F_2e^{i\bm \phi}}^2_2)\bm1_{[-\pi,\pi)}(\bm \phi)\\
	&=C\exp(-\frac{1}{\sigma^2}\left(e^{-i\bm \phi}F_2^HF_2e^{i\bm \phi}-\bm y^HF_2e^{i\bm \phi}-e^{-i\bm \phi}F_2^H\bm y+\bm y^H\bm y\right))\bm1_{[-\pi,\pi)}(\bm \phi)\\
	&=C\exp(-\frac{1}{\sigma^2}\left(\bm g^T\bm g-2\bm q^T\cos(\bm \phi-\bm\varphi)+\bm y^H\bm y\right))\bm1_{[-\pi,\pi)}(\bm \phi),
\end{align*}
where $C$ is a normalization constant. To find $C$, we integrate the above expression over $\mathbb{R}^n$ as follows:
\begin{align*}
	\frac1C&=\int_{\mathbb{R}^n}\exp(-\frac{1}{\sigma^2}\left(\bm g^T\bm g-2\bm q^T\cos(\bm \phi-\bm\varphi)+\bm y^H\bm y\right))\bm1_{[-\pi,\pi)}(\bm \phi)\mathrm{d}\bm \phi\\
	&=\int_{[-\pi,\pi)^n}\exp(-\frac{1}{\sigma^2}\left(\bm g^T\bm g-2\bm q^T\cos(\bm \phi-\bm\varphi)+\bm y^H\bm y\right))\mathrm{d}\bm \phi\\
	&=\exp(-\frac{1}{\sigma^2}\left(\bm g^T\bm g+\bm y^H\bm y\right))\int_{[-\pi,\pi)^n}\exp(\frac{2}{\sigma^2}\bm q^T\cos(\bm \phi-\bm\varphi))\mathrm{d}\bm \phi\\
	&=\exp(-\frac{1}{\sigma^2}\left(\bm g^T\bm g+\bm y^H\bm y\right))\prod_{i=1}^n2\pi I_0\left(\frac{2}{\sigma^2}q_{i}\right).
\end{align*}
Thus, we have
\begin{align*}
	f_{\Phi|\mathcal{Y},\mathcal{G},\mathcal{T}^2}(\bm \phi|\bm y,\bm g,\bm\tau^2)&=\frac{\exp(-\frac{1}{\sigma^2}\left(\bm g^T\bm g-2\bm q^T\cos(\bm \phi-\bm\varphi)+\bm y^H\bm y\right))\bm1_{[-\pi,\pi)}(\bm \phi)}{\exp(-\frac{1}{\sigma^2}\left(\bm g^T\bm g+\bm y^H\bm y\right))\prod_{i=1}^n2\pi I_0\left(\frac{2}{\sigma^2}q_{i}\right)}\\
	&=\prod_{j=1}^n\frac{\exp(\frac{2}{\sigma^2}q_{j}\cos(\phi_{j}-\varphi_{j}))}{2\pi I_0\left(\frac{2}{\sigma^2}q_{j}\right)}\bm1_{[-\pi,\pi)}(\phi_{j})\\
	&=\prod_{j=1}^n\pi_{vM}\left(\phi_{j}\Big|\varphi_{j},\frac{2}{\sigma^2}q_{j}\right),
\end{align*}
which completes the proof.
\end{proof}




\bibliographystyle{siamplain}
\bibliography{literature}

\end{document}